\documentclass{amsart} 

\usepackage{amsmath}
\usepackage{amsthm}
\usepackage{amssymb}
\usepackage{amsfonts}
\usepackage{amscd}

\begin{document}

\newcommand{\ba}{{\bold a}}
\newcommand{\balpha}{{\boldsymbol \alpha}}
\newcommand{\blambda}{{\boldsymbol \lambda}}
\newcommand{\bmu}{{\boldsymbol \mu}}
\newcommand{\bnu}{{\boldsymbol \nu}}
\newcommand{\bbeta}{{\boldsymbol \beta}}
\newcommand{\bt}{{\bold t}}
\newcommand{\bH}{{\mathbf H}}
\newcommand{\Z}{{\bold Z}}
\newcommand{\R}{{\bold R}}
\newcommand{\C}{{\bold C}}
\newcommand{\BP}{{\mathbf P}}
\newcommand{\cA}{{\mathcal A}}
\newcommand{\cG}{{\mathcal G}}
\newcommand{\cO}{{\mathcal O}}
\newcommand{\cM}{{\mathcal M}}
\newcommand{\cC}{{\mathcal C}}
\newcommand{\cF}{{\mathcal F}}
\newcommand{\cE}{{\mathcal E}}
\newcommand{\V}{{\mathbf V}}
\newcommand{\cL}{{\mathcal L}}

\newcommand\rank{\mathop{\rm rank}\nolimits}
\newcommand\im{\mathop{\rm Im}\nolimits}
\newcommand\coker{\mathop{\rm coker}\nolimits}
\newcommand{\Aut}{\mathop{\rm Aut}\nolimits}
\newcommand{\Tr}{\mathop{\rm Tr}\nolimits}
\newcommand\Hom{\mathop{\rm Hom}\nolimits}
\newcommand\Ext{\mathop{\rm Ext}\nolimits}
\newcommand\End{\mathop{\rm End}\nolimits}
\newcommand\Pic{\mathop{\rm Pic}\nolimits}
\newcommand\Spec{\mathop{\rm Spec}\nolimits}
\newcommand\Hilb{\mathop{\rm Hilb}\nolimits}
\newcommand\RH{\mathop{\bf RH}\nolimits}
\newcommand{\length}{\mathop{\rm length}\nolimits}
\newcommand{\res}{\mathop{\sf res}\nolimits}
\newcommand\Quot{\mathop{\rm Quot}\nolimits}
\newcommand\Grass{\mathop{\rm Grass}\nolimits}
\newcommand\Proj{\mathop{\rm Proj}\nolimits}
\newcommand{\codim}{\mathop{\rm codim}\nolimits}

\newcommand\lra{\longrightarrow}
\newcommand\ra{\rightarrow}
\newcommand\la{\leftarrow}
\newtheorem{Theorem}{Theorem}[section]
\newtheorem{Lemma}{Lemma}[section]
\newtheorem{Remark}{Remark}[section]
\newtheorem{Corollary}{Corollary}[section]
\newtheorem{Proposition}{Proposition}[section]
\newtheorem{Example}{Example}[section]
\newtheorem{Definition}{Definition}[section]

\title{Moduli of parabolic connections on a curve
and Riemann-Hilbert correspondence}
\author{Michi-aki Inaba}

\address{{\rm Michi-aki Inaba} \\
Department of Mathematics, Kyoto University, \\
Kyoto, 606-8502, Japan}
\email{inaba@math.kyoto-u.ac.jp}

\begin{abstract}
Let $(C,\bt)$ ($\bt=(t_1,\ldots,t_n)$)
be an $n$-pointed smooth projective curve of genus $g$
and take an element $\blambda=(\lambda^{(i)}_j)\in\C^{nr}$
such that $-\sum_{i,j}\lambda^{(i)}_j=d\in\mathbf{Z}$.
For a weight $\balpha$, let $M_C^{\balpha}(\bt,\blambda)$
be the moduli space of $\balpha$-stable $(\bt,\blambda)$-parabolic
connections on $C$ and
let $RP_r(C,\bt)_{\ba}$ be the moduli space of representations
of the fundamental group $\pi_1(C\setminus\{t_1,\ldots,t_n\},*)$
with the local monodromy data $\ba$ for a certain $\ba\in\C^{nr}$.
Then we prove that the morphism
$\RH:M_C^{\balpha}(\bt,\blambda)\rightarrow RP_r(C,\bt)_{\ba}$
determined by the Riemann-Hilbert correspondence
is a proper surjective bimeromorphic morphism.
As a corollary, we prove the geometric Painlev\'e property of
the isomonodromic deformation defined on the moduli space of
parabolic connections.
\end{abstract}

\maketitle

\section{Introduction}

Let $C$ be a smooth projective curve over $\C$
and $t_1,\ldots,t_n$ be distinct points of $C$.
For an algebraic vector bundle $E$ on $C$
and a logarithmic connection
$\nabla: E \rightarrow E\otimes\Omega^1_C(t_1+\cdots+t_n)$,
$\ker\nabla^{an}|_{C\setminus\{t_1,\ldots,t_n\}}$ becomes
a local system on $C\setminus\{t_1,\ldots,t_n\}$
and corresponds to a representation of the fundamental group
$\pi_1(C\setminus\{t_1,\ldots,t_n\},*)$,
where $\nabla^{an}$ is the analytic connection corresponding to $\nabla$.
The correspondence 
$(E,\nabla)\mapsto
\ker\nabla^{an}|_{C\setminus\{t_1,\ldots,t_n\}}$
is said to be the Riemann-Hilbert correspondence.
If $l\subset E|_{t_i}$ is a subspace satisfying
$\res_{t_i}(\nabla)(l)\subset l$,
then $\nabla$ induces a connection
$\nabla':E'\ra E'\otimes\Omega^1_C(t_1+\cdots+t_n)$,
where $E':=\ker(E\ra (E|_{t_i}/l))$.
We say $(E',\nabla')$ the elementary transform of
$(E,\nabla)$ along $t_i$ by $l$.
Note that
$\ker\nabla^{an}|_{C\setminus\{t_1,\ldots,t_n\}}\cong
\ker(\nabla')^{an}|_{C\setminus\{t_1,\ldots,t_n\}}$.

We put
\[
 \Lambda^{(n)}_r(d):=
 \left\{ (\lambda^{(i)}_j)^{1\leq i\leq n}_{0\leq j\leq r-1}
 \in \C^{nr}
 \left|
 d+\sum_{i,j}\lambda^{(i)}_j=0
 \right\}\right.
\]
for integers $d,r,n$ with $r>0$, $n>0$.
We write $\bt=(t_1,\ldots,t_n)$ and take
an element $\blambda\in\Lambda^{(n)}_r(d)$.
\begin{Definition}\rm
We say $(E,\nabla,\{l^{(i)}_*\}_{1\leq i\leq n})$
a $(\bt,\blambda)$-parabolic connection of rank $r$ if
\begin{enumerate}
\item $E$ is a rank $r$ algebraic vector bundle on $C$,
\item $\nabla: E \ra E\otimes\Omega_C^1(t_1+\cdots+t_n)$
 is a connection,
and
\item for each $t_i$,
 $l^{(i)}_*$ is a filtration
 $E|_{t_i}=l^{(i)}_0\supset l^{(i)}_1
 \supset\cdots\supset l^{(i)}_{r-1}\supset l^{(i)}_r=0$
 such that $\dim(l^{(i)}_j/l^{(i)}_{j+1})=1$ and
 $(\res_{t_i}(\nabla)-\lambda^{(i)}_j\mathrm{id}_{E|_{t_i}})
 (l^{(i)}_j)\subset l^{(i)}_{j+1}$
 for $j=0,\ldots,r-1$.
\end{enumerate}
\end{Definition} 
The filtration $\l^{(i)}_*$ ($1\leq i\leq n$) is said to be
a parabolic structure of the vector bundle $E$.
For a parabolic connection $(E,\nabla,\{l^{(i)}_j\})$,
the elementary transform of $(E,\nabla)$ along $t_i$ by
$l^{(i)}_j$ determines another parabolic connection.
Then an elementary transform gives a transformation
$\mathrm{Elm}^{(i)}_j$ on the
set of isomorphism classes of parabolic connections.
See section \ref{elementary-transform-m}, (\ref{elementary-transform})
for the precise definition of
$\mathrm{Elm}^{(i)}_j$.
Then we should say that the Riemann-Hilbert correspondence gives
a bijection $\mathrm{RH}$ between the set
\[
 (\dag)\quad
 \left.\left\{ (E,\nabla,\{l^{(i)}_j\}):
 \text{parabolic connection of rank $r$}
 \right\}\right/\sim
\]
and the set
\[
 (\dag\dag)\quad
 \left\{
 \rho:\pi_1(C\setminus\{t_1,\ldots,t_n\},*) \rightarrow GL_r(\C):
 \text{representation}\right\}/\cong,
\]
where $\sim$ is the equivalence relation generated by
$x\sim \mathrm{Elm}_{l^{(i)}_j}(x)$,
$(E,\nabla,\{l^{(i)}_j\})\sim(E,\nabla,\{l^{(i)}_j\})\otimes{\mathcal O}_C(t_k)$
and $(E,\nabla,\{l^{(i)}_j\})\sim(E,\nabla,\{l'^{(i)}_j\})$.
Then the bijectivity of $\mathrm{RH}$ immediately follows from the theory
of Deligne (\cite{Deligne}).
Indeed take any representation
$\rho:\pi_1(C\setminus\{t_1,\ldots,t_n\},*)\rightarrow GL_r(\mathbf{C})$.
Then $\rho$ corresponds to a locally constant sheaf
$\mathbf{V}$ on $C\setminus\{t_1,\ldots,t_n\}$.
By [\cite{Deligne}, II, Proposition 5.4], there is a unique
pair $(E,\nabla)$ of  a vector bundle $E$ and
a connection
$\nabla:E\rightarrow E\otimes\Omega^1_C(t_1+\cdots+t_n)$
such that all the eigenvalues of $\res_{t_i}(\nabla)$ lie in
$\{\lambda\in\mathbf{C}|0\leq\mathrm{Re}(\lambda)<1\}$
and $\ker\nabla^{an}|_{C\setminus\{t_1,\ldots,t_n\}}\cong\mathbf{V}$.
We can take a parabolic structure $\{l^{(i)}_j\}$ on
$E$ compatible with connection $\nabla$.
Then $(E,\nabla,\{l^{(i)}_j\})$ becomes a parabolic connection
and $\mathrm{RH}([(E,\nabla,\{l^{(i)}_j\})])=\rho$.
So $\mathrm{RH}$ is surjective.
Let $(E,\nabla,\{l^{(i)}_j\})$ be a $(\bt,\blambda)$-parabolic connection
and $(E',\nabla',\{(l')^{(i)}_j\})$ be a $(\bt,\blambda')$-parabolic connection
such that
$\mathrm{RH}([(E,\nabla,\{l^{(i)}_j)])=\mathrm{RH}([(E',\nabla',\{(l')^{(i)}_j\})])$,
namely
$\ker\nabla^{an}|_{C\setminus\{t_1,\ldots,t_n\}}\cong
\ker(\nabla')^{an}|_{C\setminus\{t_1,\ldots,t_n\}}$.
By Proposition \ref{prop:elementary-transform},
There are a $(\bt,\bmu)$-parabolic connection
$(E_1,\nabla_1,\{(l_1)^{(i)}_j\})$
and a $(\bt,\bmu')$-parabolic connection
$(E'_1,\nabla'_1,\{(l'_1)^{(i)}_j\})$
such that $(E,\nabla,\{l^{(i)}_j\})\sim(E_1,\nabla_1,\{(l_1)^{(i)}_j\})$,
$(E',\nabla',\{(l')^{(i)}_j\})\sim(E'_1,\nabla'_1,\{(l'_1)^{(i)}_j\})$
and
$0\leq\mathrm{Re}(\mu^{(i)}_j)<1$, $0\leq\mathrm{Re}((\mu')^{(i)}_j)<1$
for any $i,j$.
Since
$\ker\nabla_1^{an}|_{C\setminus\{t_1,\ldots,t_n\}}\cong
\ker(\nabla'_1)^{an}|_{C\setminus\{t_1,\ldots,t_n\}}$,
we can see by [\cite{Deligne}, II, Proposition 5.4] that
$(E_1,\nabla_1)\cong(E'_1,\nabla'_1)$.
So we have
$(E,\nabla,\{l^{(i)}_j\})\sim(E_1,\nabla_1,\{(l_1)^{(i)}_j\})
\sim(E'_1,\nabla'_1,\{(l'_1)^{(i)}_j\})
\sim(E',\nabla',\{(l')^{(i)}_j\})$.
Thus we obtain the injectivity of $\mathrm{RH}$.

Unfortunately, we cannot expect an appropriate algebraic structure
on the moduli space of the elements
$[(E,\nabla,\{l^{(i)}_j\})]$ in the set $(\dag)$.
So we can recognize that it is natural to consider
the moduli space of isomorphism classes of
parabolic connections in the moduli theoretic
description of the Riemann-Hilbert correspondence.
However, we must consider a stability when we construct the moduli
of parabolic connections as an appropriate space.
So we set $M^{\balpha}_C(\bt,\blambda)$ as the moduli space of
$\balpha$-stable parabolic connections.
See Theorem \ref{moduli-exists} and Definition \ref{def-stability}
for the precise definition of $M^{\balpha}_C(\bt,\blambda)$.
Next we consider the moduli space
$RP_r(C,\bt)_{\ba}$ of certain equivalence classes of representations
of the fundamental group $\pi_1(C\setminus\{t_1,\ldots,t_n\},*)$.
Here two representations are equivalent if their semisimplifications
are isomorphic.
There is not an appropriate moduli space of isomorphism classes of
the representations of the fundamental group.
For a construction of a good moduli space containing all the representations,
we must consider such an equivalence relation.
See section \ref{main-thm} for the precise definition
of $RP_r(C,\bt)_{\ba}$.
The most crucial part of the main result (Theorem \ref{main-thm})
is that the morphism
$\RH:M^{\balpha}_C(\bt,\blambda) \lra RP_r(C,\bt)_{\ba}$
determined by the Riemann-Hilbert correspondence is a proper surjective
bimeromorphic morphism.
This theorem was proved in \cite{IIS-1} for $C=\mathbf{P}^1$ and $r=2$.
However there were certain difficulties to generalize this fact
to the case of general $C$ and $r$.

One of the most important motivation to consider such a theorem
is to consider an application to the geometric description of
the differential equation determined by the isomonodromic deformation.
See section \ref{isomonodromic deformation}
for the precise definition of the isomonodromic deformation.
This differential equation is said to be the Schlesinger equation
for $C=\mathbf{P}^1$, the Garnier equation for $C=\mathbf{P}^1$ and $r=2$,
and the Painlev\'e equation of type sixth for $C=\mathbf{P}^1$, $r=2$ and $n=4$.
M.~Jimbo, T.~Miwa and K.~Ueno give in \cite{MJU} and \cite{MJ} 
an explicit description of the Schlesinger equation.
In this paper, we will constrct a space $M$ with a morphism
$\pi:M\rightarrow T$ such that there is a differential equation
on $M$ with respect to a time variable $t\in T$
which is determined by the isomonodromic deformation.
If we fix a point $t_0\in T$, $\pi^{-1}(t_0)$ should become a
space of initial conditions of the differential equation
determined by the isomonodromic deformation.
Take any point $x\in\pi^{-1}(t_0)$ and consider
the analytic continuation $\gamma$ starting at
$x$ and $\pi\gamma$ coming back to the initial point $t_0$ such that
$\gamma$ is a solution of the differential equation
determined by the isomonodromic deformation.
Then $\gamma$ should come back to a point of
$\pi^{-1}(t_0)$.
We will generalize this property to the
geometric Painlev\'e property.
Roughly speaking, the geometric Painlev\'e property is the property of the space
where the differential equation is defined and
any analytic continuation of a solution of
the differential equation stays in the space.
See Definition \ref{def-gpp} for the precise definition of
the geometric Painlev\'e property.
In fact we will take $M$ as the moduli space
$M^{\balpha}_{\cC/T}(\tilde{\bt},r,d)$
of $\balpha$-stable parabolic connections
and $\pi:M^{\balpha}_{\cC/T}(\tilde{\bt},r,d)\rightarrow T$
as the structure morphism.
For example let us assume that $C=\BP^1$
and take an $\balpha$-stable parabolic connection
$x=(E,\nabla,\{l^{(i)}_j\})$ such that $E\cong{\mathcal O}_{\BP^1}^{\oplus r}$.
Let $\gamma:[0,1]\rightarrow M^{\balpha}_{\BP\times T/T}(\tilde{\bt},r,0)$
be a path starting at $x$ and $\pi\gamma$ ends at $\pi(x)$ such that
$\gamma$ is a solution of the differential equation
determined by the isomonodromic deformation.
Then the ending point $\gamma(1)$ corresponds to an $\balpha$-stable
parabolic connection $(E',\nabla',\{(l')^{(i)}_j\})$,
but $E'$ may not be trivial.
So we recognize that it is not enough to consider only
trivial vector bundle with a connection.
We should also consider a non-trivial vector bundle with a connection.

Our aim here is to construct a space where the isomonodromic deformation
is defined and satisfies the geometric Painlev\'e property over every
value of $\blambda$.
(It is not difficult to construct such a space over generic $\blambda$
but it is difficult to construct over special $\blambda$.)
In fact that space is nothing but the moduli space of $\balpha$-stable
parabolic connections and the result is given in Theorem \ref{gpp-thm},
which is essentially a corollary of Theorem \ref{main-thm}.
Here note that the properness of the morphism $\RH$ is very essential
in the proof of Theorem \ref{gpp-thm}.
The geometric Painlev\'e property immediately deduces the usual analytic
Painlev\'e property.
So we can say that Theorem \ref{gpp-thm} gives a most clear proof of
the Painlev\'e property of the isomonodromic deformation.
As is well-known, the solutions of the Painlev\'e equation have
the Painlev\'e property, which is in some sense the property characterizing
the Painlev\'e equation (and there were many proof of the Painlev\'e property).
So the usual analytic Painlev\'e property plays an important role
in the theory of the Painlev\'e equation,
but we can see from the definition that the `` geometric Painlev\'e property''
is much more important from the view point of the description of
the geometric picture of the isomonodromic deformation.

To prove the main results, this paper consists of several sections.
In section \ref{moduli of s-p-c}, we prove the existence of
the moduli space of stable parabolic connections.
An algebraic moduli space of parabolic connections was essentially
considered by D.~Arinkin and S.~Lysenko in \cite{Arinkin-Lysenko1},
\cite{Arinkin-Lysenko2} and \cite{Arinkin1} and they showed that
the moduli space of parabolic connections on $\mathbf{P}^1$ of rank $2$
with $n=4$ for generic $\blambda$ is isomorphic to the space of
initial conditions of the Painlev\'e equation of type sixth constructed
by K.~Okamoto (\cite{Okamoto}).
For special $\blambda$, we should consider a certain stability condition
to construct an appropriate moduli space of parabolic connections.
In the case of $C=\mathbf{P}^1$ and $r=2$,
K.~Iwasaki, M.-H.~Saito and the author already considered in 
\cite{IIS-1} the moduli space of stable parabolic connections
and they proved in \cite{IIS-2} that the moduli space of
stable parabolic connections on $\mathbf{P}^1$ of rank $2$ with
$n=4$ is isomorphic to the space of initial conditions of
the Painlev\'e VI equation constructed by K.~Okamoto for all $\blambda$.
An analytic construction of the moduli space of stable parabolic connections
for general $C$ and $r=2$ was given by H.~Nakajima in \cite{N}.
However, in our aim, the algebraic construction of the moduli space
is necessary.
The morphism $\RH$ determined by the Riemann-Hilbert correspondence
is quite transcendental, which is explicitly shown in the case of
$C=\mathbf{P}^1$, $r=2$ and $n=4$ in \cite{IIS-2}.
This statement make sense only if we construct the moduli spaces
$M^{\balpha}_C(\bt,\blambda)$ and $RP_r(C,\bt)_{\ba}$
algebraically.
Moduli of logarithmic connections without parabolic structure
is constructed in \cite{Nitsure}.
For the proof of Theorem \ref{moduli-exists}, we do not use the
method in \cite{Simp-I}.
We construct $M^{\balpha}_{\cC/T}(\tilde{\bt},r,d)$ as a subscheme of
the moduli space of parabolic $\Lambda^1_{D}$-triples
constructed in \cite{IIS-1}.
We also use this embedding in the proof of the properness of $\RH$
in Theorem \ref{main-thm}.

In section \ref{irreducibility of m-s}, we prove that the moduli space
of stable parabolic connections is an irreducible variety.
For the proof, we need certain complicated calculations and
this part is a new difficulty, which did not appear in \cite{IIS-1}.

In section \ref{riemann-hilbert}, we consider the morphism $\RH$
determined by the Riemann-Hilbert correspondence and prove the
surjectivity and properness of $\RH$.
This is the essential part of the proof of Theorem \ref{main-thm}.
Notice that for generic $\blambda$, parabolic connection is irreducible
and so all parabolic connections are stable.
Moreover we can easily see that the moduli space
$M^{\balpha}_C(\bt,\blambda)$
of parabolic connections for generic $\blambda$ is analytically isomorphic
via the Riemann-Hilbert correspondence to the moduli space $RP_r(C,\bt)_{\ba}$
of representations of the fundamental group
$\pi_1(C\setminus\{t_1,\ldots,t_n\},*)$.
However, for special $\blambda$, several parabolic connections
may not be stable and the stability condition depends on $\balpha$.
So we can see that the surjectivity of $\RH$ is not trivial at all
for special $\blambda$, because for a point $[\rho]\in RP_r(C,\bt)_{\ba}$,
we must find a parabolic connection corresponding to $[\rho]$ which is 
{\it $\balpha$-stable}.
Moreover, the moduli space $M^{\balpha}_C(\bt,\blambda)$
of $\balpha$-stable parabolic connections is smooth, but the moduli space
$RP_r(C,\bt)_{\ba}$ of representations of the fundamental group
becomes singular.
So the morphism
$\RH:M^{\balpha}_C(\bt,\blambda)\rightarrow RP_r(C,\bt)_{\ba}$
is more complicated in the case of special $\blambda$ than 
the case of generic $\blambda$.
We first prove the surjectivity of $\RH$ in Proposition \ref{surjectivity}.
In this proof, we use the Langton's type theorem in the case
of parabolic connections.
This idea was already used in \cite{IIS-1}.
Secondly, we prove in Proposition \ref{compact}
that every fiber of $\RH$ is compact by using
an embedding of $M^{\balpha}_C(\bt,\blambda)$ to a certain compact
moduli space.
In this proof, we can not use the idea given in \cite{IIS-1} and
make new technique again.
Finally we obtain the properness of $\RH$
by the lemma given by A.~Fujiki.

In section \ref{symplectic structure}, we construct a canonical symplectic
form on the moduli space $M^{\balpha}_C(\bt,\blambda)$
of $\balpha$-stable parabolic connections.
Combined with the the fact that $\RH$ gives an analytic resolution of
singularities of $RP_r(C,\bt)_{\ba}$,
we can say that $RP_r(C,\bt)_{\ba}$ has symplectic singularities
(for special $\ba$) and $\RH$ gives a symplectic resolution
of singularities.
For the case of $r=2$, H.~Nakajima constructed the moduli space
$M^{\balpha}_C(\bt,\blambda)$ as a hyper-K\"{a}hler manifold and
it obviously has a holomorphic symplectic structure.
Such a construction for general $r$ is also an important problem,
though we do not treat it here.

In section \ref{isomonodromic deformation},
we first give in Proposition \ref{splitting} an algebraic construction
of the differential equation on $M^{\balpha}_{\cC/T}(\tilde{\bt},r,d)$
determined by the isomonodromic deformation.
Finally we complete the proof of Theorem \ref{gpp-thm}.

{\bf Acknowledgments.}
The author would like to thank Professors Masa-Hiko Saito and Katsunori Iwasaki
for giving him the problem solved in this paper and also for valuable
long discussions.
He also would like to thank Professors Takuro Mochizuki and Takeshi Abe
for valuable comments and discussions on the proof given
in section \ref{riemann-hilbert}.

\section{Main results}\label{main results}

Let $C$ be a smooth projective curve of genus $g$.
We put
\[
 T_n :=\left.\left\{ (t_1,\ldots,t_n)\in\overbrace{C\times\cdots\times C}^n
 \right| \text{$t_i\neq t_j$ for $i\neq j$} \right\}
\]
for a positive integer $n$.
For integers $d,r$ with $r>0$, we put
\[
 \Lambda^{(n)}_r(d):=
 \left\{ (\lambda^{(i)}_j)^{1\leq i\leq n}_{0\leq j\leq r-1}
 \in \C^{nr}
 \left|
 d+\sum_{i,j}\lambda^{(i)}_j=0
 \right\}\right..
\]
Take a member $\bt=(t_1,\ldots,t_n)\in T_n$ and
$\blambda=(\lambda^{(i)}_j)_{1\leq i\leq n,0\leq j\leq r-1}
\in\Lambda^{(n)}_r(d)$.

\begin{Definition}\rm
We say $(E,\nabla,\{l^{(i)}_*\}_{1\leq i\leq n})$
a $(\bt,\blambda)$-parabolic connection of rank $r$ if
\begin{enumerate}
\item $E$ is a rank $r$ algebraic vector bundle on $C$,
\item $\nabla: E \ra E\otimes\Omega_C^1(t_1+\cdots+t_n)$
 is a connection, that is, $\nabla$ is a homomorphism of sheaves
 satisfying $\nabla(fa)=a\otimes df +f\nabla(a)$ for
 $f\in\cO_C$ and $a\in E$,
and
\item for each $t_i$,
 $l^{(i)}_*$ is a filtration
 $E|_{t_i}=l^{(i)}_0\supset l^{(i)}_1
 \supset\cdots\supset l^{(i)}_{r-1}\supset l^{(i)}_r=0$
 such that $\dim(l^{(i)}_j/l^{(i)}_{j+1})=1$ and
 $(\res_{t_i}(\nabla)-\lambda^{(i)}_j\mathrm{id}_{E|_{t_i}})
 (l^{(i)}_j)\subset l^{(i)}_{j+1}$
 for $j=0,\ldots,r-1$.
\end{enumerate}
\end{Definition}

\begin{Remark}\rm
By condition (3) above, we have
\[
 \deg E=\deg(\det(E))=-\sum_{i=1}^n\res_{t_i}(\nabla_{\det E})
 =-\sum_{i=1}^n\sum_{j=0}^{r-1}\lambda^{(i)}_j=d.
\]
\end{Remark}

Take rational numbers
\[
 0<\alpha^{(i)}_1<\alpha^{(i)}_2<\cdots<\alpha^{(i)}_r<1
\]
for $i=1,\ldots,n$ satisfying
$\alpha^{(i)}_j\neq\alpha^{(i')}_{j'}$ for $(i,j)\neq(i',j')$.
We choose $\balpha=(\alpha^{(i)}_j)$ sufficiently generic.

\begin{Definition}\label{def-stability}\rm
A parabolic connection $(E,\nabla,\{l^{(i)}_*\}_{1\leq i\leq n})$
is $\balpha$-stable (resp.\ $\balpha$-semistable) if
for any proper nonzero subbundle $F\subset E$
satisfying $\nabla(F)\subset F\otimes\Omega_C^1(t_1+\cdots+t_n)$,
the inequality
\begin{gather*}
 \frac{\deg F + \sum_{i=1}^n\sum_{j=1}^r
 \alpha^{(i)}_j\dim((F|_{t_i}\cap l^{(i)}_{j-1})/(F|_{t_i}\cap l^{(i)}_j))}
 {\rank F} \\
 \genfrac{}{}{0pt}{}{<}{(\text{resp. $\leq$})}
 \frac{\deg E + \sum_{i=1}^n\sum_{j=1}^r
 \alpha^{(i)}_j \dim(l^{(i)}_{j-1}/l^{(i)}_j)}{\rank E}
\end{gather*}
holds.
\end{Definition}

\begin{Remark}\label{generic-polarization}
Assume that $\balpha=(\alpha^{(i)}_j)$ satisfies the condition:
For any integer $r'$ with $0<r'<\rank(E)$ and for any
$(\epsilon^{(i)}_j)^{i=1,\ldots,n}_{j=1,\ldots,\rank(E)}$
with $\epsilon^{(i)}_j\in\{0,1\}$
and $\sum_{j=1}^r\epsilon^{(i)}_j=r'$ for any $i$,
\[
 \sum_{i=1}^n\sum_{j=1}^{\rank(E)} \alpha^{(i)}_j(r'-\rank(E)\epsilon^{(i)}_j) \notin \mathbf{Z}.
\]
Then a parabolic connection $(E,\nabla,\{l^{(i)}_j\})$ is $\balpha$-stable
if and only if $(E,\nabla,\{l^{(i)}_j\})$ is $\balpha$-semistable.
\end{Remark}

Let $T$ be a smooth algebraic scheme which is a certain covering 
of the moduli stack of $n$-pointed smooth projective curves of genus $g$
over $\C$ and take a universal family $(\cC,\tilde{t}_1,\ldots,\tilde{t}_n)$
over $T$.

\begin{Theorem}\label{moduli-exists}
There exists a relative fine moduli scheme
\[
 M^{\balpha}_{\cC/T}(\tilde{\bt},r,d)\rightarrow T\times\Lambda^{(n)}_r(d)
\]
of $\balpha$-stable parabolic connections of rank $r$ and degree $d$,
which is smooth and quasi-projective.
The fiber $M^{\balpha}_{\cC_x}(\tilde{\bt}_x,\blambda)$
over $(x,\blambda)\in T\times\Lambda^{(n)}_r(d)$ is the moduli space of
$\balpha$-stable $(\tilde{\bt}_x,\blambda)$-parabolic connections
whose dimension is $2r^2(g-1)+nr(r-1)+2$ if it is non-empty.
\end{Theorem}

\begin{Definition}\rm
Take an element $\blambda\in\Lambda^{(n)}_r(d)$.
We call $\blambda$ special if
\begin{enumerate}
\item $\lambda^{(i)}_j-\lambda^{(i)}_k\in\Z$ for some $i$ and $j\neq k$, or
\item there exist an integer $s$ with $1<s<r$ and a subset
 $\{j^{i}_1,\ldots,j^{i}_s\}\subset\{0,\ldots,r-1\}$
 for each $1\leq i\leq n$ such that
 \[
  \sum_{i=1}^n\sum_{k=1}^s \lambda^{(i)}_{j^i_k}\in\Z.
 \]
\end{enumerate}
We call $\blambda$ resonant if it satisfies the condition (1) above.
Let $(E,\nabla,\{l^{(i)}_j\})$ be a $(\bt,\blambda)$-parabolic connection
such that $\blambda$ does not satisfy the condition (2) above.
Then $(E,\nabla,\{l^{(i)}_j\})$ is irreducible.
Here we say a $(\tilde{\bt}_x,\blambda)$-connection
$(E,\nabla_E,\{l^{(i)}_j\})$ reducible if there is a non-trivial
subbundle $0\neq F\subsetneq E$ such that
$\nabla_E(F)\subset F\otimes\Omega^1_{\cC_x}
((\tilde{t}_1)_x+\cdots+(\tilde{t}_n)_x)$.
We say $(E,\nabla_E,\{l^{(i)}_j\})$ irreducible if it is not reducible.
We call $\blambda\in\Lambda^{(n)}_r(d)$ generic if it is not special.
\end{Definition}

Fix a point $x\in T$.
Then the fundamental group
$\pi_1(\cC_x\setminus\{(\tilde{t}_1)_x,\ldots,(\tilde{t}_n)_x\},*)$
is generated by cycles $\alpha_1,\beta_1,\ldots,\alpha_g,\beta_g$
and loops $\gamma_i$ around $t_i$ for $1\leq i\leq n$
whose relation is given by
\[
 \prod_{j=1}^g \alpha_j^{-1}\beta_j^{-1}\alpha_j\beta_j
 \prod_{i=1}^n \gamma_i =1.
\]
So the fundamental group is isomorphic to a free group generated by
$2g+n-1$ free generators.
Then the space
\[
 \Hom\left(\pi_1(\cC_x\setminus\{(\tilde{t}_1)_x,\ldots,
 (\tilde{t}_n)_x\},*),GL_r(\C)\right)
\]
of representations of the fundamental group becomes an affine variety
isomorphic to $GL_r(\C)^{2g+n-1}$
and $GL_r(\C)$ acts on this space by the adjoint action.
We define
\[
 RP_r((\cC)_x,\tilde{\bt}_x):=
 \Hom\left(\pi_1(\cC_x\setminus\{(\tilde{t}_1)_x,\ldots,
 (\tilde{t}_n)_x\},*),GL_r(\C)\right)//GL_r(\C)
\]
as a categorical quotient.
If we put
\[
 \cA^{(n)}_r:=\left\{
 \ba=(a^{(i)}_j)^{1\leq i\leq n}_{0\leq j\leq r-1} \in \C^{nr}
 \left| a^{(1)}_0 a^{(2)}_0 \cdots a^{(n)}_0=(-1)^{rn}
 \right\}\right.,
\]
then we can define a morphism
\begin{gather*}
 RP_r((\cC)_x,\tilde{\bt}_x)\lra
 \cA^{(n)}_r \\
 \quad [\rho]\mapsto \ba=(a^{(i)}_j)
\end{gather*}
by the relation
\[
 \det(XI_r-\rho(\gamma_i))=X^r+a^{(i)}_{r-1}X^{r-1}+\cdots+a^{(i)}_0,
\]
where $X$ is an indeterminate
and $I_r$ is the identity matrix of size $r$.
Note that for any $GL_r(\C)$-representation $\rho$ of
$\pi_1(\cC_x\setminus\{(\tilde{t}_1)_x,\ldots,(\tilde{t}_n)_x\},*)$,
\begin{align*}
 &\det(\rho(\gamma_1))\cdots\det(\rho(\gamma_n)) \\
 &=
 \det(\rho(\beta_g))^{-1}\det(\rho(\alpha_g))^{-1}\det(\beta_g)\det(\alpha_g)
 \cdots
 \det(\rho(\beta_1))^{-1}\det(\rho(\alpha_1))^{-1}
 \det(\beta_1)\det(\alpha_1) \\
 &=1
\end{align*}
and so the equation
$a^{(1)}_0 a^{(2)}_0 \cdots a^{(n)}_0=(-1)^{rn}$
should be satisfied.
 
We can define a relative moduli space
$RP_r(\cC,\tilde{\bt})=\coprod_{x\in T}RP_r(\cC_x,\tilde{\bt}_x)$
of representations.
Let $\tilde{T}\rightarrow T$
be the universal covering.
Then $RP_r(\cC,\tilde{\bt})\times_T\tilde{T}$
becomes a trivial fibration over $\tilde{T}$.
So we can define a morphism
\[
 RP_r(\cC,\tilde{\bt})\times_T\tilde{T}
 \lra \tilde{T}\times\cA^{(n)}_r.
\]
We define a morphism
\begin{equation}\label{local-riemann-hilbert}
 rh:\Lambda^{(n)}_r(d) \ni\blambda\mapsto\ba\in \cA^{(n)}_r
\end{equation}
by
\[
 \prod_{j=0}^{r-1}\left(X-\exp(-2\pi\sqrt{-1}\lambda^{(i)}_j)\right)
 =X^r+a^{(i)}_{r-1}X^{r-1}+\cdots+a^{(i)}_0.
\]
For each member
$(E,\nabla,\{l^{(i)}_j\})\in
M^{\balpha}_{\cC_x}(\tilde{\bt}_x,\blambda)$,
$\ker(\nabla^{an}|_{\cC_x\setminus\{(\tilde{t}_1)_x,\ldots,(\tilde{t}_n)_x\}})$
becomes a local system on
$\cC_x\setminus\{(\tilde{t}_1)_x,\ldots,(\tilde{t}_n)_x\}$,
where $\nabla^{an}$ means the analytic connection corresponding to $\nabla$.
The local system
$\ker(\nabla^{an}|_{\cC_x\setminus\{(\tilde{t}_1)_x,\ldots,(\tilde{t}_n)_x\}})$
corresponds to a representation of
$\pi_1(\cC_x\setminus\{(\tilde{t}_1)_x,\ldots,(\tilde{t}_n)_x\},*)$.
So we can define a morphism
\[
 \RH_{(x,\blambda)}:M^{\balpha}_{\cC_x}(\tilde{\bt}_x,\blambda)\lra
 RP_r(\cC_x,\tilde{\bt}_x)_{\ba},
\]
where $\ba=rh(\blambda)$.
$\{\RH_{(x,\blambda)}\}$ induces a morphism
\begin{equation}\label{global-RH}
 \RH: M^{\balpha}_{\cC/T}(\tilde{\bt},r,d)\times_T\tilde{T}
 \lra RP_r(\cC,\tilde{\bt})\times_T\tilde{T}
\end{equation}
which makes the diagram
\[
 \begin{CD}
  M^{\balpha}_{\cC/T}(\tilde{\bt},r,d))\times_T\tilde{T} @> \RH >>
  RP_r(\cC,\tilde{\bt})\times_T\tilde{T} \\
  @VVV @VVV \\
  \tilde{T}\times\Lambda^{(n)}_r(d) @> \mathrm{id}\times rh >>
  \tilde{T}\times\cA^{(n)}_r
 \end{CD}
\]
commute.

The following is the main theorem whose proof is given in section 
\ref{riemann-hilbert}.

\begin{Theorem}\label{main-thm}
 Assume $\balpha$ is so generic that $\balpha$-stable $\Leftrightarrow$
 $\balpha$-semistable (see Remark \ref{generic-polarization}).
 Moreover we assume $rn-2r-2>0$ if $g=0$, $n>1$ if $g=1$
 and $n\geq 1$ if $g\geq 2$.
 Then the morphism
\[
 \RH:M^{\balpha}_{\cC/T}(\tilde{\bt},r,d)\times_T\tilde{T} \lra
 (RP_r(\cC,\tilde{\bt})\times_T\tilde{T})\times_{\cA^{(n)}_r}\Lambda^{(n)}_r(d)
\]
induced by (\ref{global-RH}) is a proper surjective
bimeromorphic morphism.
Combined with Proposition \ref{symplectic-form}, we can say that
for each $(x,\blambda)\in T\times\Lambda^{(n)}_r(d)$,
\begin{enumerate}
\item
 $\RH_{(x,\blambda)}: M^{\balpha}_{\cC_x}(\tilde{\bt}_x,\blambda)\lra
 RP_r(\cC_x,\tilde{\bt}_x)_{\ba}$
 is an analytic isomorphisms of symplectic varieties
 for generic $\blambda$ and
\item
 $\RH_{(x,\blambda)}: M^{\balpha}_{\cC_x}(\tilde{\bt}_x,\blambda)\lra
 RP_r(\cC_x,\tilde{\bt}_x)_{\ba}$
gives an analytic symplectic resolution of singularities of
$RP_r(\cC_x,\tilde{\bt}_x)_{\ba}$ for special $\blambda$.
\end{enumerate}
Here we put $\ba=rh(\blambda)$.
\end{Theorem}

\begin{Remark}\rm
(1) For the case of $r=1$, $(\RH)_{(x,\blambda)}$ is an isomorphism
for any $(x,\blambda)\in T\times\Lambda^{(n)}_1(d)$
and $RP_1(\cC_x,\tilde{\bt}_x)_{\ba}$ is smooth for any $\ba$. \\
(2) In Theorem \ref{main-thm} we consider the case $n\geq 1$.
As is stated in [\cite{Simp-II}, Proposition 7.8],
the Riemann-Hilbert correspondence for the case of $n=0$ gives an analytic
isomorphism between the moduli space of integrable connections and
the moduli space of representations of the fundamental group,
but in that case the moduli space of integrable connections may be singular.
\end{Remark}

Since $RP_r(\cC,\tilde{\bt})\times_T\tilde{T}\ra\tilde{T}$
is a trivial fibration, we can consider the set of constant sections
\[
 {\mathcal F}_R:=
 \left\{
 \sigma:\tilde{T}\ra RP_r(\cC,\tilde{\bt})\times_T\tilde{T}
 \right\}.
\]
As is stated in Remark \ref{singular-locus},
\[
 RP_r(C,\bt)_{\ba}^{\mathrm{sing}}=
 \left\{ [\rho]\in RP_r(C,\bt)_{\ba} \left|
 \begin{array}{l}
  \text{$\rho$ is reducible or} \\
  \text{$\dim\left(
  \ker\left(\rho(\gamma_i)-\exp(-2\pi\sqrt{-1}\lambda^{(i)}_j)I_r
  \right)\right)
  \geq 2$ for some $i,j$}
 \end{array}
 \right\}\right.
\]
is just the singular locus of $RP_r(C,\bt)_{\ba}$.
If we put
$RP_r(\cC_x,\tilde{\bt}_x)_{\ba}^{\sharp}:=
RP_r(\cC_x,\tilde{\bt}_x)_{\ba}\setminus
RP_r(C,\bt)_{\ba}^{\mathrm{sing}}$,
then
$\RH|_{M^{\balpha}_{\cC_x}(\tilde{\bt}_x,\blambda)^{\sharp}}:
M^{\balpha}_{\cC_x}(\tilde{\bt}_x,\blambda)^{\sharp}\rightarrow
RP_r(\cC_x,\tilde{\bt}_x)_{\ba}^{\sharp}$
becomes an isomorphism, where we put
$M^{\balpha}_{\cC_x}(\tilde{\bt}_x,\blambda)^{\sharp}:=
\RH_{(x,\blambda)}^{-1}(RP_r(\cC_x,\tilde{\bt}_x)_{\ba}^{\sharp})$.
Put
\[
 (RP_r(\cC,\tilde{\bt})\times_T\tilde{T})^{\sharp}
 =\coprod_{(x,\ba)\in \tilde{T}\times\cA^{(n)}_r}
 RP_r(\cC_x,\tilde{\bt}_x)_{\ba}^{\sharp}.
\]
Then the restriction ${\mathcal F}_R^{\sharp}$ of ${\mathcal F}_R$
to $(RP_r(\cC,\tilde{\bt})\times_T\tilde{T})^{\sharp}$ gives a foliation on
$(RP_r(\cC,\tilde{\bt})\times_T\tilde{T})^{\sharp}$.
The pull back
$\tilde{\mathcal F}^{\sharp}_M:=\RH^{-1}({\mathcal F}^{\sharp}_R)$
determines a foliation on
$(M^{\balpha}_{\cC/T}(\tilde{\bt},r,d)\times_T\tilde{T})^{\sharp}$
where $(M^{\balpha}_{\cC/T}(\tilde{\bt},r,d)\times_T\tilde{T})^{\sharp}=
\RH^{-1}((RP_r(\cC,\tilde{\bt})\times_T\tilde{T})^{\sharp})$.
This foliation corresponds to a subbundle of the analytic tangent bundle
$\Theta^{an}_{(M^{\balpha}_{\cC/T}(\tilde{\bt},r,d)
\times_T\tilde{T})^{\sharp}}$
determined by a splitting
\[
 D^{\sharp}:
 \tilde{\pi}^*(\Theta_{\tilde{T}}^{an})
 \ra \Theta^{an}_{(M^{\balpha}_{\cC/T}
 (\tilde{\bt},r,d)\times_T\tilde{T})^{\sharp}}
\]
of the analytic tangent map
\[
 \Theta^{an}_{(M^{\balpha}_{\cC/T}
 (\tilde{\bt},r,d)\times_T\tilde{T})^{\sharp}}\lra
 \tilde{\pi}^*\Theta^{an}_{\tilde{T}}  \ra 0,
\]
where
$\tilde{\pi}:(M^{\balpha}_{\cC/T}(\tilde{\bt},r,d)
\times_T\tilde{T})^{\sharp} \ra \tilde{T}$
is the projection.
We will show in Proposition \ref{splitting} that this splitting
$D^{\sharp}$ is in fact induced by a splitting
\[
 D:\pi^*(\Theta_T)
 \ra \Theta_{M^{\balpha}_{\cC/T}(\tilde{\bt},r,d)}
\]
of the algebraic tangent bundle,
where $\pi:M^{\balpha}_{\cC/T}(\tilde{\bt},r,d)\ra T$ is the projection.

In order to explain the concept of the geometric Painlev\'e property,
let us review here several terminologies introduced by K.~Iwasaki
in \cite{IISA}.

\begin{Definition}{\bf (\cite{IISA})}
 A time-dependent dynamical system $(M,\cF)$ is a smooth fibration
$\pi:M\ra T$ of complex manifolds together with a complex foliation
$\cF$ on $M$ that is transverse to each fiber $M_t=\pi^{-1}(t)$,
$t\in T$.
The total space $M$ is referred to as the phase space,
while the base space $T$ is called the space of time-variables.
Moreover, the fiber $M_t$ is called the space of initial conditions
at time $t$.
\end{Definition}

\begin{Definition}\label{def-gpp}{\bf (\cite{IISA})}
 A time-dependent dynamical system $(M,\cF)$ is said to have
the geometric Painlev\'e property if for any path $\gamma$ in $T$ and
any point $p\in M_t$, where $t$ is the initial point of $\gamma$,
there exists a unique $\cF$-horizontal lift $\tilde{\gamma_p}$ of $\gamma$
with initial point $p$.
Here a curve in $M$ is said to be $\cF$-horizontal if it lies in
a leaf of $\cF$.
\end{Definition}

$D^{\sharp}(\tilde{\pi}^*\Theta^{an}_{\tilde{T}})$ obviously satisfies
the integrability condition
\[
 [D^{\sharp}(\tilde{\pi}^*\Theta^{an}_{\tilde{T}}),
 D^{\sharp}(\tilde{\pi}^*\Theta^{an}_{\tilde{T}})]
 \subset D^{\sharp}(\tilde{\pi}^*\Theta^{an}_{\tilde{T}})
\]
because it corresponds to the foliation
$\tilde{\mathcal F}^{\sharp}_M$.
Since
\[
 \codim_{(M^{\balpha}_{\cC/T}(\tilde{\bt},r,d)\times_T\tilde{T})^{\sharp}}
 \left((M^{\balpha}_{\cC/T}(\tilde{\bt},r,d)\times_T\tilde{T})\setminus
 (M^{\balpha}_{\cC/T}(\tilde{\bt},r,d)\times_T\tilde{T})^{\sharp}\right)\geq 2,
\]
we can see that
$D(\pi^*(\Theta_T))$ also satisfies the integrability condition
\[
 [D(\pi^*(\Theta_T)),D(\pi^*(\Theta_T))]
 \subset D(\pi^*(\Theta_T)).
\]
Then the foliation $\tilde{\mathcal F}^{\sharp}_M$
extends to a foliation $\tilde{\mathcal F}_M$
on $M^{\balpha}_{\cC/T}(\tilde{\bt},r,d)\times_T\tilde{T}$.
We can see that $\tilde{\mathcal F}_M$ descends to a foliation
${\mathcal F}_M$ on $M^{\balpha}_{\cC/T}(\tilde{\bt},r,d)$
which corresponds to the subbundle
$D(\pi^*(\Theta_T))\subset
\Theta_{M^{\balpha}_{\cC/T}(\tilde{\bt},r,d)}$.
By construction, we can see that
$\left(M^{\balpha}_{\cC/T}(\tilde{\bt},r,d),\mathcal F_M\right)$
becomes a time-dependent dynamical system with base space $T$.
As a corollary of Theorem \ref{main-thm}, we can prove the following
theorem whose proof is given in section \ref{isomonodromic deformation}.

\begin{Theorem}\label{gpp-thm}
 Under the assumption of Theorem \ref{main-thm},
 the time-dependent dynamical system
 $\left(M^{\balpha}_{\cC/T}(\tilde{\bt},r,d),\mathcal F_M\right)$
 determined by the differential system $D(\pi^*\Theta_T)$
 has the geometric Painlev\'e property.
\end{Theorem}

\section{Elementary transform of parabolic connections}\label{elementary-transform-m}

For a parabolic connection $(E,\nabla,\{l^{(i)}_j\})$ on $(C,\bt)$,
Put
\[
 E':=\ker(E\longrightarrow E|_{t_p}/l^{(p)}_{q}).
 \]
Then $\nabla$ induces a connection
$\nabla':E'\rightarrow E'\otimes\Omega^1_C(t_1+\cdots+t_n)$
such that the diagram
\[
 \begin{CD}
  E' @>\nabla'>> E'\otimes\Omega^1_C(t_1+\cdots+t_n) \\
  @VVV @VVV \\
  E@>\nabla>> E\otimes\Omega^1_C(t_1+\cdots+t_n)
 \end{CD}
\]
commutes.
Let $E'|_{t_p}=(l')^{(p)}_0\supset (l')^{(p)}_1\supset\cdots\supset(l')^{(p)}_{r-q}$
be the inverse image of the filtration
$l^{(p)}_q\supset\cdots\supset l^{(p)}_r=0$
by the linear map $E'|_{t_p}\rightarrow E|_{t_p}$
and $(l')^{(p)}_{r-q}\supset\cdots\supset(l')^{(p)}_r=0$
be the image of the filtration
$E(-t_p)|_{t_p}=l^{(p)}_0\otimes{\mathcal O}(-t_p)
\supset\cdots\supset l^{(p)}_q\otimes{\mathcal O}(-t_p)$
by the linear map
$E(-t_p)|_{t_p}\rightarrow E'|_{t_p}$.
Then $(l')^{(p)}_*$ becomes a filtration of $E'|_{t_p}$.
For $i\neq p$, we put $(l')^{(i)}_j=l^{(i)}_j$ for any $j$.
Then $(E',\nabla',\{(l')^{(i)}_j\})$ becomes a parabolic connection
on $(C,\bt)$.
We say $(E',\nabla',\{(l')^{(i)}_j\})$ the elementary transform of
$(E,\nabla,\{l^{(i)}_j\})$ along $t_p$ by $l^{(p)}_q$.

\begin{Definition}\rm
 We define a functor ${\mathcal M}_{{\mathcal C}/T}(\tilde{\bf t},r,d)$
 of the category of locally noetherian schemes over $T$ to the category of sets by
 \[
  {\mathcal M}_{{\mathcal C}/T}(\tilde{\bf t},r,d)(S)
  =\left\{(E,\nabla,\{l^{(i)}_j\})\right\}/\sim
 \]
for a locally noetherian scheme $S$ over $T$, where
\begin{enumerate}
 \item $E$ is a vector bundle on ${\mathcal C}_S$ of rank $r$,
 \item $\nabla:E\rightarrow E\otimes
 \Omega^1_{{\mathcal C}_S/S}((\tilde{t}_1)_S+\cdots+(\tilde{t}_n)_S)$
 is a relative connection, that is, $\nabla$ is a morphism of sheaves such that
 $\nabla(fa)=a\otimes df+ f\nabla(a)$ for $f\in{\mathcal O}_{{\mathcal C}_S}$
 and $a\in E$,
 \item $E|_{(\tilde{t}_i)_S}=l^{(i)}_0\supset\cdots\supset l^{(i)}_r=0$
 is a filtration by subbundles on $(\tilde{t}_i)_S$ such that
 $\res_{(\tilde{t}_i)_S}(\nabla)(l^{(i)}_j)\subset l^{(i)}_j$ for $i=1,\ldots,n$
 and $j=0,1,\ldots,r-1$, and
 \item for any $s\in S$, $E\otimes k(s)$ is of degree $d$ and
 $\dim\left(\left(l^{(i)}_j/l^{(i)}_{j+1}\right)\otimes k(s)\right)=1$ for any $i,j$.
\end{enumerate}
Here $(E,\nabla,\{l^{(i)}_j\})\sim(E',\nabla',\{(l')^{(i)}_j\})$
if there is a line bundle ${\mathcal L}$ on $S$ such that
$(E,\nabla,\{l^{(i)}_j\})\cong(E',\nabla',\{(l')^{(i)}_j\})\otimes{\mathcal L}$.
\end{Definition}

\begin{Remark}\rm
 Take any member
 $(E,\nabla,\{l^{(i)}_j\})\in{\mathcal M}_{{\mathcal C}/T}(\tilde{\bt},r,d)(S)$.
 Then $\res_{(\tilde{t}_i)_S}(\nabla)$ induces a morphism
 $\lambda^{(i)}_j:l^{(i)}_j/l^{(i)}_{j+1}\rightarrow l^{(i)}_j/l^{(i)}_{j+1}$
 and $\blambda=(\lambda^{(i)}_j)\in\Lambda^{(n)}_r(d)(S)$.
 So we obtain a canonical morphism
 \[
  \tilde{\blambda}:
  {\mathcal M}_{{\mathcal C}/T}(\tilde{\bt},r,d)
  \longrightarrow h_{\Lambda^{(n)}_r(d)}.
 \]
 Since there is a canonical morphism
 ${\mathcal M}_{{\mathcal C}/T}(\tilde{\bt},r,d)\rightarrow h_T$,
 we obtain a morphism
 \[
  (\tilde{x},\tilde{\blambda}): {\mathcal M}_{{\mathcal C}/T}(\tilde{\bt},r,d)
  \longrightarrow h_T\times h_{\Lambda^{(n)}_r(d)}=h_{T\times\Lambda^{(n)}_r(d)}.
 \]
\end{Remark}

\begin{Definition}\rm
For $\mathbf{C}$-valued points $\blambda\in\Lambda^{(n)}_r(d)(\mathbf{C})$,
$x\in T(\mathbf{C})$,
$(x,\blambda)$ can be considered as a morphism
$h_{\Spec\mathbf{C}} \rightarrow h_T\times h_{\Lambda^{(n)}_r(d)}$.
Then we put
\[
 {\mathcal M}_{{\mathcal C}_x}(\tilde{\bt}_x,\blambda)=
 {\mathcal M}_{{\mathcal C}/T}(\tilde{\bt},r,d)
 \times_{h_{T\times\Lambda^{(n)}_r(d)}}h_{\Spec\mathbf{C}}.
\]
${\mathcal M}_{{\mathcal C}_x}(\tilde{\bt}_x,\blambda)$
can be considered as a functor of the category
of locally noetherian schemes over $\mathbf{C}$ to the category of sets.
\end{Definition}

For a locally noetherian scheme $S$ over $T$,
take any member
$(E,\nabla,\{l^{(i)}_j\})\in{\mathcal M}_{{\mathcal C}/T}(\tilde{\bt},r,d)(S)$.
For $(p,q)$ with $1\leq p\leq n$ and $0\leq q\leq r$, we put
\[
 E':=\ker\left( E\longrightarrow (E|_{(\tilde{t}_p)_S}/l^{(p)}_q)\right).
\]
Then $\nabla$ induces a relative connection
\[
 \nabla':E'\longrightarrow E'\otimes
 \Omega^1_{{\mathcal C}_S/S}((\tilde{t}_1)_S+\cdots+(\tilde{t}_n)_S)
\]
such that the diagram
\[
 \begin{CD}
  E' @>\nabla'>> 
  E'\otimes\Omega^1_{{\mathcal C}_S/S}((\tilde{t}_1)_S+\cdots+(\tilde{t}_n)_S) \\
  @VVV @VVV \\
  E @>\nabla>>
  E\otimes\Omega^1_{{\mathcal C}_S/S}((\tilde{t}_1)_S+\cdots+(\tilde{t}_n)_S) 
 \end{CD}
\]
is commutative.
Note that the sequence
\[
 E(-(\tilde{t}_p)_S)|_{(\tilde{t}_p)_S}\stackrel{\pi}\longrightarrow
 E'|_{(\tilde{t}_p)_S}\stackrel{\iota}\longrightarrow E|_{(\tilde{t}_p)_S}
\]
is an exact sequence.
$\im\iota=l^{(p)}_q$ is a subbundle of $E|_{(\tilde{t}_p)_S}$ and
$\im\pi=\ker\iota$ is a subbundle of $E'|_{(\tilde{t}_p)_S}$.
Moreover
$\ker\pi=l^{(p)}_q\otimes{\mathcal O}_{{\mathcal C}_S}(-(\tilde{t}_p)_S)$
is a subbundle of $E(-(\tilde{t}_p)_S)|_{(\tilde{t}_p)_S}$.
Let $E'|_{(\tilde{t}_p)_S}=(l')^{(p)}_0\supset\cdots\supset(l')^{(p)}_{r-q}$
be the inverse image by $\iota:E'|_{(\tilde{t}_p)_S}\rightarrow E|_{(\tilde{t}_p)_S}$
of the filtration
$l^{(p)}_q\supset\cdots\supset l^{(p)}_r$.
Then each $(l')^{(p)}_j$ is a subbundle of $E'|_{(\tilde{t}_p)_S}$ for $0\leq j\leq r-q$.
Let $(l')^{(p)}_{r-q}\supset\cdots\supset(l')^{(p)}_r$ be the image by
$\pi:E(-(\tilde{t}_p)_S)|_{(\tilde{t}_p)_S}\rightarrow E'|_{(\tilde{t}_p)_S}$
of the filtration
$l^{(p)}_0\otimes{\mathcal O}_{{\mathcal C}_S}(-(\tilde{t}_p)_S)\supset\cdots\supset
l^{(p)}_q\otimes{\mathcal O}_{{\mathcal C}_S}(-(\tilde{t}_p)_S)$.
For each $j$ with $0\leq j\leq q$,
$l^{(p)}_j\otimes{\mathcal O}_{{\mathcal C}_S}(-(\tilde{t}_p)_S)/
l^{(p)}_q\otimes{\mathcal O}_{{\mathcal C}_S}(-(\tilde{t}_p)_S)
\stackrel{\sim}\rightarrow (l')^{(p)}_{r-q+j}$
is a subbundle of
$E(-(\tilde{t}_p)_S)|_{(\tilde{t}_p)_S}/\ker\pi \cong \im\pi=\ker\iota$.
So $(l')^{(p)}_{r-q+j}$ is a subbundle of $E'|_{(\tilde{t}_p)_S}$ for any $j$
with $0\leq j\leq q$.

For $i\neq p$ with $1\leq i\leq n$, we put
$(l')^{(i)}_j:=l^{(i)}_j$ for $0\leq j\leq r$.
Then $(E',\nabla',\{(l')^{(i)}_j\})$ becomes a member of
${\mathcal M}_{{\mathcal C}/T}(\tilde{\bt},r,d-q)(S)$.
We say $(E',\nabla',\{(l')^{(i)}_j\})$ the elementary transform of
$(E,\nabla,\{l^{(i)}_j\})$ along $(\tilde{t}_p)_S$ by $l^{(p)}_q$.
Elementary transform induces a morphism
\begin{gather}\label{elementary-transform}
 \mathrm{Elm}^{(p)}_q:{\mathcal M}_{{\mathcal C}/T}(\tilde{\bt},r,d)
 \longrightarrow {\mathcal M}_{{\mathcal C}/T}(\tilde{\bt},r,d-q) \\
 (E,\nabla,\{l^{(i)}_j\}) \mapsto (E',\nabla',\{(l')^{(i)}_j\})
\end{gather}
of functors.
For any member 
$(E,\nabla,\{l^{(i)}_j\})\in{\mathcal M}_{{\mathcal C}/T}(\tilde{\bt},r,d)(S)$,
we have
\begin{align*}
 \mathrm{Elm}^{(p)}_{r-q}\circ\mathrm{Elm}^{(p)}_q(E,\nabla,\{l^{(i)}_j\})
 &=\mathrm{Elm}^{(p)}_{r-q}(E',\nabla',\{(l')^{(i)}_j\}) \\
 &=(E(-(\tilde{t}_p)_S),\nabla|_{E(-(\tilde{t}_p)_S)}, 
 \{l^{(i)}_j\otimes{\mathcal O}_{{\mathcal C}_S}(-(\tilde{t}_p)_S)\})
\end{align*}
If we put
\begin{gather}\label{twist-line-bundle}
 b_p:{\mathcal M}_{{\mathcal C}/T}(\tilde{\bt},r,d-r)\stackrel{\sim}\longrightarrow
 {\mathcal M}_{{\mathcal C}/T}(\tilde{\bt},r,d) \\
 (E,\nabla,\{l^{(i)}_j\})\mapsto
 (E\otimes{\mathcal O}_{\mathcal C}(\tilde{t}_p),
 \nabla\otimes\mathrm{id}_{{\mathcal O}_{\mathcal C}(\tilde{t}_p)}
 +\mathrm{id}_E\otimes d_{{\mathcal O}_{\mathcal C}(\tilde{t}_p)},
 \{l^{(i)}_j\otimes{\mathcal O}_{\mathcal C}(\tilde{t}_p)\}),
\end{gather}
then we have
$b_p\circ\mathrm{Elm}^{(p)}_{r-q}\circ\mathrm{Elm}^{(p)}_q
=\mathrm{id}_{{\mathcal M}_{{\mathcal C}/T}(\tilde{\bt},r,d)}$.
We can also check that
\begin{align*}
 \mathrm{Elm}^{(p)}_q\circ b_p\circ\mathrm{Elm}^{(p)}_{r-q}
 &=b_p\circ\mathrm{Elm}^{(p)}_q\circ\mathrm{Elm}^{(p)}_{r-q} \\
 &=\mathrm{id}_{{\mathcal M}_{{\mathcal C}/T}(\tilde{\bt},r,d-q)}.
\end{align*}
Thus $\mathrm{Elm}^{(p)}_q$ is an isomorphism:
\[
 \mathrm{Elm}^{(p)}_q:{\mathcal M}_{{\mathcal C}/T}(\tilde{\bt},r,d)
 \stackrel{\sim}\longrightarrow {\mathcal M}_{{\mathcal C}/T}(\tilde{\bt},r,d-q).
\]

Take any $\mathbf{C}$-valued points
$x\in T(\mathbf{C})$, $\blambda\in\Lambda^{(n)}_r(d)(\mathbf{C})$.
We put
\[
 ((\lambda')^{(p)}_0,\ldots,(\lambda')^{(p)}_{r-1})
 =(\lambda^{(p)}_q,\lambda^{(p)}_{q+1},\ldots,\lambda^{(p)}_{r-1},
 \lambda^{(p)}_0+1,\lambda^{(p)}_1+1,\ldots,\lambda^{(p)}_{q-1}+1)
\]
and $(\lambda')^{(i)}_j:=\lambda^{(i)}_j$ for $i\neq p$ with
$1\leq i\leq n$ and $0\leq j\leq r-1$.
Then
$\blambda':=((\lambda')^{(i)}_j)\in\Lambda^{(n)}_r(d-q)(\mathbf{C})$.
For a member
$(E,\nabla,\{l^{(i)}_j\})\in{\mathcal M}_{{\mathcal C}_x}(\tilde{\bt}_x,\blambda)(S)$,
$\mathrm{Elm}^{(p)}_q(E,\nabla,\{l^{(i)}_j\})\in
{\mathcal M}_{{\mathcal C}_x}(\tilde{\bt}_x,\blambda')(S)$.
So $\mathrm{Elm}^{(p)}_q$ induces an isomorphism
\begin{equation}\label{e-transform}
 \mathrm{Elm}^{(p)}_q:{\mathcal M}_{{\mathcal C}_x}(\tilde{\bt}_x,\blambda)
 \stackrel{\sim}\longrightarrow
 {\mathcal M}_{{\mathcal C}_x}(\tilde{\bt}_x,\blambda').
\end{equation}

Take a point $\blambda\in\Lambda^{(n)}_r(d)(\mathbf{C})$.
We can write
\[
 \prod_{j=0}^{r-1}(x-\lambda^{(p)}_j)=
 \prod_{k=1}^u(x-\mu^{(p)}_k)^{m^{(p)}_k},
\]
where
$\mathrm{Re}(\mu^{(p)}_1)\geq \mathrm{Re}(\mu^{(p)}_2)
\geq\cdots\geq\mathrm{Re}(\mu^{(p)}_u)$,
$\mathrm{Im}(\mu^{(p)}_j)>\mathrm{Im}(\mu^{(p)}_{j+1})$ if
$\mathrm{Re}(\mu^{(p)}_j)=\mathrm{Re}(\mu^{(p)}_{j+1})$
and $x$ is an indeterminate.
For any member
$(E,\nabla_E,\{l^{(i)}_j\})\in {\mathcal M}_{{\mathcal C}_x}(\tilde{\bt}_x,\blambda)(S)$,
we can find $j^{(k)}_0,j^{(k)}_1,\ldots,j^{(k)}_{m^{(p)}_k}$ for $1\leq k\leq u$
satisfying
$0\leq j^{(k)}_0<j^{(k)}_1<\cdots<j^{(k)}_{m^{(p)}_k}\leq r$
such that
\[
 \left(
 \ker\left(\res_{(\tilde{t}_p)_x}(\nabla_E)-\mu^{(p)}_k\cdot\mathrm{id}\right)^r
 \cap l^{(p)}_{j^{(k)}_m}\right)\left/
 \left(
 \ker\left(\res_{(\tilde{t}_p)_x}(\nabla_E)-\mu^{(p)}_k\cdot\mathrm{id}\right)^r
 \cap l^{(p)}_{j^{(k)}_{m+1}}\right)\right.
 \stackrel{\sim}\longrightarrow
 l^{(p)}_{j^{(k)}_m}/l^{(p)}_{j^{(k)}_m+1}
\]
for $m=0,\ldots,m^{(p)}_k-1$.
Let $E|_{(\tilde{t}_p)_x}=(l')^{(p)}_0\supset(l')^{(p)}_1\supset\cdots
\supset(l')^{(p)}_{r-1}\supset(l')^{(p)}_r=0$
be the filtration by subbundles satisfying
$\rank (l')^{(p)}_j/(l')^{(p)}_{j+1}=1$
for $j=0,\ldots,r-1$ and
\[
 (l')^{(p)}_q=
 \left(\ker\left(\res_{(\tilde{t}_p)_x}(\nabla_E)-\mu^{(p)}_k\cdot\mathrm{id}\right)^r
 \cap l^{(p)}_{j^{(k)}_m}\right)
 +(l')^{(p)}_{q+1}
\]
for $q=0,\ldots,r-1$, where
$q=r-(m^{(p)}_k-m)-\sum_{k<k'\leq u} m^{(p)}_{k'}$,
$1\leq k\leq u$ and $0\leq m<m^{(p)}_k$.
We put 
\[
 (\nu^{(p)}_0,\ldots,\nu^{(p)}_{r-1}):=
 (\overbrace{\mu^{(p)}_1,\ldots,\mu^{(p)}_1}^{m^{(p)}_1},\ldots,
 \overbrace{\mu^{(p)}_k,\ldots,\mu^{(p)}_k}^{m^{(p)}_k},\ldots,
 \overbrace{\mu^{(p)}_u,\ldots,\mu^{(p)}_u}^{m^{(p)}_u}),
\]
$\nu^{(i)}_j=\lambda^{(i)}_j$ for $i\neq p$
and $\bnu:=(\nu^{(i)}_j)_{i,j}$.
Then the correspondence
$(E,\nabla_E,\{l^{(i)}_j\})\mapsto
(E,\nabla_E,\{(l')^{(i)}_j\})$
determines an isomorphism
\begin{equation}\label{permutation}
 a_p:\cM_{{\mathcal C}_x}(\tilde{\bt}_x,\blambda)\stackrel{\sim}\lra
 \cM_{{\mathcal C}_x}(\tilde{\bt}_x,\bnu)
\end{equation}
of functors.

Looking at the change of $\blambda$ by the transformations
$\mathrm{Elm}^{(i)}_j$ and $b_i$,
we can easily see the following proposition,
which is useful in considering the Riemann-Hilbert correspondence
in section \ref{riemann-hilbert}.

\begin{Proposition}\label{prop:elementary-transform}
 Composing the isomorphisms $a_i,b_i$ and $\mathrm{Elm}^{(i)}_j$
 for $i=1,\ldots,n$, $j=1,\ldots,r$ given in (\ref{permutation}),
 (\ref{e-transform}) and (\ref{twist-line-bundle}),
 we obtain an isomorphism
 \begin{equation}
  \sigma: \cM_{{\mathcal C}_x}(\tilde{\bt}_x,\blambda)\stackrel{\sim}\lra
  \cM_{{\mathcal C}_x}(\tilde{\bt}_x,\bmu)
 \end{equation}
 of functors,
 where $0\leq \mathrm{Re}(\mu^{(i)}_j)<1$
 for any $i,j$.
\end{Proposition}

As a corollary, we obtain the following:

\begin{Corollary}
 Take $\blambda\in\Lambda^{(n)}_r(d)$ and $\bmu\in\Lambda^{(n)}_r(d')$
 which satisfy $rh(\blambda)=rh(\bmu)\in\cA^{(n)}_r$,
 where $rh$ is the morphism defined in (\ref{local-riemann-hilbert}).
 Then we can obtain by composing $a_i,a_i^{-1},b_i,b_i^{-1}$
 and $\mathrm{Elm}^{(i)}_j$
 an isomorphism
 \[
  \sigma:\cM_{{\mathcal C}_x}(\tilde{\bt}_x,\blambda)\stackrel{\sim}\lra
 \cM_{{\mathcal C}_x}(\tilde{\bt}_x,\bmu)
 \]
 of functors.
\end{Corollary}

\begin{proof}
Applying Proposition \ref{prop:elementary-transform},
we obtain an isomorphism
\[
 \sigma: \cM_{{\mathcal C}_x}(\tilde{\bt}_x,\blambda)\stackrel{\sim}\lra
 \cM_{{\mathcal C}_x}(\tilde{\bt}_x,\bnu)
\]
of functors, where
$0\leq\mathrm{Re}(\nu^{(i)}_j)<1$ for any $i,j$.
By composing $a_i$, we may also assume that
$\mathrm{Re}(\nu^{(i)}_0)\geq\mathrm{Re}(\nu^{(i)}_1)
\geq\cdots\geq\mathrm{Re}(\nu^{(i)}_{r-1})$ for any $i$ and
$\mathrm{Im}(\nu^{(i)}_j)\geq \mathrm{Im}(\nu^{(i)}_{j+1})$
if $\mathrm{Re}(\nu^{(i)}_j)=\mathrm{Re}(\nu^{(i)}_{j+1})$.
Applying Proposition \ref{prop:elementary-transform} again,
we have an isomorphism
\[
 \sigma': \cM_{{\mathcal C}_x}(\tilde{\bt}_x,\bmu)\stackrel{\sim}\lra
 \cM_{{\mathcal C}_x}(\tilde{\bt}_x,\bnu).
\]
Then
\[
 (\sigma')^{-1}\circ\sigma:\cM_{{\mathcal C}_x}(\tilde{\bt}_x,\blambda)
 \stackrel{\sim}\lra \cM_{{\mathcal C}_x}(\tilde{\bt}_x,\bmu)
\]
is a desired isomorphism.
\end{proof}

\section{Moduli of stable parabolic connections}\label{moduli of s-p-c}

Before proving Theorem \ref{moduli-exists},
we recall the definition of parabolic $\Lambda^1_{D}$-triple
defined in \cite{IIS-1}.
Let $D$ be an effective divisor on a curve $C$.
We define $\Lambda^1_{D}$ as $\cO_C\oplus\Omega^1_C(D)^{\vee}$
with the bimodule structure given by
\begin{align*}
 f(a,v)&=(fa,fv) \quad (f,a\in\cO_C, v\in\Omega^1_C(D)^{\vee}) \\
 (a,v)f&=(fa+v(f),fv)\quad (f,a\in\cO_C, v\in\Omega^1_C(D)^{\vee}). 
\end{align*}

\begin{Definition}\rm
We say $(E_1,E_2,\Phi,F_*(E_1))$ a parabolic
$\Lambda^1_{D}$-triple on $C$ of rank $r$ and degree $d$ if
\begin{enumerate}
\item $E_1$ and $E_2$ are vector bundles on $C$ of rank $r$ and degree $d$,
\item $\Phi:\Lambda^1_{D}\otimes E_1 \ra E_2$
 is a left $\cO_C$-homomorphism,
\item $E_1=F_1(E_1)\supset F_2(E_1)\supset\cdots\supset F_l(E_1)\supset
 F_{l+1}(E_1)=E_1(-D)$
 is a filtration by coherent subsheaves.
\end{enumerate}
\end{Definition}

Note that to give a left ${\mathcal O}_C$-homomorphism
$\Phi:\Lambda^1_D\otimes E_1\rightarrow E_2$
is equivalent to give an ${\mathcal O}_C$-homomorphism
$\phi:E_1\rightarrow E_2$ and a morphism
$\nabla:E_1\rightarrow E_2\otimes\Omega^1_C(D)$
such that
$\nabla(fa)=\phi(a)\otimes df+f\nabla(a)$
for $f\in{\mathcal O}_C$ and $a\in E_1$.
We also denote the parabolic $\Lambda^1_D$-triple
$(E_1,E_2,\Phi,F_*(E_1))$ by
$(E_1,E_2,\phi,\nabla,F_*(E_1))$.

We take positive integers
$\beta_1,\beta_2,\gamma$ and rational numbers
$0<\alpha'_1<\cdots<\alpha'_l<1$.
We assume $\gamma\gg 0$.

\begin{Definition}\rm
A parabolic $\Lambda^1_{D}$-triple $(E_1,E_2,\Phi,F_*(E_1))$
is $(\balpha,\bbeta,\gamma)$-stable
(resp.\ $(\balpha,\bbeta,\gamma)$-semistable)
if for any subbundles $(F_1,F_2)\subset(E_1,E_2)$ satisfying
$(0,0)\neq(F_1,F_2)\neq(E_1,E_2)$ and
$\Phi(\Lambda^1_{D(\bt)}\otimes F_1)\subset F_2$, the inequality
\begin{gather*}
 \frac{\beta_1\deg F_1(-D)+\beta_2(\deg F_2-\gamma\rank F_2)
 +\beta_1\sum_{j=1}^l
 \alpha'_j\length((F_j(E_1)\cap F_1)/(F_{j+1}(E_1)\cap F_1))}
 {\beta_1\rank F_1+\beta_2\rank F_2} \\
 \genfrac{}{}{0pt}{}{<}{\text{(resp.\  $\leq$)}}
 \frac{\beta_1\deg E_1(-D)+\beta_2(\deg E_2-\gamma\rank E_2)
 +\beta_1\sum_{j=1}^l
 \alpha'_j\length(F_j(E_1)/F_{j+1}(E_1))}
 {\beta_1\rank E_1+\beta_2\rank E_2}
\end{gather*}
holds.
\end{Definition}

\begin{Theorem}\label{moduli-p-triple}{\bf (\cite{IIS-1}, Theorem 7.1)}
Let $S$ be an algebraic scheme over $\mathbf{C}$, $\cC$ be a flat family of
smooth projective curves of genus $g$ and $D$ be an effective Cartier divisor
on $\cC$ flat over $S$.
Then there exists a coarse moduli scheme
$\overline{M^{D,\balpha',\bbeta,\gamma}_{\cC/S}}
(r,d,\{d_i\}_{1\leq i\leq nr})$
of $(\balpha',\bbeta,\gamma)$-stable
parabolic $\Lambda^1_{D}$-triples $(E_1,E_2,\Phi,F_*(E_1))$
on $\cC$ over $S$ such that
$r=\rank E_1=\rank E_2$,
$d=\deg E_1=\deg E_2$ and $d_i=\length(E_1/F_{i+1}(E_1))$.
If $\balpha'$ is generic, it is projective over $S$.
\end{Theorem}

\begin{Definition}\label{def-moduli-functor}\rm
We denote the pull-back of $\cC$ and $\tilde{\bt}$ by the morphism
$T\times\Lambda^{(n)}_r(d)\rightarrow T$
by the same characters $\cC$ and
$\tilde{\bt}=(\tilde{t}_1,\ldots,\tilde{t}_n)$.
Then $D(\tilde{\bt}):=\tilde{t}_1+\cdots+\tilde{t}_n$
becomes an effective Cartier divisor on $\cC$ flat over
$T\times\Lambda^{(n)}_r(d)$.
We also denote by $\tilde{\blambda}$ the pull-back of the
universal family on $\Lambda^{(n)}_r(d)$ by the morphism
$T\times\Lambda^{(n)}_r(d)\rightarrow \Lambda^{(n)}_r(d)$.
We define a functor
$\cM^{\balpha}_{\cC/T}
(\tilde{\bt},r,d)$
of the category of locally noetherian schemes to the category of sets by
\[
 \cM^{\balpha}_{\cC/T}
 (\tilde{\bt},r,d)(S):=
 \left\{ (E,\nabla,\{l^{(i)}_j\}) \right\}/\sim,
\]
for a locally noetherian scheme $S$ over
$T\times\Lambda^{(n)}_r(d)$, where
\begin{enumerate}
\item $E$ is a vector bundle on $\cC_S$ of rank $r$,
\item $\nabla:E\rightarrow E\otimes\Omega^1_{\cC_S/S}(D(\tilde{\bt})_S)$
 is a relative connection,
\item $E|_{(\tilde{t}_i)_S}=l^{(i)}_0\supset l^{(i)}_1
 \supset\cdots\supset l^{(i)}_{r-1}\supset l^{(i)}_r=0$
 is a filtration by subbundles such that
 $(\res_{(\tilde{t}_i)_S}(\nabla)-(\tilde{\lambda}^{(i)}_j)_S)(l^{(i)}_j)
 \subset l^{(i)}_{j+1}$
 for $0\leq j\leq r-1$, $i=1,\ldots,n$,
\item for any geometric point $s\in S$,
$\dim (l^{(i)}_j/l^{(i)}_{j+1})\otimes k(s)=1$ for any $i,j$ and
$(E,\nabla,\{l^{(i)}_j\})\otimes k(s)$
is $\balpha$-stable.
\end{enumerate}
Here $(E,\nabla,\{l^{(i)}_j\})\sim
(E',\nabla',\{l'^{(i)}_j\})$ if
there exist a line bundle $\cL$ on $S$ and
an isomorphism $\sigma:E\stackrel{\sim}\ra E'\otimes\cL$ 
such that $\sigma|_{t_i}(l^{(i)}_j)=l'^{(i)}_j\otimes{\mathcal L}$
for any $i,j$ and the diagram
\[
 \begin{CD}
  E @>\nabla>> E\otimes\Omega^1_{\cC/T}(D(\tilde{\bt})) \\
  @V \sigma VV  @V\sigma\otimes\mathrm{id} VV \\
  E'\otimes\cL @>\nabla'\otimes{\mathcal L}>>
  E'\otimes\Omega^1_{\cC/T}(D(\tilde{\bt}))\otimes\cL \\
 \end{CD}
\]
commutes.
\end{Definition}

\noindent
{\bf Proof of Theorem \ref{moduli-exists}.}
Fix a weight $\balpha$ which determines the stability of
parabolic connections.
We take positive integers
$\beta_1,\beta_2,\gamma$ and rational numbers
$0<\tilde{\alpha}^{(i)}_1<\tilde{\alpha}^{(i)}_2<
\cdots<\tilde{\alpha}^{(i)}_r<1$
satisfying
$(\beta_1+\beta_2)\alpha^{(i)}_j=\beta_1\tilde{\alpha}^{(i)}_j$
for any $i,j$.
We assume $\gamma\gg 0$.
We can take an increasing sequence
$0<\alpha'_1<\cdots<\alpha'_{nr}<1$ such that
$\{\alpha'_i|1\leq i\leq nr \}=
\left.\left\{\tilde{\alpha}^{(i)}_j \right|
1\leq i\leq n, 1\leq j\leq r\right\}$.
Take any member $(E,\nabla,\{l^{(i)}_j\})\in
\cM^{\balpha}_{\cC/T}(\tilde{\bt},r,d)(S)$.
For each $1\leq p \leq rn$, exist $i,j$ satisfying
$\tilde{\alpha}^{(i)}_j=\alpha'_p$.
We put $F_1(E):=E$ and define inductively
$F_{p+1}(E):=\ker(F_{p}(E)\rightarrow E|_{(\tilde{t}_i)_S}/l^{(i)}_j)$
for $p=1,\ldots,rn$.
We also put
$d_p:=\length((E/F_{p+1}(E))\otimes k(s))$ for $p=1,\ldots,rn$
and $s\in S$.
Then $(E,\nabla,\{l^{(i)}_j\})\mapsto
(E,E,\mathrm{id}_E,\nabla,F_*(E))$ determines a morphism
\[
 \iota:\cM^{\balpha}_{\cC/T}(\tilde{\bt},r,d)\rightarrow
 \overline{\cM^{D(\tilde{\bt}),\balpha',\bbeta,\gamma}
 _{\cC/T\times\Lambda^{(n)}_r(d)}}
 (r,d,\{d_i\}_{1\leq i\leq nr}),
\]
where $\overline{\cM^{D(\tilde{\bt}),\balpha',\bbeta,\gamma}_{\cC/T}}
(r,d,\{d_i\}_{1\leq i\leq nr})$
is the moduli functor of $(\balpha',\bbeta,\gamma)$-stable parabolic
$\Lambda^1_{D(\tilde{\bt})}$-triples
whose coarse moduli scheme exists by Theorem \ref{moduli-p-triple}.
Note that a parabolic connection $(E,\nabla,\{l^{(i)}_j\})$
is $\balpha$-stable if and only if the corresponding parabolic
$\Lambda^1_{D(\bt)}$-triple $(E,E,\mathrm{id},\nabla,F_*(E))$
is $(\balpha',\bbeta,\gamma)$-stable since $\gamma\gg 0$.
We can see that $\iota$ is representable by an immersion.
So we can prove in the same way as [\cite{IIS-1} Theorem 2.1]
that a certain subscheme $M^{\balpha}_{\cC/T}(\tilde{\bt},r,d)$
of $\overline{M^{D(\bt),\balpha',\bbeta,\gamma}_C}
(r,d,\{d_i\}_{1\leq i\leq nr})$
is just the coarse moduli scheme of
$\cM^{\balpha}_{\cC/T}(\tilde{\bt},r,d)$.

 Applying a certain elementary transformation, we obtain an isomorphism
\[
 \sigma:\cM_{\cC/T}(\tilde{\bt},r,d)\stackrel{\sim}\lra
 \cM_{\cC/T}(\tilde{\bt},r,d')
\]
of functors,
where $d'$ and $r$ are coprime.
Then $\sigma(M^{\balpha}_{\cC/T}(\tilde{\bt},r,d))$ can be
considered as the moduli scheme of parabolic connections of rank $r$
and degree $d'$ satisfying a certain stability condition.
We can take a vector space $V$ such that
there exists a surjection $V\otimes{\mathcal O}_{\cC_s}(-m)\rightarrow E$
such that $V\otimes k(s)\rightarrow H^0(E(m))$ is isomorphic
and $h^i(E(m))=0$ for $i>0$ for any
member $(E,\nabla,\{l^{(i)}_j\})\in
\sigma(M^{\balpha}_{\cC/T}(\tilde{\bt},r,d))$
over $s\in T$.
We put $P(n):=\chi(E(n))$.
Let ${\mathcal E}$ be the universal family on 
$\cC\times_{T\times\Lambda^{(n)}_r(d)}
\Quot^P_{V\otimes{\mathcal O}_{\mathcal C}(-m)/\cC/T\times\Lambda^{(n)}_r(d)}$.
Then we can construct a scheme $R$ over
$\Quot^P_{V\otimes{\mathcal O}_{\mathcal C}(-m)/\cC/T\times\Lambda^{(n)}_r(d)}$
which parametrizes connections
$\nabla:{\mathcal E}_s\rightarrow
{\mathcal E}_s\otimes\Omega^1_{\cC_s}(D(\tilde{\bt}_s))$
and parabolic structures
${\mathcal E}_s|_{(\tilde{t}_i)_s}=l^{(i)}_0\supset\cdots\supset
l^{(i)}_{r-1}\supset l^{(i)}_r=0$
such that
$(\res_{(\tilde{t}_i)_s}(\nabla)-(\tilde{\lambda}^{(i)}_j)_s)
(l^{(i)}_j)\subset l^{(i)}_{j+1}$ for any $i,j$.
Then $PGL(V)$ canonically acts on $R$ and for some open subscheme
$R^s$ of $R$, there is a canonical morphism
$R^s\rightarrow \sigma(M^{\balpha}_{\cC/T}(\tilde{\bt},r,d))$
which becomes a principal $PGL(V)$-bundle.
So we can prove in the same manner as [\cite{HL}, Theorem 4.6.5]
that for a certain line bundle ${\mathcal L}$ on $R^s$,
${\mathcal E}\otimes{\mathcal L}$ descends to a vector bundle on
$\cC\times_{T\times\Lambda^{(n)}_r(d)}
\sigma(M^{\balpha}_{\cC/T}(\tilde{\bt},r,d))$,
since $r$ and $d'$ are coprime.
We can easily see that the universal families
$\tilde{\nabla}$, $\{\tilde{l}^{(i)}_j\}$
of connections and parabolic structures on
${\mathcal E}\otimes{\mathcal L}$
also descend.
So we obtain a universal family for the moduli space
$\sigma(M^{\balpha}_{\cC/T}(\tilde{\bt},r,d))$
and $M^{\balpha}_{\cC/T}(\tilde{\bt},r,d)$
in fact becomes a fine moduli scheme.

Let $(\tilde{E},\tilde{\nabla},\{\tilde{l}^{(i)}_j\})$
be a universal family on
$\cC\times_T M^{\balpha}_{\cC/T}(\tilde{\bt},r,d)$.
We define a complex $\cF^{\bullet}$ by
\begin{align*}
 \cF^0&:=\left\{ s\in {\mathcal End}(\tilde{E}) \left|
 \text{$s|_{\tilde{t}_i\times M^{\balpha}_{\cC/T}(\tilde{\bt},r,d)}
 (\tilde{l}^{(i)}_j)\subset\tilde{l}^{(i)}_j$ for any $i,j$}
 \right\}\right. \\
 \cF^1&:=\left\{ s\in {\mathcal End}(\tilde{E})
 \otimes\Omega^1_{\cC/T}(D(\tilde{\bt})) \left|
 \text{$\res_{\tilde{t}_i\times M^{\balpha}_{\cC/T}(\tilde{\bt},r,d)}(s)
 (\tilde{l}^{(i)}_j)\subset \tilde{l}^{(i)}_{j+1}$ for any $i,j$}
 \right\}\right. \\
 \nabla_{\cF^{\bullet}} &: \cF^0 \lra \cF^1 ; \quad
 \nabla_{\cF^{\bullet}}(s)=\tilde{\nabla}\circ s - s\circ\tilde{\nabla}.
\end{align*}
Let $(A,m)$ be an artinian local ring over $T\times\Lambda^{(n)}_r(d)$,
$I$ an ideal of $A$ satisfying $mI=0$.
Take any member $(E,\nabla,\{l^{(i)}_j\})\in
M^{\balpha}_{\cC/T}(\tilde{\bt},r,d)(A/I)$.
We denote the image of
$\Spec A/m\rightarrow
M^{\balpha}_{\cC/T}(\tilde{\bt},r,d)$
by $x$.
The structure morphism
$\Spec A\rightarrow T\times\Lambda^{(n)}_r(d)$
corresponds to an $A$-valued point
$(\bt,\blambda)\in T(A)\times\left(\Lambda^{(n)}_r(d)\right)(A)$.
Let $\cC\times \Spec A=\bigcup_{\alpha}U_{\alpha}$
be an affine open covering
such that
$E|_{U_{\alpha}\otimes A/I}\cong
{\mathcal O}_{U_{\alpha}\otimes A/I}^{\oplus r}$,
$\sharp\{i|t_i\in U_{\alpha}\}\leq 1$ for any $\alpha$ and
$\sharp\{\alpha|t_i\in U_{\alpha}\}\leq 1$ for any $i$.
We can obviously take a vector bundle
$E_{\alpha}$ on $U_{\alpha}$
such that there is an isomorphism
$E_{\alpha}\otimes A/I\xrightarrow[\sim]{\phi_{\alpha}}
E|_{U_{\alpha}\otimes A/I}$.
If $t_i\in U_{\alpha}$, we take a basis
$e_1,e_2,\ldots,e_r$ of $E_{\alpha}$
such that
$l^{(i)}_{r-1}= \langle e_1|_{t_i\otimes A/I} \rangle,
l^{(i)}_{r-2}=\langle e_1|_{t_i\otimes A/I},e_2|_{t_i\otimes A/I} \rangle,\ldots,
l^{(i)}_1=\langle e_1|_{t_i\otimes A/I},\ldots,e_{r-1}|_{t_i\otimes A/I}\rangle$.
The connection matrix of $\nabla|_{E|_{U_{\alpha}\otimes A/I}}$
with respect to the basis $e_1\otimes A/I,\ldots,e_r\otimes A/I$
can be given by a matrix
$\Omega_{\alpha}\in M_r(\Omega^1_{\cC_{A/I}/(A/I)}((t_i)_{A/I})|
_{U_{\alpha}\otimes A/I})$
such that
\[
 \Omega_{\alpha}|_{t_i\otimes A/I}=
 \begin{pmatrix}
  \lambda^{(i)}_{r-1}\otimes A/I & \omega_{12}  & \cdots & \omega_{1r} \\
  0 & \lambda^{(i)}_{r-2}\otimes A/I & \cdots & \omega_{2r} \\
  \vdots & \vdots & \ddots & \vdots \\
  0 & 0 & \cdots & \lambda^{(i)}_0\otimes A/I
 \end{pmatrix}.
\]
We can take a lift
$\tilde{\Omega}_{\alpha}\in M_r(\Omega^1_{\cC_A/A}(t_i)|_{U_{\alpha}})$
of $\Omega_{\alpha}$ such that
\[
 \tilde{\Omega}_{\alpha}|_{t_i}:=
 \begin{pmatrix}
  \lambda^{(i)}_{r-1} & \tilde{\omega}_{12} & \cdots & \tilde{\omega}_{1r} \\
  0 & \lambda^{(i)}_{r-2} & \cdots & \tilde{\omega}_{2r} \\
  \vdots & \vdots & \ddots & \vdots \\
  0 & 0 & \cdots & \lambda^{(i)}_0
 \end{pmatrix},
\]
which induces a connection
$\nabla_{\alpha}:E_{\alpha}\rightarrow
E_{\alpha}\otimes\Omega^1_{\cC_A/A}((t_1+t_2+\cdots+t_n)_A)$.
We define a parabolic structure $\{(l_{\alpha})^{(i)}_j\}$ on
$E_{\alpha}$ by
$(l_{\alpha})^{(i)}_{r-1}:=\langle e_1|_{t_i} \rangle$,
$(l_{\alpha})^{(i)}_{r-2}:=\langle e_1|_{t_i},e_2|_{t_i}\rangle,
\ldots,(l_{\alpha})^{(i)}_1:=\langle e_1|_{t_i},\ldots,e_{r-1}|_{t_i}\rangle$.
Then
$(E_{\alpha},\nabla_{\alpha},\{(l_{\alpha})^{(i)}_j\})$
become a lift $(E,\nabla,\{l^{(i)}_j\})|_{U_{\alpha}}$.
If $t_i\notin U_{\alpha}$ we can easily take a lift
$(E_{\alpha},\nabla_{\alpha},\{(l_{\alpha})^{(i)}_j\})$
of  $(E,\nabla,\{l^{(i)}_j\})|_{U_{\alpha}}$.
We put $U_{\alpha\beta}:=U_{\alpha}\cap U_{\beta}$ and
$U_{\alpha\beta\gamma}:=U_{\alpha}\cap U_{\beta}\cap U_{\gamma}$.
Take a lift
$\theta_{\beta\alpha}:E_{\alpha}|_{U_{\alpha\beta}}
\stackrel{\sim}\rightarrow E_{\beta}|_{U_{\alpha\beta}}$
of the canonical isomorphism
$E_{\alpha}|_{U_{\alpha\beta}}\otimes A/I\xrightarrow[\sim]{\phi_{\alpha}}
E|_{U_{\alpha\beta}\otimes A/I}\xrightarrow[\sim]{\phi_{\beta}^{-1}}
E_{\beta}|_{U_{\alpha\beta}}\otimes A/I$.
We put
\begin{align*}
 u_{\alpha\beta\gamma}:&=\phi_{\alpha}|_{U_{\alpha\beta\gamma}}\circ
 \left(\theta_{\gamma\alpha}^{-1}|_{U_{\alpha\beta\gamma}}\circ
 \theta_{\gamma\beta}|_{U_{\alpha\beta\gamma}}\circ
 \theta_{\beta\alpha}|_{U_{\alpha\beta\gamma}}
 -\mathrm{id}_{E_{\alpha}|_{U_{\alpha\beta\gamma}}}\right)
 \circ\phi_{\alpha}^{-1}|_{U_{\alpha\beta\gamma}}, \\
 v_{\alpha\beta}:&=\phi_{\alpha}\circ
 \left(\nabla_{\alpha}|_{U_{\alpha\beta}}-
 \theta_{\beta\alpha}^{-1}\circ\nabla_{\beta}\circ\theta_{\beta\alpha}
 \right)\circ\phi_{\alpha}^{-1}|_{U_{\alpha\beta\gamma}}.
\end{align*}
Then we have
$\{u_{\alpha\beta\gamma}\}\in
C^2(\{U_{\alpha}\},{\mathcal F}^0\otimes I)$
and
$\{v_{\alpha\beta}\}\in
C^1(\{U_{\alpha}\},{\mathcal F}^1\otimes I)$.
We can easily see that
\[
 d\{v_{\alpha\beta}\}=
 -\nabla_{{\mathcal F}^{\bullet}}\{u_{\alpha\beta\gamma}\}
 \quad \text{and} \quad
 d\{u_{\alpha\beta\gamma}\}=0.
\]
So we can define an element
\[
 \omega(E,\nabla,\{l^{(i)}_j\}):=
 [(\{u_{\alpha\beta\gamma}\},\{v_{\alpha\beta}\})]
 \in \check{\mathbf{H}}^2(\{U_{\alpha}\},{\mathcal F}^{\bullet}\otimes I)
 =\mathbf{H}^2({\mathcal F}^{\bullet}\otimes k(x))\otimes I.
\]
We can easily see that
$(E,\nabla,\{l^{(i)}_j\})$
can be lifted to an $A$-valued point of
$M^{\balpha}_{\cC/T}(\tilde{\bt},r,d)$
if and only if $\omega(E,\nabla,\{l^{(i)}_j\})=0$
in the hypercohomology
$\bH^2(\cF^{\bullet}\otimes k(x))\otimes I$.

Consider the homomorphism
\begin{gather*}
 \Tr: \bH^2(\cF^{\bullet}\otimes k(x))\otimes I \longrightarrow
 \bH^2(\Omega^{\bullet}_{\cC_x})\otimes I \hspace{20pt} \\
 \hspace{50pt}
 [\{u'_{\alpha\beta\gamma}\},\{v'_{\alpha\beta}\}]\mapsto
 [\{\Tr(u'_{\alpha\beta\gamma})\},\{\Tr(v'_{\alpha\beta})\}],
\end{gather*}
Then we can see that $\Tr(\omega(E,\nabla,\{l^{(i)}_j\}))$
is just the obstruction class for lifting $(\det E,\wedge^r\nabla)$
to a pair of a line bundle and a connection with the fixed residue
$\det(\tilde{\blambda}_A)$ over $A$.
Here $\det(\tilde{\blambda}_A)=
(\sum_{j=0}^{r-1}(\tilde{\lambda}^{(i)}_j)_A)^{1\leq i\leq n}$.
Since $\Pic^d_{\cC/T}$ is smooth over $T$, there is a line bundle
$L$ on $\cC_A$ such that $L\otimes A/I\cong \det(E)$.
We can construct a connection
$\tilde{\nabla}^0:L\rightarrow
L\otimes\Omega^1_{\cC_A/A}(\tilde{t}_1+\cdots+\tilde{t}_n)$
because
$H^1({\mathcal End}(L)\otimes
\Omega^1_{\cC_A/A}((\tilde{t}_1+\cdots+\tilde{t}_n)_A))=0$.
We can find a member
$u\in H^0({\mathcal End}(L)\otimes
\Omega^1_{\cC_A/A}((\tilde{t}_1+\cdots+\tilde{t}_n)_A))$
such that $u\otimes A/I=\tilde{\nabla}^0\otimes A/I-\wedge^r\nabla$.
Then $(L,\tilde{\nabla}^0-u)$ becomes a lift of
$(\det(E),\wedge^r\nabla)$.
If we put
$\mu^{(i)}:=\res_{\tilde{t}_i\otimes A}(\tilde{\nabla}^0-u)
-\sum_j(\tilde{\lambda}^{(i)}_j)_A$,
then we have $\mu^{(i)}\in I$ for any $i$ and
$\sum_{i=1}^n\mu^{(i)}=0$.
Note that there is an exact sequence
\begin{gather*}
 H^0(I\otimes{\mathcal End}(L)\otimes\Omega^1_{\cC_A/A}
 ((\tilde{t}_1+\cdots+\tilde{t}_n)_A) \longrightarrow
 \bigoplus_{i=1}^nI\otimes{\mathcal End}(L)\otimes\Omega^1_{\cC_A/A}
 ((\tilde{t}_1+\cdots+\tilde{t}_n)_A)|_{(\tilde{t}_i)_A}
 \cong\bigoplus_{i=1}^n I  \\
 \stackrel{f}\longrightarrow I\cong
 H^1(I\otimes{\mathcal End}(L)\otimes\Omega^1_{\cC_A/A})
 \longrightarrow 0,
\end{gather*}
where $f$ is given by $f((\nu_i)_{1\leq i\leq n})=\sum_{i=1}^n\nu_i$.
Then we can find an element
$\omega\in I\otimes
H^0({\mathcal End}(L)\otimes
\Omega^1_{\cC_A/A}((\tilde{t}_1+\cdots+\tilde{t}_n)_A))$
such that
$\res_{\tilde{t}_i\otimes A}(\omega)=\mu^{(i)}$ for any $i$.
Then $(L,\tilde{\nabla}^0-u-\omega)$ is a lift of
$(\det(E),\wedge^r\nabla)$ with the residue
$(\det(\tilde{\blambda})_A)$.
Thus we have $\Tr(\omega(E,\nabla,\{l^{(i)}_j\}))=0$.
Since $(\cF^1)^{\vee}\otimes\Omega^1_{\cC_A/A}\cong\cF^0$ and
$(\cF^0)^{\vee}\otimes\Omega^1_{\cC_A/A}\cong\cF^1$, we have
\begin{align*}
 \bH^2(\cF^{\bullet}\otimes k(x)) &\cong
 \coker\left(H^1(\cF^0\otimes k(x))\ra H^1(\cF^1\otimes k(x))\right) \\
 &\cong
 \ker\left(H^0((\cF^1)^{\vee}\otimes\Omega^1_{\cC/T}\otimes k(x))
 \ra H^0((\cF^0)^{\vee}\otimes\Omega^1_{\cC/T}\otimes k(x))
 \right)^{\vee} \\
 &\cong \ker\left(H^0(\cF^0\otimes k(x))\ra
 H^0(\cF^1\otimes k(x))\right)^{\vee}.
\end{align*}
and a commutative diagram
\[
 \begin{CD}
  \bH^2(\cF^{\bullet}\otimes k(x)) @>\Tr>> \bH^2(\Omega_{\cC_x}^{\bullet}) \\
  @VV\cong V   @VV\cong V \\
  \ker\left(H^0(\cF^0\otimes k(x))\stackrel{\nabla_{\cF^{\bullet}}}
  \longrightarrow H^0(\cF^1\otimes k(x))\right)^{\vee}
  @>\iota^{\vee}>>
  \ker\left(H^0(\cO_{\cC_x})\stackrel{d}\rightarrow
  H^0(\Omega^1_{\cC_x})\right)^{\vee}.
 \end{CD}
\]
We can see by the stability that the endomorphisms of
$(\tilde{E},\tilde{\nabla},\{\tilde{l}^{(i)}_j\})\otimes k(x)$
are scalar multiplications.
So we have
$\ker\left(H^0(\cF^0\otimes k(x))\ra H^0(\cF^1\otimes k(x))\right)=k(x)$
and the canonical inclusion
\[
 \iota:\ker\left(H^0(\cO_{\cC_x})\stackrel{d}\rightarrow
 H^0(\Omega^1_{\cC_x})\right)
 \longrightarrow
 \ker\left(H^0(\cF^0\otimes k(x))\stackrel{\nabla_{\cF^{\bullet}}}
 \longrightarrow H^0(\cF^1\otimes k(x))\right)
\]
is an isomorphism.
Hence $\Tr$ is an isomorphism and
$\omega(E,\nabla,\{l^{(i)}_j\})=0$
because $\Tr(\omega(E,\nabla,\{l^{(i)}_j\}))=0$.
So we have proved that
$M^{\balpha}_{\cC/T}(\tilde{\bt},r,d)$
is smooth over $T\times\Lambda^{(n)}_r(d)$.

Take an affine open set $M\subset M^{\balpha}_{\cC/T}(\tilde{\bt},r,d)$
and an affine open covering $\cC_M=\bigcup_{\alpha}U_{\alpha}$
such that $\tilde{E}|_{U_{\alpha}}\cong{\mathcal O}_{U_{\alpha}}^{\oplus r}$
for any $\alpha$,
$\sharp\{i|\: t_i|_M\cap U_{\alpha}\neq\emptyset\}\leq 1$
for any $\alpha$ and
$\sharp\{\alpha | \: t_i|_M\cap U_{\alpha}\neq\emptyset\}\leq 1$
for any $i$.
Take a relative tangent vector field
$v\in\Theta_{M^{\balpha}_{\cC/T}(\tilde{\bt},r,d)
/T\times\Lambda^{(n)}_r(d)}(M)$.
$v$ corresponds to a member
$(E_{\epsilon},\nabla_{\epsilon},\{(l_{\epsilon})^{(i)}_j\})
\in M^{\balpha}_{\cC/T}(\tilde{\bt},r,d)(\Spec{\mathcal O}_M[\epsilon])$
such that
$(E_{\epsilon},\nabla_{\epsilon},\{(l_{\epsilon})^{(i)}_j\})\otimes
{\mathcal O}_M[\epsilon]/(\epsilon)\cong
(\tilde{E},\tilde{\nabla},\{\tilde{l}^{(i)}_j\})|_{\cC_M}$,
where ${\mathcal O}_M[\epsilon]={\mathcal O}_M[t]/(t^2)$.
There is an isomorphism
\[
 \varphi_{\alpha}:E_{\epsilon}|
 _{U_{\alpha}\times\Spec{\mathcal O}_M[\epsilon]}
 \stackrel{\sim}\rightarrow
 {\mathcal O}_{U_{\alpha}\times\Spec{\mathcal O}_M[\epsilon]}^{\oplus r}
 \stackrel{\sim}\rightarrow
 \tilde{E}|_{U_{\alpha}}\otimes{\mathcal O}_M[\epsilon]
\]
such that
$\varphi_{\alpha}\otimes{\mathcal O}_M[\epsilon]/(\epsilon):
E_{\epsilon}\otimes{\mathcal O}_M[\epsilon]/(\epsilon)|_{U_{\alpha}}
\stackrel{\sim}\rightarrow\tilde{E}|_{U_{\alpha}}
\otimes{\mathcal O}_M[\epsilon]/(\epsilon)
=\tilde{E}|_{U_{\alpha}}$
is the given isomorphism and that
$\varphi_{\alpha}|_{t_i\otimes{\mathcal O}_M[\epsilon]}((l_{\epsilon})^{(i)}_j)
=\tilde{l}^{(i)}_j|_{U_{\alpha}\times\Spec{\mathcal O}_M[\epsilon]}$
if $t_i|_M\cap U_{\alpha}\neq\emptyset$.
We put
\[
 u_{\alpha\beta}:=\varphi_{\alpha}\circ\varphi_{\beta}^{-1}
 -\mathrm{id}_{\tilde{E}|_{U_{\alpha\beta}\otimes{\mathcal O}_M[\epsilon]}},
 \quad
 v_{\alpha}:=(\varphi_{\alpha}\otimes\mathrm{id})\circ
 \nabla_{\epsilon}|_{U_{\alpha}\otimes{\mathcal O}_M[\epsilon]}
 \circ\varphi_{\alpha}^{-1}
 -\tilde{\nabla}|_{U_{\alpha}\otimes{\mathcal O}_M[\epsilon]}.
\]
Then $\{u_{\alpha\beta}\}\in C^1((\epsilon)\otimes{\mathcal F}^0_M)$,
$\{v_{\alpha}\}\in C^0((\epsilon)\otimes{\mathcal F}^1_M)$ and
\[
 d\{u_{\alpha\beta}\}=\{u_{\beta\gamma}-u_{\alpha\gamma}+u_{\alpha\beta}\}=0,
 \quad
 \nabla_{{\mathcal F}^{\bullet}}\{u_{\alpha\beta}\}=\{v_{\beta}-v_{\alpha}\}
 =d\{v_{\alpha}\}.
\]
So $[(\{u_{\alpha\beta}\},\{v_{\alpha}\})]$
determines an element
$\sigma_M(v)$ of $\mathbf{H}^1({\mathcal F}^{\bullet}_M)$.
We can check that $v\mapsto \sigma(v)$ determines
an isomorphism
\[
 \sigma_M:
 \Theta_{M^{\balpha}_{\cC/T}(\tilde{\bt},r,d)/T\times\Lambda^{(n)}_r(d)}(M)
 \stackrel{\sim}\longrightarrow
 \mathbf{H}^1({\mathcal F}^{\bullet}_M).
\]
These isomorphisms $\sigma_M$ induce a canonical isomorphism
\[
 \sigma:\Theta_{M^{\balpha}_{\cC/T}(\tilde{\bt},r,d)/T\times\Lambda^{(n)}_r(d)}
 \stackrel{\sim}\longrightarrow
 {\bf R}^1(\pi_{T\times\Lambda^{(n)}_r(d)})_*(\cF^{\bullet}),
\]
where
$\pi_{T\times\Lambda^{(n)}_r(d)}:\cC\ra T\times\Lambda^{(n)}_r(d)$
is the projection.
From the hypercohomology spectral sequence
$H^q(\cF^p\otimes k(x))\Rightarrow \bH^{p+q}(\cF^{\bullet}\otimes k(x))$,
we obtain an exact sequence
\[
 0\ra k(x)\ra H^0(\cF^0(x)) \ra H^0(\cF^1(x)) \ra \bH^1(\cF^{\bullet}(x))
 \ra H^1(\cF^0(x)) \ra H^1(\cF^1(x))\ra k(x)\ra 0
\]
for any $x\in M^{\balpha}_{\cC/T}(\tilde{\bt},r,d)$.
Here we denote ${\mathcal F}^0\otimes k(x)$ by
${\mathcal F}^0(x)$ and so on.
So we have $\dim \bH^1(\cF^{\bullet}(x))=-2\chi(\cF^0(x))+2$,
because $H^0(\cF^1(x))\cong H^1(\cF^0(x))^{\vee}$ and
$H^1(\cF^1(x))\cong H^0(\cF^0(x))^{\vee}$.
We define ${\mathcal E}_1$ by the exact sequence
\[
 0\longrightarrow{\mathcal E}_1\longrightarrow
 {\mathcal End}(\tilde{E}(x))\longrightarrow
 \bigoplus_{i=1}^n\Hom
 \left(\tilde{l}^{(i)}_1(x),\left(\tilde{l}^{(i)}_0/\tilde{l}^{(i)}_1\right)(x)\right)
 \longrightarrow 0.
\]
Inductively we define ${\mathcal E}_k$ by the exact sequence
\[
 0\longrightarrow{\mathcal E}_k\longrightarrow
 {\mathcal E}_{k-1}\longrightarrow
 \bigoplus_{i=1}^n\Hom\left(\tilde{l}^{(i)}_k(x),
 \left(\tilde{l}^{(i)}_{k-1}/\tilde{l}^{(i)}_k\right)(x)\right)
 \longrightarrow 0.
\]
Then we have ${\mathcal E}_{r-1}={\mathcal F}^0(x)$.
So we have
\[
 \chi(\cF^0(x))=\chi({\mathcal End}(\tilde{E}(x)))-nr(r-1)/2
 =r^2(1-g)-nr(r-1)/2.
\]
Thus we can see that every fiber of
$M^{\balpha}_{\cC/T}(\tilde{\bt},r,d)$
over $T\times\Lambda^{(n)}_r(d)$ is equidimensional and its dimension is
$2r^2(g-1)+nr(r-1)+2$ if it is non-empty.
\hfill $\square$

\begin{Proposition}\label{openness}
 Let $(E,\nabla_E,\{l^{(i)}_j\})\in{\mathcal M}_{\cC/T}(\tilde{\bt},r,d)(S)$
 be a flat family of parabolic connections over a noetherian scheme $S$
 over $\mathbf{C}$.
 Then the subset
 \[
  U^s:=\{s\in S| \text{$(E,\nabla_E,\{l^{(i)}_j\})\otimes \overline{k(s)}$
  is $\balpha$-stable}\}
 \]
 is a Zariski open subset of $S$.
 Take a flat family $(E,\nabla_E,\{l^{(i)}_j\})\in{\mathcal M}_C(\bt,\blambda)(S)$ of
 $(\bt,\blambda)$-parabolic connections over a noetherian scheme $S$
 over $\mathbf{C}$. Then
 \[
  U^{irr}:=\{s\in S| \text{$(E,\nabla_E,\{l^{(i)}_j\})\otimes \overline{k(s)}$
  is irreducible}\}
 \]
 is a Zariski open subset of $S$.
 Note that a parabolic connection $(E,\nabla_E,\{l^{(i)}_j\})$
is reducible if there is a subbundle
$0\neq F\subsetneq E$ such that 
$\nabla_E(F)\subset F\otimes\Omega^1_C(t_1+\ldots+t_n)$.
A parabolic connection $(E,\nabla_E,\{l^{(i)}_j\})$ is irreducible
if it is not reducible.
\end{Proposition}

\begin{proof}
Take $(E,\nabla_E,\{l^{(i)}_j\})\in{\mathcal M}(\tilde{\bt},r,d)(S)$.
As in the proof of Theorem \ref{moduli-exists},
$(E,E,\mathrm{id},\nabla,F_*(E))$ determines a member of
$\overline{\cM^{D(\tilde{\bt})}_{\cC/T\times\Lambda^{(n)}_r(d)}}
 (r,d,\{1\}_{1\leq i\leq nr})(S)$,
where $\overline{\cM^{D(\tilde{\bt})}_{\cC/T\times\Lambda^{(n)}_r(d)}}
(r,d,\{1\}_{1\leq i\leq nr})$
is the moduli functor of parabolic $\Lambda^1_{D(\tilde{\bt})/T}$-triples
without stability condition.
By [\cite{IIS-1}, Proposition 5.3],
 \[
  S^s=\{s\in S|\text{$(E,E,\mathrm{id},\nabla,F_*(E))\otimes \overline{k(s)}$
  is $(\balpha',\beta,\gamma)$-stable}\}
 \]
is a Zariski open subset of $S$.
As in the proof of Theorem \ref{moduli-exists},
$(E,\nabla,\{l^{(i)}_j\})\otimes \overline{k(s)}$ is $\balpha$-stable if and only if
the corresponding parabolic $\Lambda^1_{D(\tilde{\bt})_s}$-triple
is $(\balpha',\beta,\gamma)$-stable.
So we have $S^s=U^s$.

Next take a flat family $(E,\nabla_E,\{l^{(i)}_j\})\in{\mathcal M}_C(\bt,\blambda)(S)$ of
$(\bt,\blambda)$-parabolic connections.
Take any subbundle $F\subset E\otimes \overline{k(s)}$ such that
$(\nabla_E\otimes \overline{k(s)})(F)\subset F\otimes\Omega^1_C(t_1+\cdots+t_n)$.
We have a commutative diagram
\[
 \begin{CD}
  F|_{t_i}  @>{\res_{t_i}(\nabla_F)}>>  F|_{t_i} \\
  @VVV  @VVV \\
  E|_{t_i}\otimes \overline{k(s)} @>\res_{t_i}(\nabla_E\otimes \overline{k(s)})>>
  E|_{t_i}\otimes \overline{k(s)}.
 \end{CD}
\]
Then we can see that the eigenvalues $\{\mu^{(i)}_j\}$
of $\res_{t_i}(\nabla_F)$ are eigenvalues of $\res_{t_i}(\nabla_E\otimes\overline{k(s)})$.
So the set
\[
 {\mathcal P}:=\left\{(\rank (F),\deg(F))\left|
 \begin{array}{l}
  \text{$F$ is a non-zero proper subbundle of $E\otimes\overline{k(s)}$ such that} \\
  (\nabla_E\otimes\overline{k(s)})(F)\subset F\otimes\Omega^1_C(t_1+\cdots+t_n)
 \end{array}
 \right\}\right.
\]
is finite, because $\deg(F)=-\sum_{i,j}\mu^{(i)}_j$.
For $(r',d')\in {\mathcal P}$, put $P_{(r',d')}(m)=r'(m-g+1)+d'$ which is a polynomial in $m$.
Then consider the Quot-scheme
$Q_{(r',d')}:=\Quot^{\chi(E(m)\otimes k(s))-P_{(r',d')}(m)}_{E/C/\mathbf{C}}$
and let ${\mathcal F}^{(r',d')}\subset E_{Q_{(r',d')}}$ be the universal subsheaf.
Let $Y_{(r',d')}$ be the maximal closed subscheme of
$Q_{(r',d')}$ such that the composite
\[
 {\mathcal F}^{(r',d')}_{Y_{(r',d')}}\hookrightarrow E_{Y_{(r',d')}}
 \xrightarrow{(\nabla_E)_{Y_{(r',d')}}}
 E_{Y_{(r',d')}}\otimes\Omega^1_C(t_1+\cdots+t_n)
 \longrightarrow  (E_{Y_{(r',d')}}/{\mathcal F}^{(r',d')}_{Y_{(r',d')}})
 \otimes\Omega^1_C(t_1+\cdots+t_n)
\]
is zero.
Since the structure morphism $f_{(r',d')}:Y_{(r',d')}\rightarrow S$ is proper,
\[
 U^{irr}=S\setminus\bigcup_{(r',d')\in{\mathcal P}}f_{(r',d')}(Y_{(r',d')})
\]
is a Zariski open subset of $S$.
\end{proof}

\begin{Proposition}\label{openness-of-irreducibility}
The set
\[
 RP_r(C,\bt)^{irr}=
 \{ [\rho]\in RP_r(C,\bt) | \text{$\rho$ is an irreducible representation} \}
\]
is a Zariski open subset of $RP_r(C,\bt)$.
\end{Proposition}

\begin{proof}
Note that for the quotient morphism
\[
 p:\Hom(\pi_1(C\setminus\{t_1,\ldots,t_n\},*),GL_r(\mathbf{C}))
 \longrightarrow RP_r(C,\bt),
\]
$Z\subset RP_r(C,\bt)$ is a Zariski closed subset if and only if
$p^{-1}(Z)$ is a Zariski closed subset of
$\Hom(\pi_1(C\setminus\{t_1,\ldots,t_n\},*),GL_r(\mathbf{C}))$.
There is a universal representation $\tilde{\rho}$ given by
\begin{gather*}
 \tilde{\rho}(\alpha_i):
 {\mathcal O}_{\Hom(\pi_1(C\setminus\{t_1,\ldots,t_n\},*),GL_r(\mathbf{C}))}^{\oplus r}
 \stackrel{\sim}\longrightarrow
 {\mathcal O}_{\Hom(\pi_1(C\setminus\{t_1,\ldots,t_n\},*),GL_r(\mathbf{C}))}^{\oplus r} \\
 \tilde{\rho}(\beta_i):
 {\mathcal O}_{\Hom(\pi_1(C\setminus\{t_1,\ldots,t_n\},*),GL_r(\mathbf{C}))}^{\oplus r}
 \stackrel{\sim}\longrightarrow
 {\mathcal O}_{\Hom(\pi_1(C\setminus\{t_1,\ldots,t_n\},*),GL_r(\mathbf{C}))}^{\oplus r} \\
 \tilde{\rho}(\gamma_j):
 {\mathcal O}_{\Hom(\pi_1(C\setminus\{t_1,\ldots,t_n\},*),GL_r(\mathbf{C}))}^{\oplus r}
 \stackrel{\sim}\longrightarrow
 {\mathcal O}_{\Hom(\pi_1(C\setminus\{t_1,\ldots,t_n\},*),GL_r(\mathbf{C}))}^{\oplus r}
\end{gather*}
for $i=1,\ldots,g$ and $i=1,\ldots,n$.
For $0<r'<r$, consider the Grassmann variety
\[
 \mathrm{Grass}_{r'}
 ({\mathcal O}_{\Hom(\pi_1(C\setminus\{t_1,\ldots,t_n\},*),GL_r(\mathbf{C}))}^{\oplus r})
\]
and let
$V\subset\mathrm{Grass}_{r'}
 ({\mathcal O}_{\Hom(\pi_1(C\setminus\{t_1,\ldots,t_n\},*),GL_r(\mathbf{C}))}^{\oplus r})$
be the universal subbundle.
Let $Y_{r'}$ be the maximal closed subscheme of
$\mathrm{Grass}_{r'}
({\mathcal O}_{\Hom(\pi_1(C\setminus\{t_1,\ldots,t_n\},*),GL_r(\mathbf{C}))}^{\oplus r})$
such that
\begin{gather*}
 \tilde{\rho}(\alpha_i)_{Y_{r'}}(V_{Y_{r'}})\subset V_{Y_{r'}} \\
 \tilde{\rho}(\beta_i)_{Y_{r'}}(V_{Y_{r'}})\subset V_{Y_{r'}} \\
 \tilde{\rho}(\gamma_j)_{Y_{r'}}(V_{Y_{r'}})\subset V_{Y_{r'}}
\end{gather*}
for $i=1,\ldots,g$ and $j=1,\ldots,n$.
Since the structure morphism
\[
 f_{r'}:Y_{r'}\rightarrow  \Hom(\pi_1(C\setminus\{t_1,\ldots,t_n\},*),GL_r(\mathbf{C}))
\]
is a proper morphism,
$\bigcup_{r'}f_{r'}(Y_{r'})$
is a Zariski closed subset of $\Hom(\pi_1(C\setminus\{t_1,\ldots,t_n\},*),GL_r(\mathbf{C}))$
which is preserved by the adjoint action of $PGL_r(\mathbf{C})$ on
$\Hom(\pi_1(C\setminus\{t_1,\ldots,t_n\},*),GL_r(\mathbf{C}))$.
So $p(\bigcup_{r'}f_{r'}(Y_{r'}))$ is a Zariski closed subset of $RP_r(C,\bt)$.
Thus 
\[
 \{ [\rho]\in RP_r(C,\bt) | \text{$\rho$ is an irreducible representation} \}
 =RP_r(C,\bt)\setminus p(\bigcup_{r'}f_{r'}(Y_{r'}))
\]
is a Zariski open subset of $RP_r(C,\bt)$.
\end{proof}

\begin{Remark}\label{existence-of-irreducible-representation}
 Assume that $n\geq 3$ if $g=0$ and $n\geq 1$ if $g\geq 1$.
 Then $RP_r(C,\bt)^{irr}$ is non-empty
\end{Remark}

\begin{proof}
Note that the fundamental group $\pi_1(C\setminus\{t_1,\ldots,t_n\})$
is a free group generated by
$\alpha_1,\beta_1,\ldots,\alpha_g,\beta_g,\\
\gamma_1,\ldots,\gamma_{n-1}$.
If $g\geq 1$, consider the $GL_r(\mathbf{C})$-representation $\rho$ of
$\pi_1(C\setminus\{t_1,\ldots,t_n\})$ given by
\begin{gather*}
 \rho(\alpha_1)=
 \begin{pmatrix}
  \lambda & 1 & 0 & 0 & \cdots & 0 \\
  0 & \lambda & 1 & 0 & \cdots & 0\\
  \vdots & \vdots & \vdots & \vdots & \ddots & \vdots \\
  0 & 0 & 0 & 0 & \cdots & 1 \\
  0 & 0 & 0 & 0 & \cdots & \lambda
 \end{pmatrix},
 \quad
 \rho(\beta_1)=
 \begin{pmatrix}
  \mu & 0 & 0 & 0 & \cdots & 0 \\
  1 & \mu & 0 & 0 & \cdots & 0\\
  0 & 1 & \mu & 0 & \cdots & 0 \\
  \vdots & \vdots & \vdots & \vdots & \ddots & \vdots \\
  0 & 0 & 0 & 0 & \cdots & 0 \\
  0 & 0 & 0 & 0 & \cdots & \mu
 \end{pmatrix}  \\
 \text{$\rho(\alpha_2),\rho(\beta_2)\ldots,\rho(\alpha_g),\rho(\beta_g),
 \rho(\gamma_1),\ldots,\rho(\gamma_{n-1})$ are arbitrary matrix.}
\end{gather*}
Then we can see that $\rho$ is an irreducible representation.
If $g=0$, then consider the $GL_r(\mathbf{C})$-representation $\rho$ given by
\begin{gather*}
 \rho(\gamma_1)=
 \begin{pmatrix}
  \lambda & 1 & 0 & 0 & \cdots & 0 \\
  0 & \lambda & 1 & 0 & \cdots & 0\\
  \vdots & \vdots & \vdots & \vdots & \ddots & \vdots \\
  0 & 0 & 0 & 0 & \cdots & 1 \\
  0 & 0 & 0 & 0 & \cdots & \lambda
 \end{pmatrix},
 \quad
 \rho(\gamma_2)=
  \begin{pmatrix}
  \mu & 0 & 0 & 0 & \cdots & 0 \\
  1 & \mu & 0 & 0 & \cdots & 0\\
  0 & 1 & \mu & 0 & \cdots & 0 \\
  \vdots & \vdots & \vdots & \vdots & \ddots & \vdots \\
  0 & 0 & 0 & 0 & \cdots & 0 \\
  0 & 0 & 0 & 0 & \cdots & \mu
 \end{pmatrix}  \\
 \text{$\rho(\gamma_3),\ldots,\rho(\gamma_{n-1})$ are arbitrary matrix.}
\end{gather*}
Then we can see that $\rho$ is an irreducible representation.
\end{proof}

\section{Irreducibility of the moduli space}
\label{irreducibility of m-s}

In this section we will prove the irreducibility of
$M^{\balpha}_{\cC_x}(\tilde{\bt}_x,\blambda)$, which is the fiber
of $M^{\balpha}_{\cC/T}(\tilde{\bt},r,d)$ over
$(x,\blambda)\in T\times\Lambda^{(n)}_r(d)$.
We simply denote $\cC_x$ by $C$,
$\tilde{\bt}_x$ by $\bt=(t_1,\ldots,t_n)$
and $t_1+\cdots+t_n$ by $D(\bt)$.
First consider the morphism
\begin{gather}\label{determinant-morphism}
 \det:M^{\balpha}_C(\bt,\blambda)\longrightarrow
 M_C(\bt,\det(\blambda)) \\
 \quad (E,\nabla,\{l^{(i)}_j\})\mapsto (\det E,\det(\nabla)),
 \notag
\end{gather}
where $\det(\nabla)$ is the logarithmic connection on $\det E$
induced by $\nabla$ and
$\det(\blambda)=(\sum_{j=0}^{r-1}\lambda^{(i)}_j)^{1\leq i\leq n}
\in \Lambda^{(n)}_1(d)$.
Note that a line bundle with a connection is always stable.
We can also construct $M_C(\bt,\det(\blambda))$ as an affine
space bundle over $\Pic_C^d$ whose fiber is of dimension
\[
 h^0({\mathcal End}(\det E)\otimes\Omega^1_C)
 =h^0(\Omega^1_C)=g.
\]
Thus $M_C(\bt,\det(\blambda))$ is a smooth irreducible
variety of dimension $2g$.

We can prove the smoothness of the morphism (\ref{determinant-morphism}).
Indeed let $A$ be an artinian local ring over
$M_C(\bt,\det(\blambda))$ with maximal ideal $m$,
$I$ an ideal of $A$ satisfying $mI=0$ and take
any $A/I$-valued point $(E,\nabla,\{l^{(i)}_j\})$ of
$M^{\balpha}_C(\bt,\blambda)$ over $M_C(\bt,\det(\blambda))$.
Let $(L,\nabla_L)$ be the line bundle with a connection on $(C,\bt)$
which corresponds to the morphism
$\Spec A\rightarrow M_C(\bt,\det(\blambda))$.
We write
$(\bar{E},\bar{\nabla},\{\bar{l}^{(i)}_j\}):=
(E,\nabla,\{l^{(i)}_j\})\otimes A/m$
and put
\begin{align*}
 \cF^0_0&:=\left\{ s \in {\mathcal End}(\bar{E}) \left|
 \text{$\Tr(s)=0$ and $s(t_i)(\bar{l}^{(i)}_j)\subset \bar{l}^{(i)}_j$
 for any $i,j$}
 \right\}\right. \\ 
 \cF^1_0&:=\left\{ s \in {\mathcal End}(\bar{E})\otimes\Omega^1_X(D(\bt))
 \left|
 \text{$\Tr(s)=0$ and $s(t_i)(\bar{l}^{(i)}_j)\subset \bar{l}^{(i)}_{j+1}$
 for any $i,j$}
 \right\}\right. 
\end{align*}
We define a complex ${\mathcal F}_0^{\bullet}$ by
\[
 \nabla_{{\mathcal F}_0^{\bullet}}:{\mathcal F}^0_0
 \ni s \mapsto \bar{\nabla}\circ s -s\circ\bar{\nabla}
 \in {\mathcal F}^1_0.
\]
Then the obstruction class for the lifting of
$(E,\nabla,\{l^{(i)}_j\})$ to an $A$-valued point
$(\tilde{E},\tilde{\nabla},\{\tilde{l}^{(i)}_j\})$ of
$M^{\balpha}_C(\bt,\blambda)$
with $(\det(E),\det(\nabla))\cong(L,\nabla_L)$ lies in 
$\mathbf{H}^2({\mathcal F}^{\bullet}_0)\otimes I$.
We can see by the same way as the proof of
Theorem \ref{moduli-exists} that
\[
 {\bf H}^2({\mathcal F}_0^{\bullet})\cong
 \left(
 \ker(H^0({\mathcal F}^0_0)\rightarrow H^0({\mathcal F}^1_0))
 \right)^{\vee}
\]
and we can see by the stability of
$(\bar{E},\bar{\nabla},\{\bar{l}^{(i)}_j\})$ that
$\ker(H^0({\mathcal F}^0_0)\rightarrow H^0({\mathcal F}^1_0))=0$.
Hence the morphism $\det$ in (\ref{determinant-morphism})
is smooth.

From the above argument, it is sufficient to prove the
irreducibility of the fibers of the morphism (\ref{determinant-morphism})
$\det:M^{\balpha}_C(\bt,\blambda)\rightarrow M_C(\bt,\det(\blambda))$
in order to obtain the irreducibility of
$M^{\balpha}_C(\bt,\blambda)$.
So we fix a line bundle $L$ on $C$
and a connection
\[
 \nabla_L: L \lra L\otimes\Omega^1_C(t_1+\cdots+t_n).
\]
We set
\[
 M^{\balpha}_C(\bt,\blambda,L):=\left.\left\{
 (E,\nabla,\{l^{(i)}_j\})\in M^{\balpha}_C(\bt,\blambda)
 \right| (\det E,\det(\nabla))\cong (L,\nabla_L)
 \right\}.
\]
Then we can see that
$M^{\balpha}_C(\bt,\blambda,L)$ is smooth of equidimension
$(r-1)(2(r+1)(g-1)+rn)$.

\begin{Proposition}\label{irred-generic}
 Assume that $g\geq 2$ and $n\geq 1$.
 Then $M^{\balpha}_C(\bt,\blambda,L)$
 is an irreducible variety of dimension $(r-1)(2(r+1)(g-1)+rn)$.
\end{Proposition}

\begin{proof}
By taking elementary transform, we can obtain an isomorphism
\[
 \sigma:\cM_C(\bt,\blambda,L)\stackrel{\sim}\longrightarrow
 \cM_C(\bt,\blambda',L'),
\]
with $r$ and $\det L'$ coprime.
Consider the open subscheme
\[
 N:=\left.\left\{
 (E,\nabla,\{l^{(i)}_j\})\in M^{\balpha}_C(\bt,\blambda,L)
 \right| \text{$\sigma(E)$ is a stable vector bundle}
 \right\}
\]
of $M^{\balpha}_C(\bt,\blambda,L)$.

\noindent
{\bf Step 1.}
We will see that $N$ is irreducible. \\
Let $M_C(r,L')$ be the moduli space of stable vector bundles
of rank $r$ with the determinant $L'$.
As is well-known, $M_C(r,L')$ is a smooth irreducible variety
and there is a universal bundle ${\mathcal E}$ on
$C\times M_C(r,L')$.
Note that there is a line bundle ${\mathcal L}$ on $M_C(r,L')$
such that $\det({\mathcal E})\cong L\otimes{\mathcal L}$.
We can parameterize the parabolic structures on 
${\mathcal E}_s$ ($s\in M_C(r,L')$) by a product of flag schemes
\[
 U:=\prod_{i=1}^n \mathrm{Flag}({\mathcal E}|_{t_i\times M_C(r,L')}),
\]
which is obviously smooth and irreducible.
Let
$\{\tilde{l}^{(i)}_j\}$ be the universal family over $U$,
where
\[
 {\mathcal E}_U|_{t_i\times U}=\tilde{l}^{(i)}_0\supset
 \tilde{l}^{(i)}_1\supset\cdots\supset
 \tilde{l}^{(i)}_{r-1}\supset\tilde{l}^{(i)}_r=0
\]
is the filtration by subbundles for $i=1,\ldots,n$
such that
$\dim (\tilde{l}^{(i)}_j)_s/(\tilde{l}^{(i)}_{j+1})_s=1$
for any $i,j$ and $s\in U$.
Consider the functor
$\mathrm{Conn}_U:(Sch/U)\rightarrow (Sets)$
defined by
\[
 \mathrm{Conn}_U(S)=
 \left\{ \nabla:\cE_S \ra \cE_S\otimes\Omega^1_C(D(\bt)) \left|
 \begin{array}{l}
  \text{$\nabla$ is a relative connection satisfying}\\
  \text{$(\res_{t_i\times S}(\nabla)-\lambda^{(i)}_j)
  (\tilde{l}^{(i)}_j)_S
  \subset (\tilde{l}^{(i)}_{j+1})_S$ for any $i,j$}\\
  \text{and $\det(\nabla)=(\nabla_L)_S\otimes{\mathcal L}_S$}
 \end{array}
 \right\}\right..
\]
We can define a morphism of functors
\begin{align*}
 &\mathrm{Conn}_U \longrightarrow
 \Quot_{\Lambda^1_{D(\bt)}\otimes{\mathcal E}/C\times U/U} \\
 & \nabla \mapsto
 \left[ \Lambda^1_{D(\bt)}\otimes{\mathcal E}_U
 \ni (a,v)\otimes e \mapsto ae+\nabla_v(e)\in {\mathcal E}_U\right],
\end{align*}
which is representable by an immersion.
So there is a locally closed subscheme $Y$ of
$\Quot_{\Lambda^1_{D(\bt)}\otimes{\mathcal E}/C\times U/U}$
which represents the functor $\mathrm{Conn}_U$.
We will show that $Y$ is smooth over $U$.
So let
$\tilde{\nabla}:{\mathcal E}_Y\rightarrow
{\mathcal E}_Y\otimes\Omega^1_C(D(\bt))$
be the universal connection and put
\begin{align*}
 \cF_0^0&:=\left\{ s \in {\mathcal End}({\mathcal E}) \left|
 \text{$\Tr(s)=0$ and $s|_{t_i\times Y}((\tilde{l}^{(i)}_j)_Y)
 \subset (\tilde{l}^{(i)}_j)_Y$ for any $i,j$}
 \right\}\right. \\
 \cF_0^1&:=\left\{
 s \in {\mathcal End}({\mathcal E})\otimes\Omega^1_C(D(\bt))
 \left|
 \text{$\Tr(s)=0$ and $\res_{t_i\times Y}(s)((\tilde{l}^{(i)}_j)_Y)
 \subset (\tilde{l}^{(i)}_{j+1})_Y$ for any $i,j$}
 \right\}\right. \\
 &\nabla_{\cF_0^{\bullet}}: \cF_0^0 \ni s\mapsto 
 \tilde{\nabla}\circ s - s\circ\tilde{\nabla} \in \cF_0^1.
\end{align*}
Let $A$ be an aritinian local ring over $U$ with maximal ideal $m$
and $I$ an ideal of $A$ satisfying $mI=0$.
Take any $A/I$-valued point $x$ of $Y$ over $U$.
We put $\bar{x}:=x\otimes A/m$.
Then the obstruction class for the lifting of $x$
to an $A$-valued point of $Y$ lies in
$H^1(\cF_0^1\otimes k(\bar{x}))\otimes I$.
We can see that
$H^1(\cF_0^1\otimes k(\bar{x}))\cong
H^0(\cF_0^0\otimes k(\bar{x}))^{\vee}=0$
because ${\mathcal E}\otimes k(\bar{x})$ is a stable bundle.
Thus $Y$ is smooth over $U$.
We can see that
the fiber $Y_s$ of any point $s\in U$ is isomorphic to the affine space
isomorphic to $H^0(\cF_0^1\otimes k(s))$.
So $Y$ is irreducible.
Consider the open subscheme
\[
 Y':=\left\{ x\in Y \left|
 \text{$\sigma^{-1}({\mathcal E}_x,\tilde{\nabla}\otimes k(x),
 \{(\tilde{l}^{(i)}_j)_x\})$ is $\balpha$-stable}
 \right\}\right.
\]
of $Y$.
Then we obtain a morphism $Y'\ra M^{\balpha}_C(\bt,\blambda,L)$
whose image is just $N$.
Thus $N$ is irreducible.

\noindent
{\bf Step 2.}
$\dim(M^{\balpha}_C(\bt,\blambda,L)\setminus N)<
\dim M^{\balpha}_C(\bt,\blambda,L)
=(r-1)(2(r+1)(g-1)+rn)$.

For the proof of Step 2, it is sufficient to show that
there is an algebraic scheme $\tilde{B}$ with $\dim\tilde{B}<(r-1)(2(r+1)(g-1)+rn)$
and a surjection
$\tilde{B}\rightarrow \sigma(M^{\balpha}_C(\bt,\blambda,L)\setminus N)$.
Take any member
$(E,\nabla_E,\{l^{(i)}_j\})\in
\sigma(M^{\balpha}_C(\bt,\blambda,L)\setminus N)$.
Let
\[
 0=E_0\subset E_1\subset E_2\subset\cdots\subset E_s=E
\]
be the Harder-Narasimhan filtration of $E$.
Note that $s>1$ because $r$ and $\deg E=\deg L'$ are coprime.
Put $\bar{E}_k:=E_k/E_{k-1}$ for $k=1,\ldots,s$.
For each $\bar{E}_k$, there is a Jordan-H\"older filtration
\[
 0=E^{(0)}_k\subset E^{(1)}_k\subset E^{(2)}_k
 \subset\cdots\subset E^{(m_k)}_k=\bar{E}_k.
\]
We put $\bar{E}^{(i)}_k:=E^{(i)}_k/E^{(i-1)}_k$ for $i=1,\ldots,m_k$.
Put $r^{(j)}_k=\rank \bar{E}^{(j)}_k$ and
$d^{(j)}_k=\deg\bar{E}^{(j)}_k$.
We can parameterize each $\bar{E}^{(j)}_k$
by the moduli space $M(r^{(j)}_k,d^{(j)}_k)$
of stable bundles of rank $r^{(j)}_k$ and degree
$d^{(j)}_k$
whose dimension is $\dim\Ext^1(\bar{E}^{(j)}_k,\bar{E}^{(j)}_k)$.
Replacing $M(r_k^{(j)},d_k^{(j)})$ by an \'etale covering,
we may assume that there is a universal bundle
$\bar{\mathcal E}^{(j)}_k$ on $C\times M(r^{(j)}_k,d^{(j)}_k)$.
If we put
\[
 X:=\left\{ (x^{(j)}_k)
 \in \prod_{1\leq k\leq s, 1\leq j\leq m_k}
 M(r^{(j)}_k,d^{(j)}_k) \left|
 \bigotimes_{j,k}
 \det\left(\bar{\mathcal E}^{(j)}_k\right)_{x^{(j)}_k} \cong L
 \right\}\right.,
\]
then $X$ is smooth of dimension
$-g+\sum_{k=1}^s\sum_{j=1}^{m_k}
\dim\Ext^1(\bar{E}^{(j)}_k,\bar{E}^{(j)}_k)$.
Replacing $X$ by its stratification,
we may assume that the relative $\Ext$-sheaf
$\Ext^1_{C\times X/X}
(\bar{\mathcal E}^{(2)}_k,\bar{\mathcal E}^{(1)}_k)$
becomes locally free.
We put
$\tilde{W}^{(2)}_k:=X$ and
${\mathcal E}^{(2)}_k:=
\bar{\mathcal E}^{(1)}_k\oplus\bar{\mathcal E}^{(2)}_k$
if the extension
\[
 0\longrightarrow \bar{E}^{(1)}_k\longrightarrow
 E^{(2)}_k \longrightarrow \bar{E}^{(2)}_k
 \longrightarrow 0
\]
splits.
Otherwise we put
\[
 \tilde{W}^{(2)}_k:=
 \mathbf{P}_*\left(
 \Ext^1_{C\times X/X}
(\bar{\mathcal E}^{(2)}_k,\bar{\mathcal E}^{(1)}_k)
 \right)
\]
and take a universal extension
\[
 0\longrightarrow (\bar{\mathcal E}^{(1)}_k)_{\tilde{W}^{(2)}_k}
 \otimes{\mathcal L}^{(2)}_k
 \longrightarrow {\mathcal E}^{(2)}_k
 \longrightarrow (\bar{\mathcal E}^{(2)}_k)_{\tilde{W}^{(2)}_k}
 \longrightarrow 0
\]
for some line bundle ${\mathcal L}^{(2)}_k$
on $\tilde{W}^{(2)}_k$.
We define $\tilde{W}^{(j)}_k$ and ${\mathcal E}^{(j)}_k$
inductively on $j$ as follows:
Replacing $\tilde{W}^{(j-1)}_k$ by its stratification,
we may assume that the relative $\Ext$-sheaf
$\Ext^1_{C\times \tilde{W}^{(j-1)}_k/\tilde{W}^{(j-1)}_k}
 ((\bar{\mathcal E}^{(j)}_k)_{\tilde{W}^{(j-1)}_k},
 {\mathcal E}^{(j-1)}_k)$
is locally free.
Then put $\tilde{W}^{(j)}_k:=\tilde{W}^{(j-1)}_k$ and
${\mathcal E}^{(j)}_k:=
{\mathcal E}^{(j-1)}_k\oplus (\bar{\mathcal E}^{(j)}_k)_{\tilde{W}^{(j)}_k}$
if the extension
\[
 0\longrightarrow E^{(j-1)}_k \longrightarrow E^{(j)}_k
 \longrightarrow \bar{E}^{(j)}_k \longrightarrow 0
\]
splits and otherwise we put
\[
 \tilde{W}^{(j)}_k:=\mathbf{P}_*\left(
 \Ext^1_{C\times \tilde{W}^{(j-1)}_k/\tilde{W}^{(j-1)}_k}
 ((\bar{\mathcal E}^{(j)}_k)_{\tilde{W}^{(j-1)}_k},
 {\mathcal E}^{(j-1)}_k)
 \right)
\]
and take a universal extension
\[
 0\longrightarrow ({\mathcal E}^{(j-1)}_k)_{\tilde{W}^{(j)}_k}
 \otimes{\mathcal L}^{(j)}_k
 \longrightarrow {\mathcal E}^{(j)}_k
 \longrightarrow (\bar{\mathcal E}^{(j)}_k)_{\tilde{W}^{(j)}_k}
 \longrightarrow 0
\]
for some line bundle ${\mathcal L}^{(j)}_k$ on $\tilde{W}^{(j)}_k$.
We can see by the construction that the relative dimension of
$\tilde{W}^{(m_k)}_k$ over $X$ at the point corresponding to the
extensions
\[
 0\longrightarrow E^{(j-1)}_k \longrightarrow E^{(j)}_k
 \longrightarrow \bar{E}^{(j)}_k \longrightarrow 0
 \quad (j=2,\ldots,m_k, \; E^{(1)}_k=\bar{E}^{(1)}_k)
\]
is at most
$1-m_k+\sum_{1\leq i<j\leq m_k}\dim\Ext^1(\bar{E}^{(j)}_k,\bar{E}^{(i)}_k)$.
We put 
\[
 W:=\tilde{W}^{(m_1)}_1\times_X\cdots\times_X\tilde{W}^{(m_s)}_s
\]
and $\bar{\mathcal E}_k:=({\mathcal E}^{(m_k)}_k)_W$
for $k=1,\ldots,s$.
Replacing $W$ by its stratification, we assume that
the relative $\Ext$-sheaf
$\Ext^1_{C\times W/W}(\bar{\mathcal E}_2,\bar{\mathcal E}_1)$
is locally free.
Then we put
$W_2:=W$ and
${\mathcal E}_2:=\bar{\mathcal E}_1\oplus\bar{\mathcal E}_2$
if the extension
\[
 0\longrightarrow \bar{E}_1 \longrightarrow E_2
 \longrightarrow \bar{E}_2 \longrightarrow 0
\]
splits and otherwise we put
\[
 W_2:=\mathbf{P}_*\left(
 \Ext^1_{C\times W/W}(\bar{\mathcal E}_2,\bar{\mathcal E}_1)
 \right)
\]
and take a universal extension
\[
 0\longrightarrow (\bar{\mathcal E}_1)_{W_2}\otimes{\mathcal L}_2
 \longrightarrow {\mathcal E}_2 \longrightarrow
 (\bar{\mathcal E}_2)_{W_2} \longrightarrow 0
\]
for some line bundle ${\mathcal L}_2$.
We define $W_k$ and ${\mathcal E}_k$ inductively as follows:
Replacing $W_{k-1}$ by its stratification,
we assume that the relative $\Ext$-sheaf
$\Ext^1_{C\times W_{k-1}/W_{k-1}}
((\bar{\mathcal E}_k)_{W_{k-1}},{\mathcal E}_{k-1})$
is locally free.
Then we put
$W_k:=W_{k-1}$ and
${\mathcal E}_k:=\bar{\mathcal E}_k\oplus{\mathcal E}_{k-1}$
if the extension
\[
 0\longrightarrow E_{k-1} \longrightarrow E_k
 \longrightarrow \bar{E}_k \longrightarrow 0
\]
splits and otherwise we put
\[
 W_k:=\mathbf{P}_*\left(
 \Ext^1_{C\times W_{k-1}/W_{k-1}}
 ((\bar{\mathcal E}_k)_{W_{k-1}},{\mathcal E}_{k-1})
 \right)
\]
and take a universal extension
\[
 0\longrightarrow ({\mathcal E}_{k-1})_{W_k}\otimes{\mathcal L}_k
 \longrightarrow {\mathcal E}_k
 \longrightarrow (\bar{\mathcal E}_k)_{W_k}
 \longrightarrow 0
\]
for some line bundle ${\mathcal L}_k$ on $W_k$.
We put $r_k:=\rank \bar{E}_k$
for $k=1,\ldots,s$.

\noindent
{\bf Claim.}
$\dim \Ext^1(\bar{E}_k,\bar{E}_{k'})<r_kr_{k'}(2g-2)$
for $k'<k$.

Take a general point $t\in C$.
We prove this claim in the following three cases.

\noindent
{\bf Case 1.}
$\deg(\bar{E}_{k'}^{\vee}\otimes\bar{E}_k(t))>0$.

By Riemann-Roch theorem, we have
\[
 \dim \Hom(\bar{E}_{k}(t),\bar{E}_{k'})
 -\dim H^0((\bar{E}_{k'})^{\vee}\otimes\bar{E}_k(t)\otimes\omega_C)
 =r_{k'}r_k(1-g)+\deg((\bar{E}_{k}(t))^{\vee}\otimes\bar{E}_{k'}).
\]
Here we have $\Hom(\bar{E}_k(t),\bar{E}_{k'})=0$
because $\bar{E}_k(t)$ and $\bar{E}_{k'}$ are semistable and
$\deg((\bar{E}_k(t))^{\vee}\otimes\bar{E}_{k'})<0$.
So we have
\begin{align*}
 \dim\Ext^1(\bar{E}_k,\bar{E}_{k'})&=
 \dim H^0(\bar{E}_{k'}^{\vee}\otimes\bar{E}_k\otimes\omega_C) \\
 &\leq\dim H^0(\bar{E}_{k'}^{\vee}\otimes\bar{E}_k(t)\otimes\omega_C) \\
 &=r_kr_{k'}(g-1)-\deg((\bar{E}_{k}(t))^{\vee}\otimes\bar{E}_{k'}) \\
 &<r_kr_{k'}(g-1)+r_kr_{k'} \\
 &\leq r_kr_{k'}(2g-2). 
\end{align*}

\noindent
{\bf Case 2.}
$\deg(\bar{E}_{k'}^{\vee}\otimes\bar{E}_k(t))<0$.

Note that there is an exact sequence
\begin{gather*}
 0\longrightarrow
 H^0(\bar{E}_{k'}^{\vee}\otimes \bar{E}_k\otimes\omega_C(-p_1-\cdots-p_{2g-3}))
 \longrightarrow 
 H^0(\bar{E}_{k'}^{\vee}\otimes\bar{E}_k\otimes\omega_C(-p_1-\cdots-p_{2g-4})) \\
 \longrightarrow
 H^0(\bar{E}_{k'}^{\vee}\otimes\bar{E}_k\otimes\omega_C(-p_1-\cdots-p_{2g-4})|_{p_{2g-3}})
 \cong\mathbf{C}^{r_kr_{k'}}.
\end{gather*}
Since $\deg(\bar{E}_{k'}^{\vee}\otimes \bar{E}_k\otimes\omega_C(-p_1-\cdots-p_{2g-3}))<0$
and $\bar{E}_k$ and $\bar{E}_{k'}$ are semistable,
we have
\[
 H^0(\bar{E}_{k'}^{\vee}\otimes \bar{E}_k\otimes\omega_C(-p_1-\cdots-p_{2g-3}))=0.
\]
So we have
$\dim H^0(\bar{E}_{k'}^{\vee}\otimes\bar{E}_k\otimes
\omega_C(-p_1-\cdots-p_{2g-4}))\leq r_kr_{k'}$.
Inductively we can see that
\[
 \dim H^0(\bar{E}_{k'}^{\vee}\otimes\bar{E}_k\otimes\omega_C(-p_1-\cdots-p_{2g-2-m}))
 \leq r_kr_{k'}(m-1)
\]
for $2\leq m\leq 2g-2$, because there are exact sequences
\begin{gather*}
 0\longrightarrow
 H^0(\bar{E}_{k'}^{\vee}\otimes\bar{E}_k\otimes\omega_C(-p_1-\cdots-p_{2g-2-m+1}))
 \longrightarrow
 H^0(\bar{E}_{k'}^{\vee}\otimes\bar{E}_k\otimes\omega_C(-p_1-\cdots-p_{2g-2-m})) \\
 \longrightarrow
 H^0(\bar{E}_{k'}^{\vee}\otimes\bar{E}_k\otimes\omega_C(-p_1-\cdots-p_{2g-2-m})|_{p_{2g-2-m+1}})
 \cong\mathbf{C}^{r_kr_{k'}}
\end{gather*}
for $2\leq m\leq 2g-2$.
Thus we have
\[
 \dim\Ext^1(\bar{E}_k,\bar{E}_{k'})=
 \dim H^0(\bar{E}_{k'}^{\vee}\otimes\bar{E}_k\otimes\omega_C)\leq r_kr_{k'}(2g-3).
\]

\noindent
{\bf Case 3.}
$\deg(\bar{E}_{k'}^{\vee}\otimes\bar{E}_k(t))=0$.

Put ${\mathcal F}:=
\bar{E}_{k'}^{\vee}\otimes\bar{E}_k\otimes\omega_C(-p_1-\cdots-p_{2g-3})$.
Then ${\mathcal F}$ is semistable of degree $0$.
Let $0={\mathcal F}_0\subset{\mathcal F}_1
\subset\cdots\subset{\mathcal F}_l={\mathcal F}$
be a Jordan-H\"{o}lder filtration of ${\mathcal F}$.
If we put $t\in C$ sufficiently general, we may assume that
${\mathcal O}_C(t-p_{2g-3})\not\cong{\mathcal F}_i/{\mathcal F}_{i-1}$
for any $i$.
Then we have
\[
 \Hom({\mathcal O}_C(t-p_{2g-3}),\bar{E}_{k'}^{\vee}\otimes\bar{E}_k
 \otimes\omega_C(-p_1-\cdots-p_{2g-3}))=0.
\]
Since there is an exact sequence
\begin{gather*}
 0=H^0(\bar{E}_{k'}^{\vee}\otimes\bar{E}_k\otimes
 \omega_C(-p_1-\cdots-p_{2g-4}-t))
 \longrightarrow
 H^0(\bar{E}_{k'}^{\vee}\otimes\bar{E}_k\otimes
 \omega_C(-p_1-\cdots-p_{2g-4})) \\
 \longrightarrow
 H^0(\bar{E}_{k'}^{\vee}\otimes\bar{E}_k\otimes
\omega_C(-p_1-\cdots-p_{2g-4})|_t)\cong\mathbf{C}^{r_kr_{k'}},
\end{gather*}
we have
$\dim H^0(\bar{E}_{k'}^{\vee}\otimes\bar{E}_k\otimes
\omega_C(-p_1-\cdots-p_{2g-4}))\leq r_kr_{k'}$.
Inductively we can see that
\[
 \dim H^0(\bar{E}_{k'}^{\vee}\otimes\bar{E}_k\otimes
 \omega_C(-p_1-\cdots-p_{2g-2-m}))\leq r_kr_{k'}(m-1)
\]
for $2\leq m\leq 2g-2$.
So we have
$\dim\Ext^1(\bar{E}_k,\bar{E}_{k'})=
\dim H^0(\bar{E}_{k'}^{\vee}\otimes\bar{E}_k\otimes\omega_C)
\leq r_kr_{k'}(2g-3)$.

Thus the claim has been proved.

We can see that the relative dimension of $W_k$ over $W_{k-1}$
at the point corresponding to the extension
\[
 0\longrightarrow E_{k-1} \longrightarrow E_k \longrightarrow
 \bar{E}_k \longrightarrow 0
\]
is at most
\[
 \dim\Ext^1(\bar{E}_k,E_{k-1})-1 \leq
 -1+\sum_{k'\leq k-1}\dim\Ext^1(\bar{E}_k,\bar{E}_{k'})
  <-1+\sum_{1\leq k'<k}\left(r_{k'}r_k(2g-2) \right)
\]
if
$0\rightarrow E_{k-1}\rightarrow E_k\rightarrow
\bar{E}_k\rightarrow 0$ does not split.
If the extension
$0\rightarrow E_{k-1}\rightarrow E_k\rightarrow
\bar{E}_k\rightarrow 0$ splits,
the relative dimension of $W_k$ over $W_{k-1}$ is at most
$-1+\sum_{1\leq k'<k}r_{k'}r_k(2g-2)$.
Then we can see that the dimension of $W_s$ at
the point corresponding to the extensions
\begin{gather*}
 0=E_0\subset E_1\subset\cdots E_s=E,
 \quad \bar{E}_k=E_k/E_{k-1} \\
 0=E^{(0)}_k\subset E^{(1)}_k\subset\cdots\subset
 E^{(m_k)}_k=\bar{E}_k,
 \quad \bar{E}^{(j)}_k=E^{(j)}_k/E^{(j-1)}_k
\end{gather*}
is at most
\begin{align*}
 &-g+\sum_{k=1}^s\sum_{j=1}^{m_k}\dim\Ext^1(\bar{E}^{(j)}_k,\bar{E}^{(j)}_k)
  +\sum_{k=1}^s\left(1-m_k+\sum_{1\leq i<j\leq m_k}
  \dim\Ext^1(\bar{E}^{(j)}_k,\bar{E}^{(i)}_k)\right) \\
  &+\sum_{k=2}^s\left(-1+\sum_{1\leq k'<k}
  \dim\Ext^1(\bar{E}_{k},\bar{E}_{k'})\right) \\
  <&-g+\sum_{k=1}^s\sum_{j=1}^{m_k}\dim\Ext^1(\bar{E}^{(j)}_k,\bar{E}^{(j)}_k)
  +\sum_{k=1}^s\left(1-m_k+\sum_{1\leq i<j\leq m_k}
  \dim\Ext^1(\bar{E}^{(j)}_k,\bar{E}^{(i)}_k)\right) \\
  &\quad +1-s+\sum_{1\leq k'<k\leq s}r_{k'}r_k(2g-2) \\
 \leq& \sum_{k=1}^s\sum_{i=1}^{m_k}((r^{(i)}_k)^2(g-1)+1)
  +\sum_{k=1}^s
  \left(1-m_k
  +\sum_{1\leq i< j\leq m_k}(r^{(i)}_kr^{(j)}_k(g-1)+1)\right) \\
  &\quad -g+1-s+\sum_{1\leq k'<k\leq s}
  \sum_{1\leq i\leq m_{k},1\leq j\leq m_{k'}}
  r^{(i)}_{k}r^{(j)}_{k'}(2g-2) \\
 \leq & \sum_{k=1}^s\sum_{i=1}^{m_k}(r^{(i)}_k)^2(g-1)
  +\sum_{k=1}^s\sum_{1\leq i< j\leq m_k}2r^{(i)}_kr^{(j)}_k(g-1) \\
  & \quad -g+1+\sum_{k<k'}\sum_{1\leq i\leq m_k,1\leq j\leq m_{k'}}
   2r^{(i)}_{k'}r^{(j)}_{k}(g-1) \\
 =& (r^2-1)(g-1)
\end{align*}
because we have $s>1$.

Consider the product of flag schemes
\[
 Z:=\prod_{i=1}^n \mathrm{Flag}({\mathcal E}_s|_{t_i\times W_s})
\]
over $W_s$ and put ${\mathcal E}:=({\mathcal E}_s)_Z$.
Then there is a universal parabolic structure
\[
 {\mathcal E}|_{t_i\times Z}=\tilde{l}^{(i)}_0\supset\tilde{l}^{(i)}_1
 \supset\cdots\supset\tilde{l}^{(i)}_{r-1}\supset\tilde{l}^{(i)}_r=0
\]
for each $i=1,\ldots,n$.
The dimension of $Z$ at the point corresponding to
the extensions
\begin{gather*}
 0=E_0\subset E_1\subset\cdots E_s=E,
 \quad \bar{E}_k=E_k/E_{k-1} \\
 0=E^{(0)}_k\subset E^{(1)}_k\subset\cdots\subset
 E^{(m_k)}_k=\bar{E}_k,
 \quad \bar{E}^{(j)}_k=E^{(j)}_k/E^{(j-1)}_k
\end{gather*}
and the parabolic structure
\[
 E|_{t_i}=l^{(i)}_0\supset l^{(i)}_1\supset\cdots\supset
 l^{(i)}_{r-1}\supset l^{(i)}_r=0
 \quad (i=1,\ldots,n)
\]
is less than $(r^2-1)(g-1)+r(r-1)n/2=(r-1)((r+1)(g-1)+rn/2)$.

Consider the functor
$\mathrm{Conn}_Z:(Sch/Z)\rightarrow (Sets)$
defined by
\[
 \mathrm{Conn}_Z(S)=\left\{
 \nabla:{\mathcal E}_S\rightarrow{\mathcal E}_S\otimes\Omega^1_C(D(\bt))
 \left|
 \begin{array}{l}
  \text{$\nabla$ is a connection such that} \\
  \text{$(\res_{t_i\times S}(\nabla)-\lambda^{(i)}_j)
  ((\tilde{l}^{(i)}_j)_S)\subset(\tilde{l}^{(i)}_{j+1})_S$
  for any $i,j$} \\
  \text{and $\det(\nabla)=\nabla_L\otimes{\mathcal L}$
  via an identification} \\
  \text{${\mathcal E}=L\otimes{\mathcal L}$
  for some line bundle ${\mathcal L}$ on $S$}
 \end{array}
 \right\}\right..
\]
We can see that $\mathrm{Conn}_Z$ can be represented by a
scheme $B$ of finite type over $Z$.
Let
\[
 \tilde{\nabla}:{\mathcal E}_B\longrightarrow
 {\mathcal E}_B\otimes\Omega^1_C(D(\bt))
\]
be the universal connection.
If we put
\begin{align*}
 \cF_0^0&:=\left\{ u \in {\mathcal End}({\mathcal E}) \left|
 \text{$\Tr(u)=0$ and $u|_{t_i\times Y}((\tilde{l}^{(i)}_j)_Y)
 \subset (\tilde{l}^{(i)}_j)_Y$ for any $i,j$}
 \right\}\right. \\
 \cF_0^1&:=\left\{
 u \in {\mathcal End}({\mathcal E})\otimes\Omega^1_C(D(\bt))
 \left|
 \text{$\Tr(u)=0$ and $\res_{t_i\times Y}(u)((\tilde{l}^{(i)}_j)_Y)
 \subset (\tilde{l}^{(i)}_{j+1})_Y$ for any $i,j$}
 \right\}\right. \\
 &\nabla_{{\mathcal F}_0^{\bullet}}:\cF_0^0:
 \ni u \mapsto \nabla_{\cF_0^{\bullet}}\circ u-u\circ\nabla_{\cF^{\bullet}}
 \in\cF_0^1,
\end{align*}
then the fiber $B_x$ of $B$ over $x\in Z$
is isomorphic to the affine space
$H^0(\cF_0^1\otimes k(x))$.
Put
\[
 B':=\left\{ x\in B \left|
 \text{$\sigma^{-1}\left({\mathcal E}_x,\tilde{\nabla}\otimes k(x),
 \left\{\tilde{l}^{(i)}_j\otimes k(x)\right\}\right)$ is $\balpha$-stable}
 \right\}\right..
\]
Then there is a canonical morphism
\[
 \psi:B'\longrightarrow M^{\balpha}_C(\bt,\blambda,L)\setminus N.
\]
and the homomorphism
\[
 H^0(\nabla_{\cF^{\bullet}_0}\otimes k(x)):
 H^0(\cF_0^0\otimes k(x))\longrightarrow
 H^0(\cF_0^1\otimes k(x))
\]
is injective for $x\in B'$, because
$\End({\mathcal E}\otimes k(x),\tilde{\nabla}\otimes k(x),\{\tilde{l}^{(i)}_j\otimes k(x)\})
\cong k(x)$.
If we put
\[
 U:=\left\{
 u:{\mathcal E}\otimes k(x)\stackrel{\sim}\longrightarrow{\mathcal E}\otimes k(x)
 \left|
 \begin{array}{l}
  \text{$u(\tilde{l}^{(i)}_j\otimes k(x))\subset\tilde{l}^{(i)}_j\otimes k(x)$
  and $\wedge^ru=\mathrm{id}$}
 \end{array}
 \right\}\right.,
\]
for $x\in Z$, then we have a morphism
\[
 \iota:U\rightarrow B'_x; \quad
 u \mapsto u\circ(\tilde{\nabla}\otimes k(x))\circ u^{-1}.
\]
We can see that the tangent map of $\iota$ at $\mathrm{id}$ is just
the homomorphism $H^0(\nabla_{\cF^{\bullet}_0}\otimes k(x))$
which is injective.
We may assume that the morphism
$\psi:B'\rightarrow\overline{\psi(B')}$
is flat by replacing $B'$ by a disjoint union of subschemes.
So we have for any $x\in B'$
\begin{align*}
 \dim_{\psi(x)}\overline{\psi(B')}&= \dim_x(B')-\dim_x\psi^{-1}(\psi(x)) \\
 &\leq \dim(Z)+\dim_xB'_x-\dim_x\iota(U) \\
 &\leq\dim(Z)+\dim H^0(\cF_0^1\otimes k(x))-\dim H^0(\cF_0^0\otimes k(x)) \\
 &<(r-1)(2(r+1)(g-1)+rn)
\end{align*}
because
$\dim\coker(H^0(\cF_0^0\otimes k(x))\rightarrow H^0(\cF_0^1\otimes k(x)))
=(r-1)((r+1)(g-1)+rn/2)$.
Note that $M^{\balpha}_C(\bt,\blambda,L)\setminus N$
can be covered by a finite union of such $\psi(B')$'s.
So we have
\[
 \dim (M^{\balpha}_C(\bt,\blambda,L)\setminus N)
 <(r-1)(2(r+1)(g-1)+rn).
\]

So Step 2 is proved.
By Step 1, Step 2 and the fact that the moduli space
$M^{\balpha}_C(\bt,\blambda,L)$ is smooth of equidimension
$(r-1)(2(r+1)(g-1)+rn)$,
the moduli space $M^{\balpha}_C(\bt,\blambda,L)$ is irreducible.
\end{proof}

\begin{Proposition}\label{irred-g=1}
 Assume that $g=1$ and $n\geq 2$.
 Then $M^{\balpha}_C(\bt,\blambda,L)$
 is an irreducible variety of dimension $r(r-1)n$.
\end{Proposition}

\begin{proof}
Composing certain elementary transforms,
we obtain an isomorphism
\[
 \sigma:\cM_C(\bt,\blambda,L) \stackrel{\sim}\lra
 \cM_C(\bt,\blambda',L'),
\]
where $r$ and $\deg L'$ are coprime.
Put
\[
 N:=\left\{ \left.
 (E,\nabla,\{l^{(i)}_j\})
 \in M^{\balpha}_C(\bt,\blambda,L) \right|
 \text{$\sigma(E)$ is a stable bundle} \right\}.
\]
As in the proof of the previous proposition,
we can show that $N$ is irreducible.
So it suffices to show that
$\dim(M^{\balpha}_C(\bt,\blambda,L)\setminus N)<r(r-1)n$
because $M^{\balpha}_C(\bt,\blambda,L)$ is smooth of
equidimension $r(r-1)n$.
Take any member
$(E,\nabla,\{l^{(i)}_j\})
\in \sigma(M^{\balpha}_C(\bt,\blambda,L)\setminus N)$.
Let $F_s\subset F_{s-1}\subset\cdots\subset F_1=E$
be the Harder-Narasimhan filtration of $E$.
Note that $s>1$.
We can inductively see that the extension
\[
 0 \longrightarrow F_p/F_{p+1} \longrightarrow
 E/F_{p+1} \longrightarrow E/F_p \longrightarrow 0
\]
must split for $p=2,\ldots,s$, where we put $F_{s+1}=0$.
Then we have a decomposition
$E\cong E_1\oplus\cdots\oplus E_s$,
where each $E_p$ is semistable and
$\mu(E_1)<\mu(E_2)<\cdots<\mu(E_s)$.
We can write
$E_p=\bigoplus_{j=1}^{r_p}F_{p,j}$,
where each $F_{p,j}$ is a successive extension of a stable
bundle $F_p^{(j)}$ and $F_p^{(j)}\not\cong F_p^{(k)}$
for $j\neq k$.
We can see that
$F_{p,j}\cong G_1\oplus\cdots\oplus G_m$,
where each $G_i$ is a successive non-split extension of
$F_p^{(j)}$.
We put $r_p^{(j)}:=\rank F_p^{(j)}$ and
$d_p^{(j)}:=\deg F_p^{(j)}$.
Let $M_C(r_p^{(j)},d_p^{(j)})$ be the moduli space
of stable bundles on $C$ of rank $r_p^{(j)}$ and degree $d_p^{(j)}$.
Replacing $M_C(r_p^{(j)},d_p^{(j)})$ by an \'etale covering,
we may assume that there is a universal bundle
${\mathcal F}_p^{(j)}$ on  $C\times M_C(r_p^{(j)},d_p^{(j)})$.
We can construct non-split extensions
${\mathcal G}_1,\ldots,{\mathcal G}_m$ of ${\mathcal F}_p^{(j)}$ such that
$F_{p,j}\cong
\left({\mathcal G}_1\oplus\cdots\oplus{\mathcal G}_m\right)\otimes k(x_0)$
for the point $x_0$ of $M_C(r_p^{(j)},d_p^{(j)})$ corresponding to
$F_p^{(j)}$.
We put ${\mathcal F}_{p,j}:={\mathcal G}_1\oplus\cdots\oplus{\mathcal G}_m$
and
\[
 W:=\left.\left\{
 (x_p^{(j)})\in\prod_{1\leq p\leq s,1\leq j\leq r_p} M_C(r_p^{(j)},d_p^{(j)})
 \right| \bigotimes_{p,j}
 \det\left({\mathcal F}_{p,j}\otimes k(x_p^{(j)})\right)\cong L
 \right\}.
\]
Then $\dim W=-1+\sum_{p=1}^s r_p$.
We put ${\mathcal E}:=\bigoplus_{j,p}{\mathcal F}_{p,j}$.
Then parabolic structure can be parameterized by the product of 
flag schemes
\[
 U:=\prod_{i=1}^n \mathrm{Flag}({\mathcal E}|_{t_i\times W})
\]
over $W$.
The relative dimension of $U$ over $W$ is
$r(r-1)n/2$.
Let $\{\tilde{l}^{(i)}_j\}$ be the universal parabolic structure
on ${\mathcal E}_U$.
Note that
\[
 \End({\mathcal E}\otimes k(x))=\bigoplus_{p=1}^s\bigoplus_{j=1}^{r_p}
 \End({\mathcal F}_{p,j}\otimes k(x))
 \oplus\bigoplus_{p<q}\bigoplus_{j,k}
 \Hom({\mathcal F}_{p,j}\otimes k(x),{\mathcal F}_{q,k}\otimes k(x))
\]
for a point $x\in W$ and the group
\[
 G_x:=\left.\left\{
 \sum_{p,j}c_{pj}\mathrm{id}_{{\mathcal F}_{p,j}\otimes k(x)}
 +\sum_{p<q}\sum_{j,k}a^{(p,q)}_{j,k}
 \in \End({\mathcal E}\otimes k(x)) \right|
 \begin{array}{l}
 c_{pj}\in\C^*, \\
 a^{(p,q)}_{j,k}\in
 \Hom\left({\mathcal F}_{p,j}\otimes k(x),
 {\mathcal F}_{q,k}\otimes k(x)\right)
 \end{array}
 \right\}
\]
acts on the fiber $U_x$ of $U$ over $x$.
Note that we can take a section
$a\in\Hom(\cF_{p,j}\otimes k(x),\cF_{q,k}\otimes k(x))$
with $a(t_1)\neq 0$ or $a(t_2)\neq 0$ if $p<q$.
So we may assume that
$\alpha(t_2)\neq 0$ for some
$\alpha\in\Hom(\cF_{p_0,1}\otimes k(x),\cF_{p_0+1,1}\otimes k(x))$.
Modulo the action of the group
\[
 \left.\left\{
 \sum_{p,j}c_{pj}\mathrm{id}_{\cF_{p,j}}\otimes k(x)\right|
 c_{pj}\in\mathbf{C}^{\times}
 \right\},
\]
$\tilde{l}^{(1)}_{r-1}\otimes k(y)$ ($y\in U_x$)
can be parameterized by an algebraic scheme whose dimension is
at most
$r-\sum_{p=1}^s r_p$.
Moreover, modulo the action of the group
\[
 \left.\left\{ \mathrm{id}_{{\mathcal E}\otimes k(x)}+\sum_{p<q}\sum_{j,k}a^{(p,q)}_{j,k}
 \right| 
 a^{(p,q)}_{j,k}\in
 \Hom\left({\mathcal F}_{p,j}\otimes k(x),
 {\mathcal F}_{q,k}\otimes k(x)\right)\right\},
\]
$\tilde{l}^{(2)}_{r-1}\otimes k(y)$ ($y\in U_x$)
can be parameterized by an algebraic scheme whose dimension is at most
$r-2$.
Then, modulo the action of $G_x$,
$\tilde{l}^{(1)}_{r-1}\otimes k(y)$ and
$\tilde{l}^{(2)}_{r-1}\otimes k(y)$ ($y\in U_x$)
can be parameterized by an algebraic scheme
whose dimension is at most
$2r-\sum_{p=1}^s r_p-2$.
So we can take a finite number of subschemes
$Y_1,\ldots,Y_l$ of $U$ such that
$\dim Y_i\leq -1+\sum_{p=1}^s r_p +r(r-1)n/2-\sum_{p=1}^s r_p
=r(r-1)n/2-1$
and
$\bigcup_{i=1}^l G_x(Y_i)_x=U_x$
for any point $x\in W$.
We put $Y:=\coprod_{i=1}^l Y_i$.
Replacing $Y$ by its stratification, we may assume that the dimension of
\[
 F^0_y:=\left\{ s\in \End({\mathcal E}_y) \left|
 \text{$\Tr(s)=0$ and
 $s(t_i)(\tilde{l}^{(i)}_j)_y\subset(\tilde{l}^{(i)}_j)_y$
 for each $i,j$}
 \right\}\right.
\]
is constant for $y\in Y$.
There exists a scheme $Z$ over $Y$ such that
for any $Y$-scheme $S$, we have
\[
 Z(S)=\left\{
 \nabla:{\mathcal E}_S\rightarrow
 {\mathcal E}_S\otimes\Omega^1_C(t_1+\cdots+t_n)
 \left|
 \begin{array}{l}
 \text{$\nabla$ is a relative connection such that} \\
 \text{$(\res_{t_i\times S}(\nabla)-\lambda^{(i)}_j)
 ((\tilde{l}^{(i)}_j)_S)\subset (\tilde{l}^{(i)}_{j+1})_S$
 for any $i,j$} \\
 \text{and $\det(\nabla))=(\nabla_L)\otimes\cL$}
 \end{array}
 \right\}\right.
\]
We can see that each fiber of $Z\rightarrow Y$
over $y\in Y$ is an affine space isomorphic to
\[
 F^1_y=\left\{
 s\in\Hom({\mathcal E}_y,{\mathcal E}_y\otimes\Omega^1_C(t_1+\cdots+t_n))
 \left| \text{$\Tr(s)=0$ and
 $\res_{t_i}(s)(\tilde{l}^{(i)}_j)_y\subset (\tilde{l}^{(i)}_{j+1})_y$
for any $i,j$}\right\}\right..
\]
Let
$\tilde{\nabla}:{\mathcal E}_Z\rightarrow {\mathcal E}_Z
\otimes\Omega^1_C(t_1+\cdots+t_n)$
be the universal relative connection.
If we put
\[
 Z^s:=\left\{y\in Z \left|
 \text{$\sigma^{-1}\left(({\mathcal E},\tilde{\nabla},
 \{\tilde{l}^{(i)}_j\})\otimes k(y)\right)$ is $\balpha$-stable}
 \right\}\right.,
\]
then there is a natural morphism
$f:Z^s\rightarrow M^{\balpha}_C(\bt,\blambda,L)$.
For each $z\in Z$ over $y\in Y$, consider the homomorphism
\[
 \nabla^1_z:F^0_y\otimes k(z)\lra F^1_y\otimes k(z);
 \quad \nabla^1_z(s):=\tilde{\nabla}_z\circ s-s\circ\tilde{\nabla}_z.
\]
Then the group $\mathrm{Stab}_y$ of automorphisms $g$
of the parabolic bundles
$({\mathcal E}_y,\{(\tilde{l}^{(i)}_j)_y\})$
which induce the identity on $\det({\mathcal E}_y)$
acts on the fiber $Z_y$ and the tangent map of
the morphism $\mathrm{Stab}_y\ni g\mapsto gz\in Z^s$ at $\mathrm{id}$
is just the homomorphism $\nabla^1_z$ which is injective.
Replacing $Z^s$ by a finite disjoint union of subschemes,
we may assume that $Z^s\stackrel{f}\rightarrow \overline{f(Z^s)}$ is flat.
Then we have
\begin{align*}
 \dim_{f(z)}\overline{f(Z^s)}&=\dim_zZ^s-\dim_zf^{-1}(f(z)) \\
 &\leq \dim_zZ^s-\dim_{\mathrm{id}}\mathrm{Stab}_y \\
 &=\dim_zZ^s-\dim F^0_y \\
 &\leq \dim_yY+\dim_zZ^s_y-\dim F^0_y \\
 &\leq\dim_yY+\dim F^1_y-\dim F^0_y \\
 &\leq \dim Y+r(r-1)n/2\leq r(r-1)n-1.
\end{align*}
Here $z\in Z^s$ and $Z^s\ni z\mapsto y\in Y$.
Since $M^{\balpha}_C(\bt,\blambda,L)\setminus N$
is a finite union of such $f(Z^s)$'s,
we have $\dim(M^{\balpha}_C(\bt,\blambda,L)\setminus N)<r(r-1)n$.
\end{proof}

\begin{Proposition}\label{irred-g=0}
 Assume that $rn-2(r+1)>0$ and $r\geq 2$.
 Then $M^{\balpha}_{\BP^1}(\bt,\blambda,L)$ is an irreducible variety
 of dimension $(r-1)(rn-2r-2)$.
\end{Proposition}

\begin{proof}
First we will show that the Zariski open set
\[
 N:=\left\{ \left. (E,\nabla,\{l^{(i)}_j\})\in
 M^{\balpha}_{\BP^1}(\bt,\blambda,L) \right|
 \dim \End(E,\{l^{(i)}_j\})=1
 \right\}
\]
is irreducible.
There exists an integer $n_0$ such that
$h^1(E(n_0))=0$ for any
$(E,\nabla,\{l^{(i)}_j\})\in
M^{\balpha}_{\BP^1}(\bt,\blambda,L)$.
Then we can easily construct a smooth irreducible variety
$X$ and a flat family of parabolic bundles
$(\tilde{E},\{\tilde{l}^{(i)}_j\})$ on $\BP^1\times X$ over $X$
whose geometric fibers is just the set of all the
simple parabolic bundles $(E,\{l^{(i)}_j\})$ satisfying $h^1(E(n_0))=0$
and $\det E=L$.
We can take an isomorphism
$\det\tilde{E}\cong L\otimes{\mathcal L}$
for some line bundle ${\mathcal L}$ on $X$.
The functor $\mathrm{Conn}_X: (Sch/X)\lra(Sets)$ defined by
\[
 \mathrm{Conn}_X(U)=\left\{
 \nabla: \tilde{E}_U\ra\tilde{E}_U\otimes\Omega^1_{\BP^1}(D(\bt))
 \left|
 \begin{array}{l}
  \text{$\nabla$ is a relative connection satisfying} \\
  \text{$(\res_{t_i}(\nabla)-\lambda^{(i)}_j)((\tilde{l}^{(i)}_j)_U)
  \subset(\tilde{l}^{(i)}_{j+1})_U$ for any $i,j$} \\
  \text{and $\det(\nabla)=(\nabla_L)_U\otimes{\mathcal L}$}
 \end{array}
 \right\}\right.
\]
is represented by a scheme $Y$ of finite type over $X$.
We can see by the same proof as Proposition \ref{irred-generic} that
$Y$ is smooth over $X$ whose fiber is an affine space of dimension
$h^0((\cF_0^1)_x)$ for $x\in X$,
where
\[
 \cF^1_0=\left\{
 s\in {\mathcal End}(\tilde{E})\otimes\Omega^1_{\BP^1}(D(\bt))
 \left|
 \text{$\Tr(s)=0$ and
 $\res_{t_i\times X}(s)(\tilde{l}^{(i)}_j)\subset\tilde{l}^{(i)}_{j+1}$
 for any $i,j$}
 \right\}\right..
\]
So $Y$ is irreducible.
Let
\[
 \tilde{\nabla}:\tilde{E}_Y\longrightarrow
 \tilde{E}_Y\otimes\Omega^1_{\BP^1}(t_1+\cdots+t_n)
\]
be the universal relative connection and put
\[
 Y':=\left\{ y\in Y \left|
 \text{$(\tilde{E},\tilde{\nabla},\{\tilde{l}^{(i)}_j\})\otimes k(y)$
 is $\balpha$-stable}
 \right\}\right..
\]
Then $Y'$ a Zariski open subset of $Y$ and it is irreducible.
A morphism $Y'\ra M^{\balpha}_{\BP^1}(\bt,\blambda,L)$
is induced whose image is just $N$.
Thus $N$ becomes irreducible.

In order to prove that $M^{\balpha}_{\BP^1}(\bt,\blambda,L)$ is irreducible,
it suffices to show that
$\dim (M^{\balpha}_{\BP^1}(\bt,\blambda,L)\setminus N)<(r-1)(rn-2r-2)$
because $M^{\balpha}_{\BP^1}(\bt,\blambda,L)$ is smooth of
equidimension $(r-1)(rn-2r-2)$.
For this, we will construct an algebraic scheme $W'$ of dimension less than
$(r-1)(rn-2r-2)$ and a surjection
$f:W'\rightarrow M^{\balpha}_{\BP^1}(\bt,\blambda,L)\setminus N$.

Take any parabolic bundle
$(E,\{l^{(i)}_j\})$ such that $\rank E=r$, $\deg E=d$,
$\dim\End(E,\{l^{(i)}_j\})\geq 2$ and $h^1(E(n_0))=0$.
We will show that such parabolic bundles can be parameterized by
an algebraic scheme of dimension less than $(r-1)(rn/2-r-1)$.
We can write
\[
 E \cong \cO_{\BP^1}(a_1)^{\oplus r_1}\oplus\cO_{\BP^1}(a_2)^{\oplus r_2}
 \oplus\cdots\oplus \cO_{\BP^1}(a_s)^{\oplus r_s}
\]
with $a_1<a_2<\cdots<a_s$.
If
$l^{(k)}_{r-1}\not\subset
\cO_{\BP^1}(a_2)^{\oplus r_2}|_{t_k}\oplus\cdots\oplus
\cO_{\BP^1}(a_s)^{\oplus r_s}|_{t_k}$ for some $k$,
then we replace $(E,\{l^{(i)}_j\})$ by its elementary transform along
$t_k$ by $l^{(k)}_{r-1}$, which is isomorphic to the bundle
\[
 \cO_{\BP^1}(a_1-1)^{\oplus r_1-1} \oplus \cO_{\BP^1}(a_1)
 \oplus \cO_{\BP^1}(a_2-1)^{\oplus r_2}
 \oplus\cdots\oplus \cO_{\BP^1}(a_s-1)^{\oplus r_s}
\]
with a certain parabolic structure.
Repeating this process, we finally obtain two cases:
\[
 \begin{cases}
 E\cong \cO_{\BP^1}(a)^{\oplus r} \\
 l^{(i)}_{r-1}\subset\cO_{\BP^1}(a_2)^{\oplus r_2}|_{t_i}\oplus\cdots\oplus
 \cO_{\BP^1}(a_s)^{\oplus r_s}|_{t_i}
 \quad \text{for all $i$ and $r_1>0$}.
 \end{cases}
\]

\noindent
{\bf Case 1.}
$E\cong\cO_{\BP^1}(a)^{\oplus r}$. \\
Tensoring $\cO_{\BP^1}(-a)$, we may assume $E\cong\cO_{\BP^1}^{\oplus r}$.
For a suitable choice of a basis of $\cO_{\BP^1}^{\oplus r}$,
we may assume that
\[
 l^{(1)}_{r-i}=k(t_1)e_1+k(t_1)e_2+\cdots+k(t_1)e_i
\]
for $i=1,\ldots,r$, where $e_j$ is the vector of size $r$
whose $j$-th component is $1$ and others are zero.
Then the group of automorphisms of $E$ fixing
$l^{(1)}_{r-1},l^{(1)}_{r-2},\ldots,l^{(1)}_1$ is
\[
 B= \left \{ (a_{ij})\in GL_r(\C) \left|
 \text{$a_{ij}=0$ for $i>j$} \right\}\right..
\]
We put
\[
 p(1):=\max \left\{ p \left|
 \text{$v_p\neq 0$ for some $(v_1,\ldots,v_r)\in l^{(2)}_{r-1}$}
 \right\}\right..
\]
Applying a certain automorphism in $B$, we may assume that
$l^{(2)}_{r-1}=k(t_2)e_{p(1)}$.
Inductively we put
\[
 p(i):= \max \left\{ p \left|
 \text{$v_p\neq 0$ for some $(v_1,\ldots,v_r)\in l^{(2)}_{r-i}$
 and $p\neq p(j)$ for any $j<i$}
 \right\}\right.
\]
for $i=1,\ldots,r$.
Applying an automorphism in $B$ which fixes
$l^{(2)}_{r-1},\ldots,l^{(2)}_{r-i+1}$, we may assume that
\[
 l^{(2)}_{r-i}=l^{(2)}_{r-i+1}+k(t_2)e_{p(i)}.
\]
Then the group of automorphisms of $E$
fixing both $\{l^{(1)}_j\}$ and $\{l^{(2)}_j\}$ is
\[
 B'=\left\{ (a_{ij})\in GL_r(\C) \left|
 \text{$a_{ij}=0$ for $i>j$ and $a_{p(i)p(j)}=0$ for $i>j$}
 \right\}\right..
\]
Applying a certain automorphism in $B'$, we may assume that
$l^{(3)}_{r-1}$ can be generated by a vector $w=(w_1,\ldots,w_r)$,
where each $w_i$ is either $1$ or $0$.
Note that we have $n\geq 3$ by the assumption of the proposition.
The group of automorphisms of $E$ fixing $\{l^{(1)}_j\}$,
$\{l^{(2)}_j\}$ and $l^{(3)}_{r-1}$ is
\[
 B''=\left\{ (a_{ij})\in GL_r(\C) \left|
 \begin{array}{l}
  \text{$a_{ij}=0$ for $i>j$, $a_{p(i)p(j)}=0$ for $i>j$} \\
  \text{and there is $c\in\C^{\times}$ satisfying
  $a_{ii}w_i+\sum_{k\neq i} a_{ik}w_k=cw_i$ for any $i$}
 \end{array}
 \right\}\right..
\]
Since the parabolic bundle $(E,\{l^{(i)}_j\})$
has a nontrivial endomorphism,
there is $(i,j)$ with $i<j$ and $p(i)<p(j)$ or
there is some $i$ satisfying $w_i=0$.
Note that $r\geq 3$ or $n\geq 4$ because we assume that
$rn-2r-2>0$.
First assume that $r\geq 3$.
If $w_{p(i)}=w_{p(j)}=1$ and
$p(i)<p(j)$ for some $i<j$, then we apply an automorphism of the form
$(a_{kl})$ such that $a_{kk}=1$ if $k\neq p(i)$, $a_{kl}=0$ if $k\neq l$
and $(k,l)\neq(p(i),p(j))$, $a_{p(i)p(i)}=1-c$ and $a_{p(i)p(j)}=c$
for some $c\in\C\setminus\{1\}$.
Then we can parameterize $l^{(3)}_{r-2}$ modulo this action
by an algebraic scheme whose dimension is less than $r-2$.
If $w_i=0$ for some $i$, then we apply an automorphism of the form
$(a_{kl})$ such that $a_{kk}=1$ for $k\neq i$, $a_{ii}=c$ ($c\in\C^{\times}$)
and $a_{kl}=0$ for $k\neq l$.
Then $l^{(3)}_{r-2}$ can be parameterized modulo this action
by an algebraic scheme of dimension less than $r-2$.
Next assume that $n\geq 4$.
If $w_{p(i)}=w_{p(j)}=1$ and
$p(i)<p(j)$ for some $i<j$, then we apply an automorphism of the form
$(a_{kl})$ such that $a_{kk}=1$ if $k\neq p(i)$, $a_{kl}=0$ if $k\neq l$
and $(k,l)\neq(p(i),p(j))$, $a_{p(i)p(i)}=1-c$ and $a_{p(i)p(j)}=c$
for some $c\in\C\setminus\{1\}$.
Then we can parameterize $l^{(4)}_{r-1}$ modulo this action
by an algebraic scheme whose dimension is less than $r-1$.
If $w_i=0$ for some $i$, then we apply an automorphism of the form
$(a_{kl})$ such that $a_{kk}=1$ for $k\neq i$, $a_{ii}=c$ ($c\in\C^{\times}$)
and $a_{kl}=0$ for $k\neq l$.
Then $l^{(4)}_{r-1}$ can be parameterized modulo this action
by an algebraic scheme of dimension less than $r-1$.
Taking account of all the above arguments,
we can parameterize parabolic bundles $(E,\{l^{(i)}_j\})$
satisfying $E\cong\cO_{\BP^1}^{\oplus r}$ and
$\dim\End(E,\{l^{(i)}_j\})>1$ by an algebraic scheme
of dimension less than $(r-1)(rn/2-(r+1))$.

\noindent
{\bf Case 2.}
$E\cong\cO_{\BP^1}(a_1)^{\oplus r_1}\oplus\cdots\oplus
\cO_{\BP^1}(a_s)^{\oplus r_s}$ with $r_1>0$ and
$l^{(i)}_{r-1} \subset
\cO_{\BP^1}(a_2)^{\oplus r_2}|_{t_i}\oplus\cdots\oplus
\cO_{\BP^1}(a_s)^{\oplus r_s}|_{t_i}$. \\
For each $p$, we consider an identification
$\cO_{\BP^1}(a_p)^{\oplus r_p}|_{t_i}=\C^{r_p}$ and
denote by $e^{(i)}_{p,j}$ the element of this vector space
whose $j$-th component is $1$ and others are zero.
We put
\[
 p(1):=\min\left\{ p \left|
 \text{$v_p\neq 0$ for some $(v_1,\ldots,v_s)\in l^{(1)}_{r-1}$
 where $v_q\in\cO_{\BP^1}(a_q)^{\oplus r_q}|_{t_1}$ for any $q$}
 \right\}\right..
\]
Applying an automorphism of $E$, we may assume that
$l^{(1)}_{r-1}$ is generated by $e^{(1)}_{p(1),j(1)}$
with $j(1)=1$.
Inductively we put
\[
 p(i):=\min\left\{ p \left|
 \text{$v^{(p)}_j\neq 0$ for some $\sum_{q,k} v^{(q)}_k e^{(1)}_{q,k}\in l^{(1)}_{r-i}$
 where $(p,j)\neq(p(k),j(k))$ for any $k<i$}
 \right\}\right.
\]
for $i=1,\ldots,r$ and we put
\[
 j(i):=1+\max\left(\{0\}\cup\{j(i')|\text{$i'<i$ and $p(i')=p(i)$}\}\right).
\]
Applying an automorphism of $E$ fixing $l^{(1)}_{r-1},\ldots,l^{(1)}_{r-i+1}$,
we may assume that $l^{(1)}_{r-i}=l^{(1)}_{r-i+1}+k(t_1)e^{(1)}_{p(i),j(i)}$.
Then each $l^{(1)}_{r-i}$ is generated by
$e^{(1)}_{p(1),j(1)},\ldots,e^{(1)}_{p(i),j(i)}$ and
the group of automorphisms of $E$ fixing $\{l^{(1)}_j\}$ is
\[
 B=\left\{ (a^{pq}_{jk})^{1\leq p,q \leq s}
 _{1\leq j\leq r_p, 1\leq k\leq r_q}
 \left|
 \begin{array}{l}
  \text{$(a^{pq}_{jk})_{1\leq j\leq r_p, 1\leq k\leq r_q}\in
  \End(\cO(a_q)^{r_q},\cO(a_p)^{r_p})$ for each $(p,q)$} \\
  \text{$(a^{(pp)}_{jk})\in\Aut(\cO(a_p)^{r_p})$ and
  $a^{p(i)p(i')}_{j(i)j(i')}(t_1)=0$ for $i>i'$}
 \end{array}
 \right\}\right..
\]
Note that we have $j(i')<j(i)$ if $i'<i$ and $p(i')=p(i)$.
We will also normalize $\{l^{(2)}_j\}$.
First we put
\[
 (p'(1),-j'(1)):=\min\left\{ (p,-j) \left|
 \text{$v^{(p)}_j\neq 0$ for some $\sum_{q,k} v^{(q)}_k e^{(2)}_{q,k}\in l^{(2)}_{r-1}$}
 \right\}\right.,
\]
where we consider the lexicographic order on the pair $(p,-j)$.
Applying an automorphism of $E$ fixing $\{l^{(1)}_j\}$,
we may assume that $l^{(2)}_{r-1}$ is generated by
$e^{(2)}_{p'(1),j'(1)}$.
Inductively we put
\[
 (p'(i),-j'(i)):=\min\left\{ (p,-j) \left|
 \begin{array}{l}
  \text{$v^{(p)}_j\neq 0$ for some
  $\sum_{q,k} v^{(q)}_k e^{(2)}_{q,k}\in l^{(2)}_{r-i}$ and} \\
  \text{$(p,j)\neq (p'(k),j'(k))$ for any $k<i$}
 \end{array}
 \right\}\right.
\]
for $i=1,\ldots,r$.
Applying an automorphism in $B$ fixing
$l^{(2)}_{r-1},\ldots,l^{(2)}_{r-i+1}$, we may assume that
$l^{(2)}_{r-i}=l^{(2)}_{r-i+1}+k(t_2)e^{(2)}_{p'(i),j'(i)}$.
Then the group of automorphisms of $E$ fixing $\{l^{(1)}_j\}$
and $\{l^{(2)}_j\}$ is
\[
 B'=\left\{ (a^{pq}_{jk})^{1\leq p,q \leq s}
 _{1\leq j\leq r_p, 1\leq k\leq r_q}\in B
 \left|
 \text{$a^{p'(i)p'(l)}_{j'(i)j'(l)}(t_2)=0$ for $i>l$}
 \right\}\right..
\]
Applying a certain automorphism in $B'$, we may assume that
$l^{(3)}_{r-1}$ can be generated by a vector
$w=\sum w^{(p)}_j e^{(3)}_{p,j}$,
where each $w_j^{(p)}$ is either $1$ or $0$.
Then the group of automorphisms of $E$ fixing $\{l^{(1)}_j\}$,
$\{l^{(2)}_j\}$ and $l^{(3)}_{r-1}$ is
\[
 B''=\left\{ (a^{pq}_{jk})^{1\leq p,q \leq s}
 _{1\leq j\leq r_p, 1\leq k\leq r_q}\in B'
 \left|
 \begin{array}{l}
  \text{for some $c\in\C^{\times}$,
  $\sum_{q,k} a^{pq}_{jk}w^{(q)}_k=cw^{(p)}_j$ for any $p,j$}
 \end{array}
 \right\}\right..
\]
Note that $p(1)\geq 2$ because
$l^{(1)}_{r-1}\subset \cO_{\BP^1}(a_2)^{\oplus r_2}|_{t_1}
\oplus\cdots\oplus \cO_{\BP^1}(a_s)^{\oplus r_s}|_{t_1}$.
If $w^{(1)}_1\neq 0$ and $w^{(p(1))}_{j(1)}\neq 0$,
then there exist elements
$(a^{pq}_{jk})$ of $B''$ satisfying
$a^{p(1)1}_{j(1)1}(t_3)w^{(1)}_{1}+a^{p(1)p(1)}_{j(1)j(1)}w^{(p(1))}_{j(1)}
=w^{(p(1))}_{j(1)}$,
$a^{pp}_{jj}=1$ for $(p,j)\neq(p(1),j(1))$
and $a^{(pq)}_{jk}=0$ for $(p,j)\neq(q,k)$ and $(p,j,q,k)\neq (p(1),j(1),1,1)$.
Applying such automorphisms,
$l^{(3)}_{r-2}$ can be parameterized modulo this action by an algebraic scheme
of dimension less than $r-2$ when $r\geq 3$
and $l^{(4)}_{r-1}$ can be parameterized modulo this action
by an algebraic scheme of dimension less than $r-1$ when $n\geq 4$.
If $w^{(p_0)}_{j_0}=0$ for some $(p_0,j_0)$, then we apply an automorphism
$(a^{pq}_{jk})$ in $B''$ satisfying $a^{pp}_{jj}=1$ for $(p,j)\neq(p_0,j_0)$,
$a^{p_0p_0}_{j_0j_0}=c$ ($c\in\C^{\times}$) and
$a^{pq}_{jk}=0$ for $(p,j)\neq(q,k)$.
Then $l^{(3)}_{r-2}$ can be parameterized modulo this action
by an algebraic scheme of dimension less than $r-2$ when $r\geq 3$ and
$l^{(4)}_{r-1}$ can be parameterized modulo this action by
an algebraic scheme of dimension less than $r-1$ when $n\geq 4$.
Therefore we can parameterize parabolic bundles $(E,\{l^{(i)}_j\})$
satisfying
$E\cong\cO_{\BP^1}(a_1)^{\oplus r_1}
\oplus\cdots\oplus \cO_{\BP^1}(a_s)^{\oplus r_s}$
and
$l^{(i)}_{r-1}\subset \cO_{\BP^1}(a_2)^{\oplus r_2}|_{t_i}
\oplus\cdots\oplus \cO_{\BP^1}(a_s)^{\oplus r_s}|_{t_i}$
for any $i$
by an algebraic scheme of dimension less than
$(r-1)(rn/2-r-1)$.

By Case 1 and Case 2, we can construct an algebraic scheme
$Z$ of dimension less than $(r-1)(rn/2-r-1)$ and a flat family
$(\tilde{E},\{\tilde{l}^{(i)}_j\})$ of parabolic bundles on
$\BP^1_Z$ over $Z$ such that the set of the geometric fibers of
$(\tilde{E},\{\tilde{l}^{(i)}_j\})$ is just the set of parabolic
bundles $(E,\{l^{(i)}_j\})$ satisfying 
$\rank E=r$, $\deg E=d$,
$\dim\End(E,\{l^{(i)}_j\})\geq 2$
and $h^1(E(n_0))=0$.
Note that the functor
$\mathrm{Conn}_Z: (Sch/Z) \lra (Sets)$ defined by
\[
 \mathrm{Conn}_Z(S)=\left\{
 \nabla:\tilde{E}_S \ra \tilde{E}_S\otimes\Omega^1_{\BP^1}(t_1+\cdots+t_n)
 \left|
 \begin{array}{l}
 \text{$\nabla$ is a relative connection satisfying} \\
 \text{$(\res_{t_i\times S}(\nabla)-\lambda^{(i)}_j)(\tilde{l}^{(i)}_j)_S
 \subset(\tilde{l}^{(i)}_{j+1})_S$ for any $i,j$} \\
 \text{and $\det(\nabla)=(\nabla_L)_S\otimes{\mathcal L}$}
 \end{array}
 \right\}\right.
\]
can be represented by a scheme $W$ of finite type over $Z$.
Let $\tilde{\nabla}:\tilde{E}_W \ra 
\tilde{E}_W\otimes\Omega^1_{\BP^1}(t_1+\cdots+t_n)$
be the universal relative connection.
If we put
\begin{gather*}
 \cF^0_0:=\left.\left\{ s\in{\mathcal End}(\tilde{E}) \right|
 \text{$\Tr(s)=0$ and
 $s(\tilde{l}^{(i)}_j)\subset\tilde{l}^{(i)}_j$
 for any $i,j$} \right\}, \\
 \cF^1_0:=\left.\left\{ s\in{\mathcal End}(\tilde{E})
 \otimes\Omega^1_{\BP^1}(t_1+\cdots+t_n) \right|
 \text{$\Tr(s)=0$ and
 $\res_{t_i\times Z}(s)(\tilde{l}^{(i)}_j)\subset\tilde{l}^{(i)}_{j+1}$
 for any $i,j$} \right\}, \\
 \nabla^{\dag}:\cF^0_0\longrightarrow \cF^1_0;
 \quad \nabla^{\dag}(s)=\tilde{\nabla}\circ s-s\circ\tilde{\nabla},
\end{gather*}
then each fiber $W_z$ of $W\rightarrow Z$ over $z$ is an affine space
isomorphic to
$H^0(\cF^1_0\otimes k(z))$.
If we put
\[
 W':=\left\{ p\in W \left|
 \text{$(\tilde{E},\tilde{\nabla},\{\tilde{l}^{(i)}_j\})\otimes k(p)$
 is $\balpha$-stable}
 \right\}\right.,
\]
then a morphism
$f:W' \ra M^{\balpha}_{\BP^1}(\bt,\blambda,L)$
is induced.
Replacing $W'$ by a finite disjoint union of subschemes,
we may assume that $W'\stackrel{f}\longrightarrow \overline{f(W')}$ is flat.
There is a canonical action of
\[
 \Aut_0(\tilde{E}\otimes k(z),\{\tilde{l}^{(i)}_j\otimes k(z)\})=
 \left\{
 g:\tilde{E}\otimes k(z)\stackrel{\sim}\longrightarrow\tilde{E}\otimes k(z)
 \left|
 \begin{array}{l}
  \text{$g(\tilde{l}^{(i)}_j\otimes k(z))\subset\tilde{l}^{(i)}_j\otimes k(z)$ for any $i,j$} \\
  \text{and $\wedge^rg=\mathrm{id}$}
 \end{array}
 \right.\right\}
\]
on the fiber $W'_z$ of $W'\rightarrow Z$ over $z\in Z$
and we have
$f(gp)=p$ for any $p\in W'_z$ and $g\in\Aut_0(E,\{l^{(i)}_j\})$.
Take $p\in W'$ lying over $z\in Z$.
Then the tangent map of
$\Aut_0(\tilde{E}\otimes k(p),\{\tilde{l}^{(i)}_j\otimes k(p)\})
\ni g \mapsto gp\in W'_z$
at $\mathrm{id}$
is the linear map $H^0(\nabla^{\dag}\otimes k(p))$,
which is injective because
$(\tilde{E}\otimes k(p),\tilde{\nabla}\otimes k(p),\{\tilde{l}^{(i)}_j\otimes k(p)\})$
is $\balpha$-stable.
So we have for any $p\in W'$,
\begin{align*}
 \dim_{f(p)}\overline{f(W')}&=\dim_pW'-\dim_pf^{-1}(f(p)) \\
 &\leq \dim_pW'-
 \dim_{\mathrm{id}}\Aut_0(\tilde{E}\otimes k(p),\{\tilde{l}^{(i)}_j\otimes k(p)\}) \\
 &=\dim_pW'-\dim H^0(\cF^0_0\otimes k(p)) \\
 &\leq \dim_zZ+\dim_pW'_z-\dim H^0(\cF^0_0\otimes k(p)) \\
 &\leq\dim_zZ+\dim H^0(\cF^1_0\otimes k(p))-\dim H^0(\cF^0_0\otimes k(p)) \\
 &=\dim_zZ+(r-1)(rn/2-r-1) \\
 &<(r-1)(rn-2(r+1)).
\end{align*}
Here $W'\ni p\mapsto z\in Z$.
Since $M^{\balpha}_{\BP^1}(\bt,\blambda,L)\setminus N$ can be covered by
such $f(W')$'s, we have
\[
 \dim\left(M^{\balpha}_{\BP^1}(\bt,\blambda,L)\setminus N\right)
 <(r-1)(rn-2(r+1)).
\]
Hence $M^{\balpha}_{\BP^1}(\bt,\blambda,L)$ is irreducible.
\end{proof}

\section{Riemann-Hilbert correspondence}
\label{riemann-hilbert}

We take a point $x\in T$ and
denote $\cC_x$ and $\tilde{\bt}_x$ simply by $C$ and $\bt$.
We also denote the fiber of $M_{\cC/T}^{\balpha}(\tilde{\bt},r,d)$
over $x\in T$ by $M_C^{\balpha}(\bt,r,d)$.
Sending $(E,\nabla,\{l^{(i)}_*\}_{1\leq i\leq n})$
to $\ker\nabla^{an}|_{C\setminus\{t_1,\ldots,t_n\}}$,
we obtain a holomorphic mapping
\[
 \RH_x: M^{\balpha}_C(\bt,r,d)\lra
 RP_r(C,\bt),
\]
which makes the diagram
\[
 \begin{CD}
  M^{\balpha}_C(\bt,r,d) @> \RH_x >> RP_r(C,\bt) \\
  @VVV  @VVV \\
  \Lambda^{(n)}_r(d) @> rh >> \cA^{(n)}_r
 \end{CD}
\]
commute.
Here $rh:\Lambda^{(n)}_r(d)\ni\blambda=(\lambda^{(i)}_j)
\mapsto\ba=(a^{(i)}_j)\in\cA^{(n)}_r$
is defined by
\[
 \prod_{j=0}^{r-1}(X-\exp(-2\pi\sqrt{-1}\lambda^{(i)}_j))
 =X^r+a^{(i)}_{r-1}X^{r-1}+\cdots+a^{(i)}_1X+a^{(i)}_0
\]
for $i=1,\ldots,n$.
If we denote the fiber of $M_C^{\balpha}(\bt,r,d)$ over
$\blambda\in\Lambda^{(n)}_r(d)$ by
$M_C^{\balpha}(\bt,\blambda)$, then
$\RH_x$ induces a holomorphic mapping
\[
 \RH_{(x,\blambda)}: M_C^{\balpha}(\bt,\blambda)
 \longrightarrow
 RP_r(C,\bt)_{\ba}
\]
for $rh(\blambda)=\ba$.
In this section we will prove the properness of
this morphism $\RH_{(x,\blambda)}$ and obtain Theorem \ref{main-thm}.
In order to prove the properness of $\RH_{(x,\blambda)}$,
it is essential to prove the surjectivity of $\RH_{(x,\blambda)}$
and compactness of every fiber of $\RH_{(x,\blambda)}$.
For the proof of the surjectivity of $\RH_{(x,\blambda)}$,
it is essential to use the following Langton's type theorem:

\begin{Lemma}\label{langton}
Let $R$ be a discrete valuation ring with residue field $k=R/m$
and quotient field $K$.
Take a flat family $(E,\nabla,\{l^{(i)}_j\})$ of parabolic connections
on $C\times\Spec R$ over $\Spec R$ whose generic fiber
$(E,\nabla,\{l^{(i)}_j\})\otimes K$ is $\balpha$-semistable.
Then we can obtain, by repeating elementary transformations
of $(E,\nabla,\{l^{(i)}_j\})$ along the special fiber $C\times\Spec k$,
a flat family
$(E',\nabla',\{{l'}^{(i)}_j\})$ of parabolic connections
on $C\times\Spec R$ over $\Spec R$ such that
$(E',\nabla',\{{l'}^{(i)}_j\})\otimes K\cong
(E,\nabla,\{l^{(i)}_j\})\otimes K$ and
$(E',\nabla',\{{l'}^{(i)}_j\})\otimes k$ is
$\balpha$-semistable.
\end{Lemma}

\begin{proof}
This is just a corollary of the proof of
[\cite{IIS-1}, Proposition 5.5].
\end{proof}

\begin{Proposition}\label{surjectivity}
 Assume that $n\geq 3$ if $g=0$ and that $n\geq 1$ if $g\geq 1$.
 Then the morphism
 $\RH_{(x,\blambda)}:M^{\balpha}_C(\bt,\blambda)\ra RP_r(C,\bt)_{\ba}$
 is surjective for any $\blambda\in\Lambda^{(n)}_r(d)$
 satisfying $rh(\blambda)=\ba$.
\end{Proposition}

\begin{proof}
By Proposition \ref{prop:elementary-transform},
we obtain an isomorphism
\begin{equation}\label{sigma}
 \sigma: \cM_C(\bt,\blambda)\stackrel{\sim}\lra
 \cM_C(\bt,\bmu)
\end{equation}
of functors,
where $0\leq \mathrm{Re}(\mu^{(i)}_j)<1$
for any $i,j$.
Note that $\sigma$ induces an isomorphism
\[
 M^{\mathrm{irr}}_C(\bt,\blambda)\xrightarrow[\sim]{\sigma}
 M^{\mathrm{irr}}_C(\bt,\bmu),
\]
where $M^{\mathrm{irr}}_C(\bt,\blambda)$
(resp.\ $M^{\mathrm{irr}}_C(\bt,\bmu)$) is the moduli space
of irreducible $(\bt,\blambda)$-parabolic connections
(resp.\ $(\bt,\bmu)$-parabolic connections).
By construction, we have
$\RH_{(x,\bmu)}\circ\sigma|_{M^{\mathrm{irr}}_C(\bt,\blambda)}
=\RH_{(x,\blambda)}|_{M^{\mathrm{irr}}_C(\bt,\blambda)}$.

Consider the restriction
\[
 \RH_{(x,\blambda)}: M^{\mathrm{irr}}_C(\bt,\blambda)
 \lra RP_r(C,\bt)^{\mathrm{irr}}_{\ba},
\]
of $\RH$, where $M^{\mathrm{irr}}_C(\bt,\blambda)$
is the moduli space of irreducible parabolic connections and
$RP_r(C,\bt)^{\mathrm{irr}}_{\ba}$ is the moduli space of
irreducible representations.
For any point $\rho\in RP_r(C,\bt)_{\ba}^{\rm irr}$,
$\rho$ corresponds to a holomorphic vector bundle
$E^{\circ}$ on $C\setminus\{t_1,\ldots,t_n\}$ and
a holomorphic connection
$\nabla^{\circ}:E^{\circ}\rightarrow
E^{\circ}\otimes\Omega^1_{C\setminus\{t_1,\ldots,t_n\}}$.
By [\cite{Deligne}, II, Proposition 5.4],
there is a unique pair $(E,\nabla)$ of
a holomorphic vector bundle $E$ on $C$
and a logarithmic connection
$\nabla:E\rightarrow E\otimes\Omega^1_C(t_1+\cdots+t_n)$
such that
$(E,\nabla)|_{C\setminus\{t_1,\ldots,t_n\}}\cong
(E^{\circ},\nabla^{\circ})$
and that all the eigenvalues of
$\res_{t_i}(\nabla)$ lies in $\{z\in \C | 0\leq\mathrm{Re}(z)<1\}$.
We can take a parabolic structure $\{l^{(i)}_j\}$ on $E$
which is compatible with the connection $\nabla$.
Then $(E,\nabla,\{l^{(i)}_j\})$ becomes a
$(\bt,\bmu)$-parabolic connection such that
$\RH_{(x,\bmu)}(E,\nabla,\{l^{(i)}_j\})=\rho$.
Thus we have
$\sigma^{-1}(E,\nabla,\{l^{(i)}_j\})\in M_C^{\rm irr}(\bt,\blambda)
\subset M_C^{\balpha}(\bt,\blambda)$
(note that an irreducible parabolic connection is automatically
$\balpha$-stable)
and $\RH_{(x,\blambda)}(\sigma^{-1}(E,\nabla,\{l^{(i)}_j\}))=\rho$.

So it is sufficient to show that any member
$[\rho]\in RP_r(C,\bt)_{\ba}$ corresponding to a reducible
representation is contained in the image of $\RH_{(x,\blambda)}$,
in order to prove the surjectivity of $\RH_{(x,\blambda)}$.
The fundamental group $\pi_1(C\setminus\{t_1,\ldots,t_n\},*)$
is generated by cycles $\alpha_1,\beta_1,\ldots,\alpha_g,\beta_g$
and loops $\gamma_i$ around $t_i$ for $i=1,\ldots,n$
whose relation is given by
\[
 \prod_{i=1}^g \alpha_i\beta_i\alpha_i^{-1}\beta_i^{-1}
 \prod_{j=1}^n \gamma_j=1.
\]
Since $n>0$, the fundamental group is isomorphic to the free group
generated by the free generators
$\alpha_1,\beta_1,\ldots,\alpha_g,\beta_g,\gamma_1,\ldots,\gamma_{n-1}$.
So the space $\Hom(\pi_1(C\setminus\{t_1,\ldots,t_n\},*),GL_r(\C))$
of representations is isomorphic to $GL_r(\C)^{2g+n-1}$
and it is irreducible.
By Proposition \ref{openness-of-irreducibility} and
Remark \ref{existence-of-irreducible-representation},
the subset of the points of $GL_r(\mathbf{C})^{2g+n-1}$ corresponding to
the irreducible representations is a non-empty Zariski open subset.
Thus there exist a smooth affine curve $T$ and a morphism
$f:T\ra\Hom(\pi_1(C\setminus\{t_1,\ldots,t_n\},*),GL_r(\C))$ such that
the image $f(p_0)$ of a special point $p_0\in T$ corresponds to the
representation $\rho$ and the image $f(\eta)$ of 
$\eta\in T\setminus\{p_0\}$ corresponds to an irreducible representation.
Since the morphism $f$ corresponds to a family of
representations of the fundamental group
$\pi_1(C\setminus\{t_1,\ldots,t_n\},*)$,
we can construct a holomorphic vector bundle
$\tilde{E}^{\circ}$ on $(C\setminus\{t_1,\ldots,t_n\})\times T$
and a holomorphic connection
$\tilde{\nabla}^{\circ}:\tilde{E}^{\circ}\rightarrow
\tilde{E}^{\circ}\otimes\Omega^1_{(C\setminus\{t_1,\ldots,t_n\})\times T/T}$
such that
$(\tilde{E}^{\circ},\tilde{\nabla}^{\circ})|
_{(C\setminus\{t_1,\ldots,t_n\})\times\eta}$
corresponds to the representation $f(\eta)$ which is irreducible
for $\eta\in T\setminus\{p_0\}$
and
$(\tilde{E}^{\circ},\tilde{\nabla}^{\circ})|
_{(C\setminus\{t_1,\ldots,t_n\})\times p_0}$
corresponds to the representation $\rho$.
By the relative version of  [\cite{Deligne}, II, Proposition 5.4],
we can construct, after replacing $T$ by an analytic open neighborhood of $p_0$,
a holomorphic vector bundle
$\tilde{E}$ on $C\times T$ and a relative connection
$\tilde{\nabla}:\tilde{E}\rightarrow
\tilde{E}\otimes\Omega^1_C(t_1+\cdots+t_n)$
such that
$(\tilde{E},\tilde{\nabla})|_{(C\setminus\{t_1,\ldots,t_n\})\times T}
\cong (\tilde{E}^{\circ},\tilde{\nabla}^{\circ})$
and that
all the eigenvalues of $\res_{t_i}(\tilde{\nabla}\otimes k(t))$
lie in $\{z\in\C|-\epsilon< \mathrm{Re}(z)<1-\epsilon\}$
for any $t\in T$,
where $\epsilon>0$ is a sufficiently small positive real number
which depends on $\bmu$.
We can take a flat family of parabolic structures
$\{\tilde{l}^{(i)}_j\}$ on $\tilde{E}$
which is compatible with $\tilde{\nabla}$
and $(\tilde{E},\tilde{\nabla},\{\tilde{l}^{(i)}_j\})$
becomes a flat family of parabolic connections over $T$.

We denote by $(An)$ the category of analytic spaces.
Let ${\mathcal M}_C(\bt,r,d)^{an}:(An) \rightarrow (Sets)$
be the functor defined by
\[
 {\mathcal M}_C(\bt,r,d)^{an}(S):=\left\{
 (E,\nabla,\{l^{(i)}_j\})\left|
 \begin{array}{l}
  \text{$E$ is a holomorphic vector bundle of rank $r$ on $C\times S$,} \\ 
  \text{$\nabla:E\rightarrow E\otimes\Omega^1_{C\times S/S}((t_1+t_2+\cdots+t_n)_S)$} \\
  \text{is a relative connection and} \\ 
  \text{$E|_{t_i\times S}=l^{(i)}_0\supset l^{(i)}_1\supset\cdots \supset l^{(i)}_{r-1}\supset l^{(i)}_r=0$} \\
  \text{is a filtration by subbundles such that} \\
  \text{$\res_{t_i\times S}(\nabla)(l^{(i)}_j)\subset l^{(i)}_j$ for any $i,j$} \\
  \text{and for any $s\in S$, $\deg(E|_{C\times\{s\}})=d$} \\
  \text{and $\dim(l^{(i)}_j/l^{(i)}_{j+1})|_s=1$ for any $i,j$} 
 \end{array}
 \right\}\right/\sim
\]
for an analytic space $S$, where $(E,\nabla,\{l^{(i)}_j\})\sim(E',\nabla',\{{l'}^{(i)}_j\})$
if there are line bundle ${\mathcal L}$ on $S$ and an isomorphism
$g:E\stackrel{\sim}\longrightarrow E'\otimes{\mathcal L}$
such that $g|_{t_i}(l^{(i)}_j)={l'}^{(i)}_j$ for any $i,j$ and the diagram
\[
\begin{CD}
 E @>\nabla>> E\otimes\Omega^1_{C\times S/S}((t_1+\cdots+t_n)_S) \\
 g @VV\cong V @Vg\otimes\mathrm{id}V\cong V \\
 E'\otimes{\mathcal L} @>\nabla'\otimes{\mathcal L}>>
 E'\otimes{\mathcal L}\otimes\Omega^1_{C\times S/S}((t_1+\cdots+t_n)_S)
\end{CD}
\]
is commutative.
We can also define a functor ${\mathcal M}_C(\bt,\blambda)^{an}:(An)\rightarrow(Sets)$
similarly.

There are canonical morphisms
$\pi_{\Lambda^{(n)}_r(d)}:{\mathcal M}_C(\bt,r,d)^{an}\rightarrow\Lambda^{(n)}_r(d)$ and
$\pi_{\Lambda^{(n)}_r(d')}:{\mathcal M}_C(\bt,r,d')^{an}\rightarrow\Lambda^{(n)}_r(d')$.
Taking a sufficiently small open neighborhood $\Delta$ of $\blambda$ in $\Lambda^{(n)}_r(d)$
and an open neighborhood $\Delta'$ of $\bmu$ in $\Lambda^{(n)}_r(d')$,
$\sigma$ extends canonically to an isomorphism
\[
 \tilde{\sigma}:\pi_{\Lambda^{(n)}_r(d)}^{-1}(\Delta)\stackrel{\sim}\longrightarrow
 \pi_{\Lambda^{(n)}_r(d')}^{-1}(\Delta')
\]
which makes the diagram
\[
\begin{CD}
 \pi_{\Lambda^{(n)}_r(d)}^{-1}(\Delta) @>>> \Delta \\
 @V\tilde{\sigma}V\cong V @V\cong VV \\
 \pi_{\Lambda^{(n)}_r(d')}^{-1}(\Delta') @>>> \Delta'
\end{CD}
\]
commute.
Replacing $T$ by an open neighborhood of $p_0$,
we have
$(\tilde{E},\tilde{\nabla},\{\tilde{l}^{(i)}_j\})\in
\pi_{\Lambda^{(n)}_r(d')}^{-1}(\Delta')(T)$.
Then we have
$\tilde{\sigma}^{-1}(\tilde{E},\tilde{\nabla},\{\tilde{l}^{(i)}_j\})\in
\pi_{\Lambda^{(n)}_r(d)}^{-1}(\Delta)(T)\subset
{\mathcal M}_C(\bt,r,d)^{an}(T)$.
For a sufficiently large integer $n_0$,
$H^1(\tilde{\sigma}^{-1}(\tilde{E})_t(n_0))=0$ and
$H^0(\tilde{\sigma}^{-1}(\tilde{E})_t(n_0))\otimes{\mathcal O}_C(-n_0)
\longrightarrow\tilde{\sigma}^{-1}(\tilde{E})_t$
is surjective for any $t\in T$.
Let ${\mathcal M}'_C(\bt,r,d)$ be the subfunctor of
${\mathcal M}_C(\bt,r,d)$ defined by
\[
 {\mathcal M}'_C(\bt,r,d)(S)=\left\{
 (E,\nabla,\{l^{(i)}\})\in{\mathcal M}_C(\bt,r,d)(S)
 \left|
 \begin{array}{l}
  \text{for any $t\in S$, $H^1(E_t(n_0))=0$ and} \\
  \text{$H^0(E_t(n_0))\otimes{\mathcal O}_C(-n_0)\rightarrow E_t$
  is surjective}
 \end{array}
 \right.\right\}.
\]
Then ${\mathcal M}'_C(\bt,r,d)\hookrightarrow{\mathcal M}_C(\bt,r,d)$
is representable by an open immersion and we have
$\tilde{\sigma}^{-1}(\tilde{E},\tilde{\nabla},\{\tilde{l}^{(i)}_j\})
\in{\mathcal M}'_C(\bt,r,d)^{an}(T)$.
By the same argument as the proof of Theorem \ref{moduli-exists},
there is an algebraic scheme $P$ on which $PGL_N$ acts and
there is a flat family $(F,\tilde{\nabla}_F,\{l^{(i)}_{F,j}\})$
of parabolic connections on $C\times P$ over $P$ which induces
$PGL_N$-equivariant morphism
$p:P\rightarrow{\mathcal M}'_C(\bt,r,d)$
such that the induced morphism
$(h_P/h_{PGL_N})^{Zar}\rightarrow {\mathcal M}'_C(\bt,r,d)^{Zar}$
is isomorphic,
where $(h_P/h_{PGL_N})^{Zar}$ (resp.\ ${\mathcal M}'_C(\bt,r,d)^{Zar}$)
is the sheafification of $h_P/h_{PGL_N}$ (resp. ${\mathcal M}'_C(\bt,r,d)$)
in Zariski topology.
We can see that the canonical morphism
$(h_{P^{an}}/h_{PGL_N^{an}})^{\sim} \longrightarrow
({\mathcal M}'_C(\bt,r,d)^{an})^{\sim}$
is also an isomorphism, where
$(h_{P^{an}}/h_{PGL_N^{an}})^{\sim}$ (resp.\ $({\mathcal M}'_C(\bt,r,d)^{an})^{\sim}$)
 is the analytic sheafification of
$h_{P^{an}}/h_{PGL_N^{an}}$ (resp.\ ${\mathcal M}'_C(\bt,r,d)^{an}$).
Note that
$\tilde{\sigma}^{-1}(\tilde{E},\tilde{\nabla},\{\tilde{l}^{(i)}_j\})
\in{\mathcal M}'_C(\bt,r,d)^{an}(T)
\subset({\mathcal M}'_C(\bt,r,d)^{an})^{\sim}(T)$.
Since the morphism $h_P\rightarrow ({\mathcal M}'_C(\bt,r,d)^{an})^{\sim}$
is surjective as sheaves in analytic topology,
there is a morphism $g:T\rightarrow P$ such that
$g^*(F,\tilde{\nabla}_F,\{l^{(i)}_{F,j}\})^{an}\cong
\tilde{\sigma}^{-1}(\tilde{E},\tilde{\nabla},\{\tilde{l}^{(i)}_j\})$
after replacing $T$ by an open neighborhood of $p_0$.
There is a ring homomorphism
$g^*:{\mathcal O}^{an}_{P,g(p_0)}\rightarrow{\mathcal O}_{T,p_0}^{an}\cong\mathbf{C}\{z\}$,
where $\mathbf{C}\{z\}$ is the convergent power series ring.
Then we have
\[
 (F,\tilde{\nabla}_F,\{l^{(i)}_{F,j}\})\otimes_{{\mathcal O}_{P,g(p_0)}}{\mathcal O}^{an}_{T,p_0}
 \in{\mathcal M}'_C(\bt,r,d)({\mathcal O}^{an}_{T,p_0})=
 {\mathcal M}'_C(\bt,r,d)(\mathbf{C}\{z\}).
 \]
By Lemma \ref{langton}, we can obtain a flat family
$(\tilde{E'},\tilde{\nabla'},\{\tilde{l'}^{(i)}_j\})$
of $\balpha$-semistable parabolic connections over
$\mathbf{C}\{z\}$ such that
$(\tilde{E'},\tilde{\nabla'},\{\tilde{l'}^{(i)}_j\})\otimes_{\mathbf{C}\{z\}}\mathbf{C}\{z\}[1/z]\cong
(F,\tilde{\nabla}_F,\{l^{(i)}_{F,j}\})\otimes_{{\mathcal O}_{P,g(p_0)}}\mathbf{C}\{z\}[1/z]$.
Since we take $\balpha$ so generic that
$\balpha$-semistable $\Leftrightarrow$ $\balpha$-stable,
we have
$(\tilde{E'},\tilde{\nabla'},\{\tilde{l'}^{(i)}_j\})\in M^{\balpha}_C(\bt,r,d)(\mathbf{C}\{z\})$.
Let $h:\Spec\mathbf{C}\{z\}\rightarrow M^{\balpha}_C(\bt,r,d)$ be the induced morphism
and
$h^*:{\mathcal O}_{M^{\balpha}_C(\bt,r,d),h(q_0)}\rightarrow\mathbf{C}\{z\}$
be the corresponding ring homomorphism,
where $q_0$ is the closed point of $\Spec\mathbf{C}\{z\}$.
Then, replacing $T$ by an open neighborhood of $p_0$, we have a morphism
$h^{an}:T\rightarrow M^{\balpha}_C(\bt,r,d)$ such that
$(h^{an})^*=h^*:{\mathcal O}_{M^{\balpha}_C(\bt,r,d),h(q_0)}
\rightarrow{\mathcal O}^{an}_{T,p_0}=\mathbf{C}\{z\}$.
By construction, we have $h^{an}(p_0)=h(q_0)$ and
the diagram
\[
 \begin{CD}
  T\setminus\{p_0\} @>f>>
  \Hom(\pi_1(C\setminus\{t_1,\ldots,t_n\},*),GL_r(\C))\\
  @Vh^{an}VV @VVV \\
   M^{\balpha}(\bt,r,d) 
  @>\RH>> RP_r(C,\bt)
 \end{CD}
\]
is commutative after replacing $T$ by an open neighborhood of $p_0$.
Since $RP_r(C,\bt)$ is separated over $\mathbf{C}$, we have
$\RH(h^{an}(p_0))=[\rho]$.
Hence $\RH_{(x,\blambda)}$ is surjective.
\end{proof}

\begin{Remark}\rm
It was not obvious in
Proposition \ref{irred-generic}, \ref{irred-g=1} and
\ref{irred-g=0} that the moduli space
$M_C^{\balpha}(\bt,\blambda,L)$ is non-empty.
However, we can see now that it is in fact non-empty
in the following way:
We first take a simple
parabolic bundle $(E,\{l^{(i)}_j\})$
of rank $r$ on $C$ such that $\det E\cong L$.
For the existence of such a parabolic bundle, we may assume,
by taking certain elementary transform, that $r$ and $\deg L$ are coprime.
Then the existence for $g\geq 1$ follows because the moduli space
of stable bundles of rank $r$ with the fixed determinant $L$
on a curve of genus $g\geq 1$ is non-empty. 
For the case of $C=\mathbf{P}^1$, we may assume that
$E=\cO_{\BP^1}^{\oplus r}$ and
we can take a parabolic structure $\{l^{(i)}_j\}$ on $E$
such that $\End(E,\{l^{(i)}_j\})\cong\C$
because $n\geq 3$.
We want to construct a connection
$\nabla:E\rightarrow E\otimes\Omega^1_C(t_1+\cdots+t_n)$
such that $(E,\nabla,\{l^{(i)}_j\})$ becomes a
$(\bt,\blambda)$-parabolic connection
with the determinant $(L,\nabla_L)$.
The obstruction class for the existence of such a connection
lies in $H^1(\cF_0^1)$, where
\[
 \cF_0^1=\left\{ s\in {\mathcal End}(E)\otimes\Omega^1_C(t_1+\cdots+t_n)
 \left| \text{$\Tr(s)=0$ and
 $\res_{t_i}(a)(l^{(i)}_j)\subset l^{(i)}_{j+1}$ for any $i,j$}
 \right\}\right..
\]
If we put
\[
 \cF_0^0:=\left\{ s\in{\mathcal End}(E) \left| \text{$\Tr(s)=0$ and
 $s|_{t_i}(l^{(i)}_j)\subset l^{(i)}_j$ for any $i,j$} \right\}\right.,
\]
then we have an isomorphism
$H^1(\cF_0^1)\cong H^0(\cF_0^0)^{\vee}$.
Since $\End(E,\{l^{(i)}_j\})\cong\C$,
we have $H^0(\cF_0^0)=0$ and so
$H^1(\cF_0^1)=0$.
Thus there is a connection
$\nabla:E\rightarrow E\otimes\Omega^1_C(t_1+\cdots+t_n)$
such that $(E,\nabla,\{l^{(i)}_j\})$ becomes a
$(\bt,\blambda)$-parabolic connection.
If this connection is irreducible, then
$(E,\nabla,\{l^{(i)}_j\})$ becomes $\balpha$-stable
and the moduli space $M_C^{\balpha}(\bt,\blambda,L)$ becomes
non-empty.
If $(E,\nabla,\{l^{(i)}_j\})$ is reducible,
then we take the representation $\rho$
of $\pi_1(C\setminus\{t_1,\ldots,t_n\},*)$
which corresponds to $\ker(\nabla^{an}|_{C\setminus\{t_1,\ldots,t_n\}})$.
Then $[\rho]$ determines a point of $RP_r(C,\bt)_{\ba}$,
where $rh(\blambda)=\ba$.
Since $\RH_{(x,\blambda)}$ is surjective, we can take an $\balpha$-stable
$(\bt,\blambda)$-parabolic connection $(E',\nabla',\{{l'}^{(i)}_j\})$
with the determinant $(L,\nabla_L)$ such that
$\RH_{(x,\blambda)}(E',\nabla',\{{l'}^{(i)}_j\})=[\rho]$.
Hence we have $M_C^{\balpha}(\bt,\blambda,L)\neq\emptyset$.
\end{Remark}

\begin{Proposition}\label{compact}
 Every fiber of
 $\RH_{(x,\blambda)}:M^{\balpha}_C(\bt,\blambda)
 \rightarrow RP_r(C,\bt)_{\ba}$
 is compact.
\end{Proposition}

\begin{proof}
Take any irreducible representation
$\rho\in RP_r(C,\bt)_{\ba}^{\mathrm{irr}}$.
We want to show that the fiber $\RH_{(x,\blambda)}^{-1}(\rho)$ in
$M^{\balpha}_C(\bt,\blambda)$ is compact.
By Proposition \ref{prop:elementary-transform},
we can obtain an isomorphism
\[
 \sigma:\cM_C(\bt,\blambda)\stackrel{\sim}\longrightarrow
 \cM_C(\bt,\bmu),
\]
of functors, where
$0\leq\mathrm{Re}(\mu^{(i)}_j)<1$ for any $i,j$.
Then $\sigma$ induces an isomorphism
\[
 \sigma:M_C^{\mathrm{irr}}(\bt,\blambda)\stackrel{\sim}\longrightarrow
 M_C^{\mathrm{irr}}(\bt,\bmu)
\]
between the moduli spaces of irreducible parabolic connections.
Note that
$M_C^{\mathrm{irr}}(\bt,\blambda)\subset M^{\balpha}_C(\bt,\blambda)$.
So $\sigma$ also induces an isomorphism
\[
 \RH_{(x,\blambda)}^{-1}(\rho)
 \stackrel{\sim}\longrightarrow
 \RH_{(x,\bmu)}^{-1}(\rho).
\]
By [\cite{Deligne}, II, Proposition 5.4], there is a unique pair
$(E,\nabla_E)$ of a bundle $E$ on $C$ and a logarithmic connection
$\nabla_E:E\rightarrow E\otimes\Omega^1_C(t_1+\cdots+t_n)$ such that
the local system $\ker(\nabla^{an}_E)|_{C\setminus\{t_1,\ldots,t_n\}}$
corresponds to the representation $\rho$
and all the eigenvalues of $\res_{t_i}(\nabla_E)$ lies in
$\{z\in\C|0\leq\mathrm{Re}(z)<1\}$.
Then the fiber
$\RH_{(x,\bmu)}^{-1}(\rho)$
is just the moduli space of parabolic structures $\{l^{(i)}_j\}$ on $E$
which satisfy
$(\res_{t_i}(\nabla_E)-\mu^{(i)}_j)(l^{(i)}_j)\subset l^{(i)}_{j+1}$
for any $i,j$.
So $\RH_{(x,\bmu)}^{-1}(\rho)$
becomes a Zariski closed subset of a product of flag varieties.
Thus $\RH_{(x,\bmu)}^{-1}(\rho)$ is
compact.
Since the fiber $\RH_{(x,\blambda)}^{-1}(\rho)$ of $\RH_{(x,\blambda)}$
over $\rho$ in
$M^{\balpha}_C(\bt,\blambda)$ is isomorphic to
$\RH_{(x,\bmu)}^{-1}(\rho)$ via $\sigma$, it is also compact.

Now take a reducible representation $[\rho]\in RP_r(C,\bt)_{\ba}$.
We may assume
$\rho=\rho_1\oplus\cdots\oplus\rho_s$
with each $\rho_p$ irreducible.
We will show that $\RH_{(x,\blambda)}^{-1}([\rho])$
is an algebraic constructible subset of
$M^{\balpha}_C(\bt,\blambda)$.
For this, we will construct an algebraic scheme $Z$
over $\mathbf{C}$ and a morphism
$f:Z\rightarrow M^{\balpha}_C(\bt,\blambda)$
such that
$f(Z)=\RH^{-1}([\rho])$.
By Proposition \ref{prop:elementary-transform},
we can get an isomorphism
\[
 \varphi:{\mathcal M}_C(\bt,\blambda)\longrightarrow
 {\mathcal M}_C(\bt,\bmu)
\]
of functors such that
$0\leq \mathrm{Re}(\mu^{(i)}_j)<1$ for any $i,j$ and
$\RH_{(x,\blambda)}(e)
=(\RH_{(x,\bmu)}\circ\varphi)(e)$
for any 
$e\in{\mathcal M}_C(\bt,\blambda)(\mathbf{C})$.
Take any member $e\in M^{\balpha}_C(\bt,\blambda)$
such that $\RH_{(x,\blambda)}(e)=[\rho]$.
Put
$\varphi(e)=(E,\nabla,\{l^{(i)}_j\})$.
Then there is  filtration
$0=F_0\subset F_1\subset F_2\subset\cdots\subset F_s=E$
by subbundles such that 
$\nabla(F_k)\subset F_k\otimes\Omega^1_C(t_1+\cdots+t_n)$
and for the induced connections
$\nabla_k:F_k/F_{k-1}\rightarrow F_k/F_{k-1}
\otimes\Omega^1_C(t_1+\cdots+t_n)$,
there is a permutation $\sigma$ of
$1,\ldots,s$ such that
$\ker\nabla^{an}_k|_{C\setminus\{t_1,\ldots,t_n\}}$
corresponds to $\rho_{\sigma(k)}$
for any $k$.
Note that $(F_k/F_{k-1},\nabla_k)$ is
uniquely determined by $\rho_{\sigma(k)}$
because $\rho_{\sigma(k)}$ is irreducible
and all the eigenvalues of $\res_{t_i}(\nabla_k)$
lie in $\{z\in\mathbf{C}|0\leq\mathrm{Re}(z)<1\}$.
Now we will construct $Z$ and $f$.
For each $k$ with $1\leq k\leq s$,
there is a unique pair $(E_k,\nabla_k)$ of vector bundle
$E_k$ with a connection
$\nabla_k:E_k\rightarrow E_k\otimes\Omega^1_C(t_1+\cdots+t_n)$
such that 
$\ker\nabla_k^{an}|_{C\setminus\{t_1,\ldots,t_n\}}$ corresponds to
$\rho_k$ and
all the eigenvalues of $\res_{t_i}(\nabla)$ lie in
$\{z\in\mathbf{C}|0\leq\mathrm{Re}(z)<1\}$
for $i=1,\ldots,n$.
Take a permutation $\sigma$ of $(1,\ldots,s)$.
Put $Y_1:=\mathbf{V}_*(\Ext^1_C(E_{\sigma(2)},E_{\sigma(1)}))$
and let
\[
 0\longrightarrow (E_{\sigma(1)})_{Y_1}\longrightarrow
 F^{(\sigma)}_2\longrightarrow (E_{\sigma(2)})_{Y_1}\longrightarrow 0
\]
be the universal extension.
We may assume that
$\Ext^1_{C\times Y_1/Y_1}((E_{\sigma(3)})_{Y_1},F^{(\sigma)}_2)$
is a locally free sheaf on $Y_1$
by replacing $Y_1$ by a disjoint union of locally closed subschemes.
Inductively we put
$Y_k:=\mathbf{V}_*(\Ext^1_{C\times Y_{k-1}/Y_{k-1}}
((E_{\sigma(k+1)})_{Y_{k-1}},F^{(\sigma)}_k))$
and let
\[
 0\longrightarrow (F^{(\sigma)}_k)_{Y_k}\longrightarrow F^{(\sigma)}_{k+1}
 \longrightarrow (E_{\sigma(k+1)})_{Y_k}\longrightarrow 0
\]
be the universal extension.
We may assume by replacing $Y_k$ by a disjoint union of
locally closed subschemes that
$\Ext^1_{C\times Y_k/Y_k}((E_{\sigma(k+2)})_{Y_k},F^{(\sigma)}_{k+1})$
is a locally free sheaf on $Y_k$.
There is a universal extension
\[
 0\longrightarrow (F^{(\sigma)}_{r-1})_{Y_{r-1}} \longrightarrow
 F^{(\sigma)}_r\longrightarrow (E_{\sigma(r)})_{Y_{r-1}}
 \longrightarrow 0
\]
and we obtain a filtration
\[
 0=:F^{(\sigma)}_0\subset F^{(\sigma)}_1:=(E_{\sigma(1)})_{Y_{r-1}}
 \subset(F^{(\sigma)}_2)_{Y_{r-1}}
 \subset\cdots\subset(F^{(\sigma)}_{r-1})_{Y_{r-1}}
 \subset F^{(\sigma)}_r
\]
Consider the functor $\mathrm{Conn}_{C\times Y_{r-1}/Y_{r-1}}$
from the category of algebraic schemes over $\mathbf{C}$
to the category of sets defined by
\[
 \mathrm{Conn}_{C\times Y_{r-1}/Y_{r-1}}(S):=\left\{
 \nabla:F^{(\sigma)}_r\rightarrow
 F^{(\sigma)}_r\otimes\Omega^1_C(D(\bt))\left|
 \begin{array}{l}
 \text{$\nabla$ is a relative connection,} \\
 \nabla((F^{(\sigma)}_k)_S)\subset(F^{(\sigma)}_k)_S
 \otimes\Omega^1_C(D(\bt)) \\
 \text{for $1\leq k\leq r-1$ and for the connection} \\
 \text{$\nabla^{(\sigma)}_k$ on
 $(F^{(\sigma)}_k)_S/(F^{(\sigma)}_{k-1})_S$ induced by $\nabla$,} \\
 \text{$\nabla^{(\sigma)}_k=(\nabla_{\sigma(k)})_S$
 for $1\leq k\leq s-1$}
 \end{array}
 \right\}\right..
\]
Then we can see that $\mathrm{Conn}_{Y_{r-1}}$
can be represented by an affine scheme
$Y^{(\sigma)}$ over $Y_{r-1}$.
Let
\[
 \nabla:(F^{(\sigma)}_r)_{Y^{\sigma}}\longrightarrow
 (F^{(\sigma)}_r)_{Y^{\sigma}}
 \otimes\Omega^1_C(t_1+\cdots+t_n)
\]
be the universal connection.
The set of parabolic structures $\{l^{(i)}_j\}$
on $(F^{(\sigma)}_r)\otimes k(y)$ ($y\in Y^{(\sigma)}$)
such that
$(\res_{t_i}(\nabla\otimes k(y))-\mu^{(i)}_j)(l^{(i)}_j)
\subset l^{(i)}_{j+1}$
for any $i,j$
can by parametrized by a closed subscheme $Z^{(\sigma)}$ of a
product of flag schemes over $Y$.
Let $\{(l^{(\sigma)})^{(i)}_j\}$ be the universal parabolic
structure on $(F^{(\sigma)}_r)_{Z^{(\sigma)}}$.
Put
\[
 \tilde{Z}^{(\sigma)}:=\{z\in Z^{(\sigma)}|
 \text{$\varphi^{-1}(F^{(\sigma)}_r\otimes k(z),\nabla\otimes k(z),
 \{\tilde{l}^{(i)}_j\otimes k(s)\})$ is $\balpha$-stable}\}.
\]
Then $\tilde{Z}^{(\sigma)}$ is a Zariski open subset of
$Z^{(\sigma)}$ and induces a morphism
$f_{\sigma}:\tilde{Z}^{(\sigma)}\rightarrow M^{\balpha}_C(\bt,\blambda)$.
By construction we have
$\RH_{(x,\blambda)}^{-1}([\rho])=
\bigcup_{\sigma}f_{\sigma}(\tilde{Z}^{(\sigma)})$.
Thus we can see that $\RH_{(x,\blambda)}^{-1}([\rho])$ is
a constructible subset of $M^{\balpha}_C(\bt,\blambda)$
by Chevalley's Theorem.

So it is sufficient to show that
$\RH_{(x,\blambda)}^{-1}([\rho])$ is stable
under specialization in the compact moduli space
$\overline{M^{D(\bt),\balpha',\bbeta,\gamma}_C}
(r,d,\{d_i\}_{1\leq i\leq nr})$,
in order to prove the compactness of
$\RH_{(x,\blambda)}^{-1}([\rho])$.
Here $\overline{M^{D(\bt),\balpha',\bbeta,\gamma}_C}
(r,d,\{d_i\}_{1\leq i\leq nr})$ is the moduli scheme of
$(\balpha',\bbeta,\gamma)$-stable parabolic
$\Lambda_{D(\bt)}^1$-triples appeared in Theorem \ref{moduli-p-triple},
which contains $M^{\balpha}_C(\bt,\blambda)$ as a locally closed subscheme.
Take any scheme point $x_1\in\RH_{(x,\blambda)}^{-1}([\rho])$ and a point
$x_0$ of the closure $\overline{\{x_1\}}$ of $\{x_1\}$ in 
$\overline{M^{D(\bt),\balpha',\bbeta,\gamma}_C}
(r,d,\{d_i\}_{1\leq i\leq nr})$.
Let $K$ be the residue field $k(x_1)$ of $x_1$.
Then there exists a discrete valuation ring $R$
with quotient field $K$ which dominates $\cO_{\overline{\{x_1\}},x_0}$.
The inclusion $\cO_{\overline{\{x_1\}},x_0}\hookrightarrow R$ induces a flat family
$(\tilde{E}_1,\tilde{E}_2,\tilde{\Phi},F_*(\tilde{E}_1))$
of parabolic $\Lambda^1_{D(\bt)}$-triples
on $C\times\Spec R$ over $R$.
Here the left $\cO_{C\times\Spec R}$-homomorphism
$\tilde{\Phi}:\Lambda^1_{D(\bt)}\otimes\tilde{E}_1
\rightarrow\tilde{E}_2$
corresponds to a pair $(\tilde{\phi},\tilde{\nabla})$ of
an $\cO_{C\times\Spec R}$-homomorphism
$\tilde{\phi}:\tilde{E}_1\rightarrow\tilde{E}_2$ and a morphism
$\tilde{\nabla}:\tilde{E}_1\rightarrow
\tilde{E}_2\otimes\Omega^1_C(t_1+\cdots+t_n)$
satisfying $\tilde{\nabla}(as)=\phi(s)\otimes da+a\tilde{\nabla}(s)$
for $a\in\cO_{C\times\Spec R}$ and $s\in\tilde{E}_1$.
We also denote $(\tilde{E}_1,\tilde{E}_2,\tilde{\Phi},F_*(\tilde{E}_1))$
by $(\tilde{E}_1,\tilde{E}_2,\tilde{\phi},\tilde{\nabla},F_*(\tilde{E}_1))$.
Let $\eta$ be the generic point of $\Spec R$ and
$\xi$ be the closed point of $\Spec R$ and put
\begin{align*}
 (E_1,E_2,\phi,\nabla,F_*(E_1))&:=
 (\tilde{E}_1,\tilde{E}_2,\tilde{\phi},\tilde{\nabla},
 F_*(\tilde{E}_1))\otimes K \\
 (E'_1,E'_2,\phi',\nabla',F_*(E'_1))&:=
 (\tilde{E}_1,\tilde{E}_2,\tilde{\phi},\tilde{\nabla},
 F_*(\tilde{E}_1))\otimes k(\xi).
\end{align*}
Note that $\phi$ is isomorphic and
$\nabla_E=(\phi)^{-1}\circ\nabla$ becomes a connection on
$E:=E_1$.
There exists a filtration
\[
 (0,0)=(F^{(0)},{\nabla}^{(0)})\subset(F^{(1)},{\nabla}^{(1)})
 \subset\cdots\subset(F^{(s)},{\nabla}^{(s)})=(E,\nabla_E)\otimes_KK'
\]
such that $K'$ is a finite extension field of $K$
and that each
$(F^{(p)}/F^{(p-1)},\overline{\nabla^{(p)}})\otimes_{K'}\overline{K'}$
is irreducible, where $\overline{{\nabla}^{(p)}}$ is the connection
on $F^{(p)}/F^{(p+1)}$ induced by ${\nabla}^{(p)}$.
We may assume that $(F^{(p)}/F^{(p-1)},\overline{\nabla^{(p)}})$
corresponds to the representation $\rho_p$.
Note that there is a discrete valuation ring $R'$ with quotient field $K'$
which dominates $R$.
Let $\xi'$ be the closed point of $\Spec R'$.
The filtration $\{(F^{(p)},{\nabla}^{(p)})\}$ induces filtrations
\begin{gather*}
 0=\tilde{F}^{(0)}_1\subset\tilde{F}^{(1)}_1
 \subset\cdots\subset\tilde{F}^{(s)}_1=\tilde{E}_1\otimes_RR' \\
 0=\tilde{F}^{(0)}_2\subset\tilde{F}^{(1)}_2
 \subset\cdots\subset\tilde{F}^{(s)}_2=\tilde{E}_2\otimes_RR'
\end{gather*}
such that each $\tilde{F}^{(p)}_j/\tilde{F}^{(p-1)}_j$
is flat over $R'$ for $j=1,2$ and
$\tilde{F}^{(p)}_1\otimes K'=F^{(p)}$,
$\tilde{F}^{(p)}_2\otimes K'=(\phi\otimes\mathrm{id})(F^{(p)})$.
By construction, we have
$(\tilde{\phi}\otimes\mathrm{id})(\tilde{F}^{(p)}_1)\subset\tilde{F}^{(p)}_2$ and
$(\tilde{\nabla}\otimes\mathrm{id})(\tilde{F}^{(p)}_1)\subset
\tilde{F}^{(p)}_2\otimes\Omega^1_C(t_1+\cdots+t_n)$.
Let
\begin{gather*}
 \overline{\tilde{\phi}^{(p)}}:\tilde{F}^{(p)}_1/\tilde{F}^{(p-1)}_1
 \lra\tilde{F}^{(p)}_2/\tilde{F}^{(p-1)}_2 \\
 \overline{\tilde{\nabla}^{(p)}}:\tilde{F}^{(p)}_1/\tilde{F}^{(p-1)}_1
 \lra\tilde{F}^{(p)}_2/\tilde{F}^{(p-1)}_2\otimes\Omega^1_C(t_1+\cdots+t_n)
\end{gather*}
be the morphisms induced by $\tilde{\phi}$ and $\tilde{\nabla}$.
We can construct a parabolic structure on
$(F^{(p)}/F^{(p-1)},\overline{\nabla^{(p)}})$
and obtain a $(\bt,\blambda^{(p)})$-parabolic connections
of rank $r_p$ and degree $d_p$.
Extending this parabolic structure,
we can obtain a flat family of parabolic structures
$F_*(\tilde{F}^{(p)}_1/\tilde{F}^{(p-1)}_1)$ on
$\tilde{F}^{(p)}_1/\tilde{F}^{(p-1)}_1$.
Repeating elementary transforms of the flat family
\[
 (\tilde{F}^{(p)}_1/\tilde{F}^{(p-1)}_1,\tilde{F}^{(p)}_2/\tilde{F}^{(p-1)}_2,
 \overline{\tilde{\phi}^{(p)}},\overline{\tilde{\nabla}^{(p)}},
 F_*(\tilde{F}^{(p)}_1/\tilde{F}^{(p-1)}_1))
\]
along the special fiber $C\times\{\xi'\}$, we obtain a flat family
$(\tilde{G'}^{(p)}_1,\tilde{G'}^{(p)}_2,\tilde{\phi'}^{(p)},
\tilde{\nabla'}^{(p)},F_*(\tilde{G'}^{(p)}_1))$
of $(\balpha'_p,\bbeta_p,\gamma_p)$-semistable
$\Lambda^1_{D(\bt)}$-triples
for some weight $(\balpha'_p,\bbeta_p,\gamma_p)$.
We may assume that
$(\balpha'_p,\bbeta_p,\gamma_p)$-semistable
$\Leftrightarrow$
$(\balpha'_p,\bbeta_p,\gamma_p)$-stable.
This flat family defines a morphism
\[
 f:\Spec R' \longrightarrow
 \overline{M^{D(\bt),\balpha'_p,\bbeta_p,\gamma_p}_C}
 (r_p,d_p,\{d^{(p)}_i\}_{1\leq i\leq nr_p}).
\]
Note that the moduli space $M_C^{\balpha_p}(\bt,\blambda^{(p)})$
of $\balpha_p$-stable
$(\bt,\blambda^{(p)})$-parabolic connections
for a certain weight $\balpha_p$
is a locally closed subscheme of
$\overline{M^{D(\bt),\balpha'_p,\bbeta_p,\gamma_p}_C}
(r_p,d_p,\{d^{(p)}_i\}_{1\leq i\leq nr_p})$
and the image $f(\Spec K)$ is contained in the moduli space
$M_C^{\balpha_p}(\bt,\blambda^{(p)})$.
Since
$(\tilde{G'}^{(p)}_1,\tilde{G'}^{(p)}_2,\tilde{\phi'}^{(p)},
\tilde{\nabla'}^{(p)})\otimes K'$
is equivalent to the irreducible connection
$(F^{(p)}/F^{(p-1)},\overline{\nabla^{(p)}})$
which corresponds to a constant family of the irreducible representation
$\rho_p$, we have
$f(\Spec K')\subset\RH_{(x,\blambda^{(p)})}^{-1}(\rho_p)$,
where $\RH_{(x,\blambda^{(p)})}$ is the morphism
\[
 \RH_{(x,\blambda^{(p)})}:
 M_C^{\balpha_p}(\bt,\blambda^{(p)})\longrightarrow
 RP_{r_p}(C,\bt)_{\ba^{(p)}}
\]
determined by the Riemann-Hilbert correspondence and
$rh(\blambda^{(p)})=\ba^{(p)}$.
Since $\RH_{(x,\blambda^{(p)})}^{-1}(\rho_p)$
is compact by the previous proof
and $f(\Spec K')$ is contained in $\RH_{(x,\blambda^{(p)})}^{-1}(\rho_p)$,
the closure of $f(\Spec K')$ must be contained in
$\RH_{(x,\blambda^{(p)})}^{-1}(\rho_p)$.
In particular, $f(\xi')$ is contained in 
$\RH_{(x,\blambda^{(p)})}^{-1}(\rho_p)\subset
M_C^{\balpha_p}(\bt,\blambda^{(p)})$.
So the special fiber
$(\tilde{G'}^{(p)}_1,\tilde{G'}^{(p)}_2,\tilde{\phi'}^{(p)},
\tilde{\nabla'}^{(p)},F_*(\tilde{G'}^{(p)}_1))\otimes k(\xi')$
becomes a $(\bt,\blambda^{(p)})$-parabolic connection
which implies that $\tilde{\phi'}^{(p)}\otimes k(\xi')$ must be isomorphic.
Then $\overline{\tilde{\phi}^{(p)}}\otimes k(\xi')$ also becomes isomorphic
and $\phi'=\tilde{\phi}\otimes k(\xi)$ must be an isomorphism.
Thus $\tilde{\phi}$ is an isomorphism.
Let $\{\tilde{l}^{(i)}_j\}$ be the flat family of parabolic structures
corresponding to $F_*(\tilde{E}_1)$.
Then we have
$(\res_{t_i}(\tilde{\nabla})-\lambda^{(i)}_j\tilde{\phi}|_{t_i\otimes R})
(\tilde{l}^{(i)}_j)\subset\tilde{\phi}(\tilde{l}^{(i)}_{j+1})$.
Thus
$(E'_1,E'_2,\phi',\nabla',F_*(E'_1))\in
M^{\balpha}_C(\bt,\blambda)(R)$
which is equivalent to $x_0\in M^{\balpha}_C(\bt,\blambda)$.
Then we have $x_0\in\RH^{-1}_{(x,\blambda)}([\rho])$.
Therefore $\RH_{(x,\blambda)}^{-1}([\rho])$ is stable under specialization in
$\overline{M^{D(\bt),\balpha',\bbeta,\gamma}_C}
(r,d,\{d_i\}_{1\leq i\leq nr})$.
Hence $\RH_{(x,\blambda)}^{-1}([\rho])$ becomes a compact subset of
$M^{\balpha}_C(\bt,\blambda)$.
\end{proof}

\noindent
{\bf Proof of Theorem \ref{main-thm}.}
We will prove that
\[
 \RH_{(x,\blambda)}: M_C^{\balpha}(\bt,\blambda)\longrightarrow
 RP_r(C,\bt)_{\ba}
\]
is an analytic isomorphism for generic $\blambda$ and
gives an analytic resolution of singularities of
$RP_r(C,\bt)_{\ba}$ for special $\blambda$.
First note that $RP_r(C,\bt)_{\ba}$ is an irreducible
variety since by Proposition \ref{surjectivity}
it is the image of the irreducible
variety $M_C^{\balpha}(\bt,\blambda)$ by $\RH_{(x,\blambda)}$.
We set
\[
 RP_r(C,\bt)^{\mathrm{sing}}_{\ba}:=
 \left\{ [\rho]\in RP_r(C,\bt)_{\ba} \left|
 \begin{array}{l}
  \text{$\rho$ is reducible or} \\
  \text{$\dim\left(
  \ker\left(\rho(\gamma_i)-\exp(-2\pi\sqrt{-1}\lambda^{(i)}_j)I_r
  \right)\right)
  \geq 2$ for some $i,j$}
 \end{array}
 \right\}\right.,
\]
where $\gamma_i$ is a loop around $t_i$ for $1\leq i\leq n$.
Then we can see that
$RP_r(C,\bt)^{\mathrm{sing}}_{\ba}$ is a Zariski closed
subset of $RP_r(C,\bt)_{\ba}$.
We will see that it is a proper closed subset.
Since
$\dim M_C^{\balpha}(\bt,\blambda)=2r^2(g-1)+nr(r-1)+2$
and $\RH_{(x,\blambda)}$ is surjective, we have
$\dim RP_r(C,\bt)_{\ba}\leq 2r^2(g-1)+nr(r-1)+2$.
On the other hand, $RP_r(C,\bt)_{\ba}$ is a Zariski closed subset
of the $r^2(2g+n-2)+1$ dimensional irreducible variety $RP_r(C,\bt)$
defined by $nr-1$ equations given by
\begin{equation}\label{local-monodromy}
 \det(XI_r-\rho(\gamma_i))=
 X^r+a^{(i)}_{r-1}X^{r-1}+\cdots+a^{(i)}_1X+a^{(i)}_0
 \quad (i=1,\ldots,n).
\end{equation}
Note that the equation
$\det(\rho(\gamma_1))\cdots\det(\rho(\gamma_n))=1=
(-1)^{nr}a^{(1)}_0\cdots a^{(n)}_0$
is automatically satisfied and so
the condition (\ref{local-monodromy}) is in fact
equivalent to $nr-1$ equations.
So we have
\[
 \dim RP_r(C,\bt)_{\ba}\geq 
 r^2(2g+n-2)+1-(nr-1)=2r^2(g-1)+nr(r-1)+2
\]
and then we have $\dim RP_r(C,\bt)_{\ba}=2r^2(g-1)+nr(r-1)+2$.

First consider the locus
$RP_r(C,\bt)_{\ba}^{\mathrm{irr}}\cap RP_r(C,\bt)_{\ba}^{\mathrm{sing}}$
in $RP_r(C,\bt)^{\mathrm{sing}}_{\ba}$ corresponding to
the irreducible representations.
By Proposition \ref{prop:elementary-transform},
we can obtain an isomorphism
\[
 \sigma: \cM_C(\bt,\blambda)\stackrel{\sim}\longrightarrow
 \cM_C(\bt,\bmu),
\]
where $0\leq \mathrm{Re}(\mu^{(i)}_j)<1$
for any $i,j$.
We note that $\sigma$ also induces an isomorphism
\[
 \sigma:M^{\mathrm{irr}}_C(\bt,\blambda)\stackrel{\sim}\longrightarrow
 M^{\mathrm{irr}}_C(\bt,\bmu).
\]
Take any point
$[\rho]\in RP_r(C,\bt)_{\ba}^{\mathrm{irr}}\cap
RP_r(C,\bt)_{\ba}^{\mathrm{sing}}$.
By [\cite{Deligne}, II, Proposition 5.4], we can take
a unique pair $(E,\nabla)$ of a vector bundle $E$ on $C$
and  a logarithmic connection
$\nabla:E\rightarrow E\otimes\Omega^1_C(t_1+\cdots+t_n)$
such that $\ker\nabla^{an}|_{C\setminus\{t_1,\ldots,t_n\}}$
corresponds to the representation $\rho$
and all the eigenvalues of $\res_{t_i}(\nabla)$
lie in $\{z\in\C|0\leq \mathrm{Re}(z)<1\}$.
Then the fiber
$\RH_{(x,\bmu)}^{-1}([\rho])$
is isomorphic to the moduli space of parabolic structures
$\{l^{(i)}_j\}$ on $E$ whose dimension should be positive
because $\mu^{(i)}_j=\mu^{(i)}_{j'}$
for some $i$ and $j\neq j'$ and
the automorphisms of $(E,\nabla)$ are only scalar multiplications.
So the fiber
$\RH_{(x,\blambda)}^{-1}([\rho])$
also has positive dimension.
For a general point $\xi$ of
$RP_r(C,\bt)_{\ba}^{\mathrm{irr}}\cap
RP_r(C,\bt)_{\ba}^{\mathrm{sing}}$,
we have 
\begin{align*}
 \dim \RH_{(x,\blambda)}^{-1}(\xi)
 + \dim_{\xi}\left(RP_r(C,\bt)_{\ba}^{\mathrm{irr}}\cap
 RP_r(C,\bt)_{\ba}^{\mathrm{sing}}\right)
 &\leq \dim M^{\mathrm{irr}}_C(\bt,\blambda) \\
 &=2r^2(g-1)+nr(r-1)+2.
\end{align*}
Thus we have
$\dim \left(RP_r(C,\bt)_{\ba}^{\mathrm{irr}}\cap
 RP_r(C,\bt)_{\ba}^{\mathrm{sing}}\right)
<2r^2(g-1)+nr(r-1)+2=\dim RP_r(C,\bt)_{\ba}$.

Next we consider the locus
$RP_r(C,\bt)_{\ba}^{\mathrm{red}}\subset RP_r(C,\bt)_{\ba}^{\mathrm{sing}}$
corresponding to the reducible representations.
For a point $[\rho]$ of $RP_r(C,\bt)_{\ba}^{\mathrm{red}}$
we may write $\rho=\rho_1\oplus\rho_2$.
If $\rho_i$ is a representation of dimension $r_i$
and its local monodromy data is $\ba_i$ for $i=1,2$,
then the set of such $[\rho]$ can be parameterized by
$RP_{r_1}(C,\bt)_{\ba_1}\times RP_{r_2}(C,\bt)_{\ba_2}$.
Note that the assumption $rn-2(r+1)>0$ for $g=0$,
$n\geq 2$ for $g=1$ and $n\geq 1$ for $g\geq 2$ implies
\begin{align*}
 \dim\left(RP_{r_1}(C,\bt)_{\ba_1}\times RP_{r_2}(C,\bt)_{\ba_2}
 \right)
 &=2r_1^2(g-1)+nr_1(r_1-1)+2 +2r_2^2(g-1)+nr_2(r_2-1)+2 \\
 &\leq  2(r_1+r_2)^2(g-1)+n(r_1+r_2)(r_1+r_2-1) \\
 &=\dim RP_r(C,\bt)_{\ba}-2.
\end{align*}
Since $RP_r(C,\bt)_{\ba}^{\mathrm{red}}$ can be covered by
a finite numbers of the images of such spaces,
we have
\[
 \dim RP_r(C,\bt)_{\ba}^{\mathrm{red}}<2r^2(g-1)+nr(r-1)+2
 =\dim RP_r(C,\bt)_{\ba}.
\]
Thus $RP_r(C,\bt)_{\ba}^{\mathrm{sing}}$ becomes a proper 
closed subset of $RP_r(C,\bt)_{\ba}$.
If we set
$RP_r(C,\bt)_{\ba}^{\sharp}:=
RP_r(C,\bt)_{\ba}\setminus RP_r(C,\bt)_{\ba}^{\mathrm{sing}}$,
then it becomes a non-empty Zariski open subset of
$RP_r(C,\bt)_{\ba}$.
We put
$M_C^{\balpha}(\bt,\blambda)^{\sharp}:=
\RH_{(x,\blambda)}^{-1}(RP_r(C,\bt)_{\ba}^{\sharp})$
and consider the restriction
\[
 \RH_{(x,\blambda)}|_{M_C^{\balpha}(\bt,\blambda)^{\sharp}}:
 M_C^{\balpha}(\bt,\blambda)^{\sharp}\longrightarrow
 RP_r(C,\bt)_{\ba}^{\sharp}.
\]
Recall that $\sigma$ induces an isomorphism
$\sigma:M_C^{\balpha}(\bt,\blambda)^{\sharp}\stackrel{\sim}\longrightarrow
M_C^{\balpha}(\bt,\bmu)^{\sharp}$
which is compatible with $\RH$.
For any point $[\rho]$ of $RP_r(C,\bt)_{\ba}^{\sharp}$,
we can see by [\cite{Deligne}, II, Proposition 5.4] that
there is a unique $(\bt,\bmu)$-parabolic connection
$(E,\nabla,\{l^{(i)}_j\})$ such that
$\RH_{(x,\bmu)}((E,\nabla,\{l^{(i)}_j\})=[\rho]$
because the parabolic structure $\{l^{(i)}_j\}$ is uniquely determined
by $(E,\nabla)$.
So $\RH_{(x,\bmu)}$ gives a one to one correspondence between
the points of
$M_C^{\balpha}(\bt,\bmu)^{\sharp}$ and the points of
$RP_r(C,\bt)_{\ba}^{\sharp}$.
We can extend this correspondence to a correspondence between 
flat families.
Thus the morphism
$M_C^{\balpha}(\bt,\bmu)^{\sharp}
\xrightarrow{\RH_{(x,\bmu)}} RP_r(C,\bt)_{\ba}^{\sharp}$
becomes an isomorphism and so
$M_C^{\balpha}(\bt,\blambda)^{\sharp}\xrightarrow{\RH_{(x,\blambda)}}
RP_r(C,\bt)_{\ba}^{\sharp}$
is also an isomorphism.
Hence the morphism
$\RH_{(x,\blambda)}:
M_C^{\balpha}(\bt,\blambda)\longrightarrow RP_r(C,\bt)_{\ba}$
becomes a bimeromorphic morphism.
If $\blambda$ is generic, we have
$M^{\balpha}_C(\bt,\blambda)^{\sharp}=M^{\balpha}_C(\bt,\blambda)$,
$RP_r(C,\bt)_{\ba}^{\sharp}=RP_r(C,\bt)_{\ba}$, and so
\[
 \RH_{(x,\blambda)}:M^{\balpha}_C(\bt,\blambda)\longrightarrow
 RP_r(C,\bt)_{\ba}
\]
is an analytic isomorphism.

In order to prove the properness of $\RH_{(x,\blambda)}$,
we will use the following lemma due to A. Fujiki:

\begin{Lemma}\label{lem:fujiki}{\bf (\cite{IIS-1}, Lemma 10.3)}
 Let $f:X\ra Y$ be a surjective holomorphic mapping of
 irreducible analytic varieties.
 Assume that an analytic closed subset $S$ of $Y$ exists such that
 $\codim_Y S\geq 2$, 
 $X^{\sharp}:=f^{-1}(Y^{\sharp})$ is dense in $X$,
 where $Y^{\sharp}=Y\setminus S$ and that the restriction
 $f|_{X^{\sharp}}:X^{\sharp}\ra Y^{\sharp}$ is an analytic isomorphism.
 Moreover assume that the fibers $f^{-1}(y)$ are compact for all $y\in Y$.
 Then $f$ is a proper mapping.
\end{Lemma}

Applying the above lemma and using Proposition \ref{surjectivity}
and Proposition \ref{compact}, it suffices to prove
that $\codim_{RP_r(C,\bt)_{\ba}}(RP_r(C,\bt)_{\ba}^{\mathrm{sing}})\geq 2$
in order to obtain the properness of $\RH_{(x,\blambda)}$.
Recall that we have
$\dim RP_r(C,\bt)_{\ba}^{\mathrm{red}}\leq \dim RP_r(C,\bt)_{\ba} -2$.
On the other hand, $M_C^{\balpha}(\bt,\blambda)^{\sharp}$ is a non-empty
Zariski open subset of $M_C^{\balpha}(\bt,\blambda)$ and so we have
$\dim\left(M_C^{\balpha}(\bt,\blambda)\setminus
M_C^{\balpha}(\bt,\blambda)^{\sharp}\right)
\leq \dim M_C^{\balpha}(\bt,\blambda)-1$.
Since the dimension of every fiber of the surjective morphism
\[
 M_C^{\mathrm{irr}}(\bt,\blambda)\setminus
 M_C^{\balpha}(\bt,\blambda)^{\sharp}
 \xrightarrow{\RH_{(x,\blambda)}}
 RP_r(C,\bt)_{\ba}^{\mathrm{sing}}\cap RP_r(C,\bt)_{\ba}^{\mathrm{irr}}
\]
has positive dimension, we have
\begin{align*}
 \dim \left(RP_r(C,\bt)_{\ba}^{\mathrm{sing}}\cap
 RP_r(C,\bt)_{\ba}^{\mathrm{irr}} \right) & \leq
 \dim\left(M_C^{\mathrm{irr}}(\bt,\blambda)\setminus
 M_C^{\balpha}(\bt,\blambda)^{\sharp}\right)-1 \\
 &\leq \dim M_C^{\balpha}(\bt,\blambda)-1-1
 =\dim RP_r(C,\bt)_{\ba}-2.
\end{align*}
Thus we have
$\codim_{RP_r(C,\bt)_{\ba}}(RP_r(C,\bt)_{\ba}^{\mathrm{sing}})\geq 2$
and obtain the properness of $\RH_{(x,\blambda)}$ by Lemma \ref{lem:fujiki}.
In particular,
$\RH_{(x,\blambda)}:M_C^{\balpha}(\bt,\blambda)\rightarrow
RP_r(C,\bt)_{\ba}$
gives an analytic resolution of singularities of
$RP_r(C,\bt)_{\ba}$ for special $\blambda$.
By the same argument, we can see that the morphism
\[
 \RH: M^{\balpha}_{\cC/T}(\tilde{\bt},r,d)\times_T\tilde{T} \lra
 (RP_r(\cC,\tilde{\bt})\times_T\tilde{T})\times_{\cA^{(n)}_r}\Lambda^{(n)}_r(d)
\]
is also a proper mapping.
The existence of a symplectic structure is given in
Proposition \ref{symplectic-form}.
\hfill $\square$

\begin{Remark}\label{singular-locus}\rm
Under the assumption of Theorem \ref{main-thm},
$RP_r(C,\bt)_{\ba}$ is a normal variety and the locus
\[
 RP_r(C,\bt)_{\ba}^{\mathrm{sing}}=
 \left\{ [\rho]\in RP_r(C,\bt)_{\ba} \left|
 \begin{array}{l}
  \text{$\rho$ is reducible or} \\
  \text{$\dim\left(
  \ker\left(\rho(\gamma_i)-\exp(-2\pi\sqrt{-1}\lambda^{(i)}_j)I_r
  \right)\right)
  \geq 2$ for some $i,j$}
 \end{array}
 \right\}\right.
\]
is just the singular locus of $RP_r(C,\bt)_{\ba}$.
\end{Remark}

\begin{proof}
Note that $RP_r(C,\bt)$ is a Cohen-Macaulay
irreducible variety by [\cite{For}, Theorem 4 (2)].
In the proof of Theorem \ref{main-thm},
we showed that
$RP_r(C,\bt)_{\ba}$ is a Zariski closed subset of $RP_r(C,\bt)$
defined by $nr-1$ equations and it is of codimension
$nr-1$ in $RP_r(C,\bt)$.
Thus $RP_r(C,\bt)_{\ba}$ is Cohen-Macaulay.
Since $RP_r(C,\bt)_{\ba}^{\sharp}$ is isomorphic to
$M_C^{\balpha}(\bt,\blambda)^{\sharp}$, $RP_r(C,\bt)_{\ba}^{\sharp}$
is smooth.
So $RP_r(C,\bt)_{\ba}$ is regular in codimension one because
$\codim_{RP_r(C,\bt)_{\ba}}(RP_r(C,\bt)_{\ba}^{\mathrm{sing}})\geq 2$.
Therefore $RP_r(C,\bt)_{\ba}$ is a normal variety by Serre's criterion
for normality.

We will prove that every fiber of $\RH_{(x,\blambda)}$ over a point of
$RP_r(C,\bt)_{\ba}^{\mathrm{sing}}$ has positive dimension.
For a point
$[\rho]\in RP_r(C,\bt)_{\ba}^{\mathrm{irr}}
\cap RP_r(C,\bt)_{\ba}^{\mathrm{sing}}$,
we have already showed in the proof of Theorem \ref{main-thm} that
$\dim \RH_{(x,\blambda)}^{-1}([\rho])\geq 1$.
So take a point $[\rho]\in RP_r(C,\bt)_{\ba}^{\mathrm{red}}$.
By Proposition \ref{prop:elementary-transform},
we can take an isomorphism
\[
 \sigma:\cM_C(\bt,\blambda)\stackrel{\sim}\longrightarrow
 \cM_C(\bt,\bmu),
\]
where $0\leq \mathrm{Re}(\mu^{(i)}_j)<1$
for any $i,j$.

First we assume that $\rho$ is an extension of two irreducible
representations $\rho_1,\rho_2$ and that
$\mu^{(i)}_j\neq\mu^{(i)}_{j'}$
for any $j\neq j'$ and any $i$.
We denote the dimension of $\rho_i$ by $r_i$.
Since $\RH_{(x,\blambda)}$ is surjective by Proposition \ref{surjectivity},
there is an $\balpha$-stable $(\bt,\blambda)$-parabolic connection
$(E,\nabla_E,\{l^{(i)}_j\})$ such that
$\RH_{(x,\blambda)}(E,\nabla_E,\{l^{(i)}_j\})=[\rho]$.
We put
$\sigma(E,\nabla_E,\{l^{(i)}_j\})=(E',\nabla_{E'},\{{l'}^{(i)}_j\})$.
Since $\mu^{(i)}_j\neq\mu^{(i)}_{j'}$
for any $j\neq j'$ and any $i$, the parabolic structure
$\{{l'}^{(i)}_j\}$ is uniquely determined by $(E',\nabla_{E'})$.
So we have
\[
 \End(E',\nabla_{E'})=\End(E',\nabla_{E'},\{{l'}^{(i)}_j\})
 \cong\End(E,\nabla_E,\{l^{(i)}_j\})
\]
which is isomorphic to $\C$ because
$(E,\nabla_E,\{l^{(i)}_j\})$ is $\balpha$-stable.
Then $(E',\nabla_{E'})$ can not split.
Thus the representation $\rho_{E'}$ corresponding to
$\ker\nabla_{E'}^{an}|_{C\setminus\{t_1,\ldots,t_n\}}$ can not split
and so we may assume that $\rho_{E'}$ is given by matrices
\begin{gather*}
 \rho_{E'}(\alpha_k)=
 \begin{pmatrix}
  \rho_1(\alpha_k) & A^{(0)}_k \\
  0 & \rho_2(\alpha_k)
 \end{pmatrix}
 \quad
 \rho_{E'}(\beta_k)=
 \begin{pmatrix}
  \rho_1(\beta_k) & B^{(0)}_k \\
  0 & \rho_2(\beta_k)
 \end{pmatrix} \\
 \rho_{E'}(\gamma_i)=
 \begin{pmatrix}
  \rho_1(\gamma_i) & C^{(0)}_i \\
  0 & \rho_2(\gamma_i)
 \end{pmatrix}
 \quad (1\leq k\leq g,\; 1\leq i\leq n-1).
\end{gather*}
Note that $\rho_{E'}(\gamma_n)$ is uniquely determined by
the above data.
Consider the family of representations over $M(r_1,r_2,\C)^{2g+n-1}$
given by matrices
\begin{gather*}
 \tilde{\rho}(\alpha_k)=
 \begin{pmatrix}
  \rho_1(\alpha_k) & A^{(0)}_k+A_k \\
  0 & \rho_2(\alpha_k)
 \end{pmatrix}
 \quad
 \tilde{\rho}(\beta_k)=
 \begin{pmatrix}
  \rho_1(\beta_k) & B^{(0)}_k+B_k \\
  0 & \rho_2(\beta_k)
 \end{pmatrix} \\
 \tilde{\rho}(\gamma_i)=
 \begin{pmatrix}
  \rho_1(\gamma_i) & C^{(0)}_i+C_i \\
  0 & \rho_2(\gamma_i)
 \end{pmatrix}
 \quad (1\leq k\leq g,\; 1\leq i\leq n-1),
\end{gather*}
where matrices $A_k,B_k (1\leq k\leq g), C_i (1\leq i\leq n-1)$
move around in $M(r_1,r_2,\C)^{2g+n-1}$.
Applying the relative version of [\cite{Deligne}, II, Proposition 5.4]
and the transform $\sigma^{-1}$,
we can obtain an analytic flat family
$(\tilde{E},\tilde{\nabla},\{\tilde{l}^{(i)}_j\})$
of $(\bt,\blambda)$-parabolic connections on
$C\times M(r_1,r_2,\C)^{2g+n-1}$ over $M(r_1,r_2,\C)^{2g+n-1}$
which corresponds to $\tilde{\rho}$.
We denote the point of $M(r_1,r_2,\C)^{2g+n-1}$ corresponding to
$A_k=B_k=C_i=0$ ($1\leq k\leq g, 1\leq i\leq n-1$) by $p_0$.
Then the fiber of $(\tilde{E},\tilde{\nabla},\{\tilde{l}^{(i)}_j\})$
over $p_0$ is just $(E,\nabla_E,\{l^{(i)}_j\})$ which is $\balpha$-stable.
So there is an analytic open neighborhood $U$ of $p_0$ such that
every fiber of $(\tilde{E},\tilde{\nabla},\{\tilde{l}^{(i)}_j\})$
over a point of $U$ is $\balpha$-stable.
Then the family $(\tilde{E},\tilde{\nabla},\{\tilde{l}^{(i)}_j\})_U$
induces a morphism
$f:U\rightarrow M_C^{\balpha}(\bt,\blambda)$
such that $f(U)\subset \RH_{(x,\blambda)}^{-1}([\rho])$.
The group
\[
 G=\left.\left\{
 \begin{pmatrix}
  uI_{r_1} & V \\
  0 & I_{r_2}
 \end{pmatrix}
 \right|
 \begin{array}{l}
  u\in\C^{\times} \\
  V\in M(r_1,r_2,\C)
 \end{array}
 \right\} 
\]
acts on $M(r_1,r_2,\C)^{2g+n-1}$ by the adjoint action.
If $f(x)=f(y)$, then there is an element
$g\in G$ such that $x=gy$.
So we have
$\dim f(U) \geq (2g+n-1)r_1r_2-(r_1r_2+1)\geq 1$,
since we assume $rn-2r-2>0$ if $g=0$,
$n\geq 2$ if $g=1$ and $n\geq 1$ if $g\geq 2$.
In particular, we have $\dim\RH_{(x,\blambda)}^{-1}([\rho])\geq 1$.

Next we consider the case $r\geq 3$ or $n\geq 2$ and
no assumption on the reducible representation $\rho$.
We may assume that $\rho$ is an extension of two representations
$\rho_1,\rho_2$ and denote the dimension of $\rho_i$ by $r_i$.
For each $p=1,2$, there is an $\balpha_p$-stable
$(\bt,\bmu_p)$-parabolic connection
$(E_p,\nabla_{E_p},\{(l_p)^{(i)}_j\})$ for some weight $\balpha_p$
such that
$\RH_{(x,\bmu_p)}(E_p,\nabla_{E_p},\{(l_p)^{(i)}_j\})=[\rho_p]$
and $0\leq\mathrm{Re}((\mu_p)^{(i)}_j)<1$ for any $i,j$.
We put $d_p:=\deg E_p$.
If we put $\Delta:=\{z\in\C| |z|<\epsilon \}$ ($\epsilon>0$),
then we can take morphisms
\[
 m_p:\Delta\longrightarrow \Lambda^{(n)}_{r_p}(d_p) \quad (p=1,2)
\]
such that $m_p(0)=\bmu_p$ and
for $m_p(t)=\{(\mu'_p)^{(i)}_j\}$ ($t\in\Delta\setminus\{0\}$),
we have $(\mu'_p)^{(i)}_j\neq(\mu'_{p'})^{(i)}_{j'}$
for $(p,i,j)\neq(p',i,j')$.
Replacing $\Delta$ by a neighborhood of $0$, we can take a morphism
$\varphi_p:\Delta\longrightarrow M^{\balpha_p}_C(\bt,r_p,d_p)$
which is a lift of $m_p$ and satisfies
$\varphi_p(0)=(E_p,\nabla_{E_p},\{(l_p)^{(i)}_j\})$.
The composite $\RH\circ\varphi_p$ determines a family $\tilde{\rho}_p$
of representations over $\Delta$ and
$\tilde{\rho}_1\oplus\tilde{\rho}_2$ determines a morphism
$\varphi:\Delta\rightarrow RP_r(C,\bt)$ such that
$\varphi(0)=[\rho]$ and for $t\in\Delta\setminus\{0\}$,
the representation corresponding to $\varphi(t)$ satisfies
the condition considered in the former case.
Then we have $\dim\RH^{-1}(\varphi(t))\geq 1$
for any $t\in\Delta\setminus\{0\}$ by the proof of the former case.
Since the base change
$\RH_{\Delta}:M^{\balpha}_C(\bt,r,d)_{\Delta}\rightarrow \Delta$
is a proper morphism,
we can see by the upper semi-continuity of the fiber dimension
of a holomorphic mapping that
$\dim(\RH_{\Delta})^{-1}(0)\geq 1$,
which means $\dim\RH_{(x,\blambda)}^{-1}([\rho])\geq 1$.

Finally assume that $r=2$ and $n=1$.
Let $\rho$ be a reducible representation such that
$[\rho]=[\rho_1\oplus\rho_2]$ in $RP_2(C,\bt)$.
In this case we have $g\geq 2$ by the assumption of Theorem \ref{main-thm}.
By Proposition \ref{surjectivity}, we can take an $\balpha$-stable
$(\bt,\blambda)$-parabolic connection $(E,\nabla,\{l^{(i)}_j\})$ such that
$\RH(E,\nabla,\{l^{(i)}_j\})=[\rho]$.
We put $\sigma(E,\nabla,\{l^{(i)}_j\})=(E',\nabla',\{(l')^{(i)}_j\})$.
Note that the local monodromy $\rho_i(\gamma_1)$ is trivial for $i=1,2$,
because $\rho_i$ is a one dimensional representation and $n=1$.
So the representation $\rho_{E'}$ corresponding to
$\ker(\nabla')^{an}|_{C\setminus\{t_1,\ldots,t_n\}}$
can be given by the following data:
\begin{gather*}
 \rho_{E'}(\alpha_k)=
 \begin{pmatrix}
  \rho_1(\alpha_k) & a^{(0)}_k \\
  0 & \rho_2(\alpha_k)
 \end{pmatrix}
 \quad
 \rho_{E'}(\beta_k)=
 \begin{pmatrix}
  \rho_1(\beta_k) & b^{(0)}_k \\
  0 & \rho_2(\beta_k)
 \end{pmatrix} 
 \quad (1\leq k\leq g).
\end{gather*}
Note that the monodromy matrix
\[
 \rho_{E'}(\gamma_1)=
 \begin{pmatrix}
  1 & c^{(0)} \\
  0 & 1
 \end{pmatrix}
\]
is given by
\[
 \rho_{E'}(\beta_g)\rho_{E'}(\alpha_g)\rho_{E'}(\beta_g)^{-1}\rho_{E'}(\alpha_g)^{-1}
 \cdots\rho_{E'}(\beta_1)\rho_{E'}(\alpha_1)\rho_{E'}(\beta_1)^{-1}\rho_{E'}(\alpha_1)^{-1}.
\]
By the upper semi-continuity of fiber dimension,
we may assume that $\rho_1\not\cong\rho_2$.
If $c^{(0)}\neq 0$, then we can prove in the same manner
as the first case that
$\dim_{(E,\nabla,\{l^{(i)}_j\})}\RH_{(x,\blambda)}^{-1}([\rho])\geq 2g-2>0$.
If $c^{(0)}=0$, then we may assume by the upper semi-continuity
of fiber dimension that
$l^{(1)}_1\neq \left\{\genfrac{(}{)}{0pt}{}{*}{0}\right\}$.
Consider the family of representations  given by matrices
\begin{gather*}
 \tilde{\rho}(\alpha_k)=
 \begin{pmatrix}
  \rho_1(\alpha_k) & a^{(0)}_k+a_k \\
  0 & \rho_2(\alpha_k)
 \end{pmatrix}
 \quad
 \tilde{\rho}(\beta_k)=
 \begin{pmatrix}
  \rho_1(\beta_k) & b^{(0)}_k+b_k \\
  0 & \rho_2(\beta_k)
 \end{pmatrix}
 \quad (1\leq k\leq g),
\end{gather*}
where $a_k,b_k$ ($1\leq k\leq g$) move around in
$\{c_k=0\}\subset\C^{2g}$.
Here $c_k$ is given by
\[
 \begin{pmatrix}
  1 & c_k \\
  0 & 1
 \end{pmatrix}
 =\tilde{\rho}(\beta_g)\tilde{\rho}(\alpha_g)
 \tilde{\rho}(\beta_g)^{-1}\tilde{\rho}(\alpha_g)^{-1}\cdots
 \tilde{\rho}(\beta_{1})\tilde{\rho}(\alpha_1)
 \tilde{\rho}(\beta_1)^{-1}\tilde{\rho}(\alpha_1)^{-1}.
\]
Adding parabolic structure, we can obtain a family
$(\tilde{E}',\tilde{\nabla}',\{(\tilde{l}')^{(i)}_j\})$
of $(\bt,\bmu)$-parabolic connections over a variety $Y$
of dimension at least $2g$.
If we put
\[
 Y':=\left\{y\in Y\left|
 \text{$\sigma^{-1}(\tilde{E}',\tilde{\nabla}',\{(\tilde{l}')^{(i)}_j\})
 \otimes k(y)$ is $\balpha$-stable}
 \right\}\right.
\]
we obtain a morphism $f':Y'\rightarrow M^{\balpha}_C(\bt,\blambda)$.
We can see that
$f'(y_1)=f'(y_2)$ if and only if there is an isomorphism
$\varphi:(\tilde{E}',\tilde{\nabla}')\otimes k(y_1)
\stackrel{\sim}\rightarrow
(\tilde{E}',\tilde{\nabla}')\otimes k(y_2)$
such that
$\varphi((\tilde{l}')^{(i)}_j\otimes k(y_1))=
(\tilde{l}')^{(i)}_j\otimes k(y_2)$.
The isomorphism $\varphi$ corresponds to the isomorphisms between
the representations
$\tilde{\rho}_{y_1}$, $\tilde{\rho}_{y_2}$
which is given by an element of the group
\[
 \left.\left\{
 \begin{pmatrix}
  c & a \\
  0 & 1
 \end{pmatrix}
 \right| c\in\C^{\times}, \; a\in\C \right\}.
\]
So we have $\dim f'(Y')= \dim Y'-2 \geq 2g-2>0$.
Since $f'(Y')\subset\RH_{(x,\blambda)}^{-1}([\rho])$, we have
$\dim\RH_{(x,\blambda)}^{-1}([\rho])>0$.

From all the above argument, we have in any case that
$\dim\RH^{-1}_{(x,\blambda)}([\rho])\geq 1$
for any $[\rho]\in RP_r(C,\bt)_{\ba}^{\mathrm{sing}}$.
Now take any point $p\in RP_r(C,\bt)_{\ba}^{\mathrm{sing}}$.
We want to show that $p$ is a singular point of $RP_r(C,\bt)_{\ba}$.
Assume that $p$ is a non-singular point of $RP_r(C,\bt)_{\ba}$.
By Proposition \ref{non-deg-form}, there is a non-degenerate
$2$-form $\omega$ on $M^{\balpha}_C(\bt,\blambda)$.
Via an isomorphism
$M^{\balpha}_C(\bt,\blambda)^{\sharp}\stackrel{\RH}\rightarrow
RP_r(C,\bt)_{\ba}^{\sharp}$,
we can obtain a non-degenerate $2$-form $\omega'$ on
$RP_r(C,\bt)_{\ba}^{\sharp}$ such that
$\RH_{(x,\blambda)}^*(\omega')=\omega|_{M^{\balpha}_C(\bt,\blambda)^{\sharp}}$.
Then $\omega'$ can be extended to a form defined also in a neighborhood of $p$
because $\codim_{RP_r(C,\bt)_{\ba}}(RP_r(C,\bt)_{\ba}^{\mathrm{sing}})\geq 2$
and we are assuming that $p$ is a non-singular point.
Since $\dim_x\RH_{(x,\blambda)}^{-1}(p)\geq 1$ for a point
$x\in\RH_{(x,\blambda)}^{-1}(p)$,
there is a tangent vector
$0\neq v\in \Theta_{M^{\balpha}_C(\bt,\blambda),x}$
such that $(\RH_{(x,\blambda)})_*(v)=0$.
For any tangent vector $w\in \Theta_{M^{\balpha}_C(\bt,\blambda),x}$,
we have
\[
 \omega(v,w)=(\RH_{(x,\blambda)}^*\omega')(v,w)
 =\omega'((\RH_{(x,\blambda)})_*(v),(\RH_{(x,\blambda)})_*(w))=0.
\]
Since $\omega$ determines a non-degenerate pairing on the tangent space
$\Theta_{M^{\balpha}_C(\bt,\blambda),x}$,
we have $v=0$ which is a contradiction.
Thus $p$ is a singular point of $RP_r(C,\bt)_{\ba}$.
\end{proof}

\section{Symplectic structure on the moduli space}\label{symplectic structure}

\begin{Proposition}\label{symplectic-form}
 Put $(C,\bt)=(\cC_x,\tilde{\bt}_x)$ for a point $x\in T$ and
 take a point $\blambda\in\Lambda_r^{(n)}(d)$.
 Then the moduli space $M_C^{\balpha}(\bt,\blambda)$
 has a holomorphic symplectic structure
 (more precisely an algebraic symplectic structure).
\end{Proposition}

We can obtain the above Proposition by the following two
propositions.

\begin{Proposition}\label{non-deg-form}
There is a nondegenerate relative 2-form
$\omega\in H^0(M^{\balpha}_{\cC/T}(\tilde{\bt},r,d),
\Omega^2_{M^{\balpha}_{\cC/T}(\tilde{\bt},r,d)/
T\times\Lambda^{(n)}_r(d)})$.
\end{Proposition}

\begin{proof}
Take a universal family $(\tilde{E},\tilde{\nabla},\{\tilde{l}^{(i)}_j\})$
on $\cC\times_T M^{\balpha}_{\cC/T}(\tilde{\bt},r,d)$.
We define a complex $\cF^{\bullet}$ by
\begin{align*}
 \cF^0 &=\left\{ s\in {\mathcal End}(\tilde{E}) \left|
 s|_{\tilde{t}_i\times M^{\balpha}_{\cC/T}(\tilde{\bt},r,d)}
 (\tilde{l}^{(i)}_j)\subset\tilde{l}^{(i)}_j \;
 \text{for any $i,j$}
 \right\}\right. \\
 \cF^1 &=\left\{
 s\in{\mathcal End}(\tilde{E})\otimes\Omega^1_{\cC/T}(D(\tilde{\bt}))
 \left|
 \text{$\res_{\tilde{t}_i\times M^{\balpha}_{\cC/T}(\tilde{\bt},r,d)}(s)
 (\tilde{l}^{(i)}_j)\subset\tilde{l}^{(i)}_{j+1}$ for any $i,j$}
 \right\}\right. \quad \text{and} \\
 & \nabla_{\cF^{\bullet}} :\cF^0\ni s\mapsto
 \tilde{\nabla}\circ s-s\circ\tilde{\nabla}\in\cF^1.
\end{align*}
Notice that there is an isomorphism of sheaves
\[
 \R^1(\pi_{M^{\balpha}_{\cC/T}(\tilde{\bt},r,d)})_*(\cF^{\bullet})
 \stackrel{\sim}\longrightarrow
 \Theta_{M^{\balpha}_{\cC/T}(\tilde{\bt},r,d)/T\times\Lambda^{(n)}_r(d)},
\]
where
$\pi_{M^{\balpha}_{\cC/T}(\tilde{\bt},r,d)}:
\cC\times_TM^{\balpha}_{\cC/T}(\tilde{\bt},r,d)\ra
M^{\balpha}_{\cC/T}(\tilde{\bt},r,d)$
is the projection and
$\Theta_{M^{\balpha}_{\cC/T}(\tilde{\bt},r,d)/T\times\Lambda^{(n)}_r(d)}$
is the relative algebraic tangent bundle on the moduli space.
For each affine open subset
$U\subset M^{\balpha}_{\cC/T}(\tilde{\bt},r,d)$,
we define a pairing
\begin{gather*}
 \bH^1(\cC\times_TU,\cF^{\bullet}_U)\otimes\bH^1(\cC\times_TU,\cF^{\bullet}_U)
 \longrightarrow
 \bH^2(\cC\times_TU,\Omega^{\bullet}_{\cC\times_TU/U})\cong H^0(\cO_U) \\
 [\{u_{\alpha\beta}\},\{v_{\alpha}\}]\otimes
 [\{u'_{\alpha\beta}\},\{v'_{\alpha}\}] \mapsto
 [\{\Tr(u_{\alpha\beta}\circ u'_{\beta\gamma})\}-\{\Tr(u_{\alpha\beta}\circ v'_{\beta})
 -\Tr(v_{\alpha}\circ u'_{\alpha\beta})\}],
\end{gather*}
where we consider in \u{C}ech cohomology with respect to an
affine open covering $\{U_{\alpha}\}$ of $\cC\times_TU$,
$\{u_{\alpha\beta}\}\in C^1(\cF^0)$,
$\{v_{\alpha}\}\in C^0(\cF^1)$ and so on.
This pairing determines a pairing
\begin{equation}
 \omega:\R^1(\pi_{M^{\balpha}_{\cC/T}(\tilde{\bt},r,d)})_*(\cF^{\bullet})
 \otimes\R^1(\pi_{M^{\balpha}_{\cC/T}(\tilde{\bt},r,d)})_*(\cF^{\bullet})
 \longrightarrow \cO_{M^{\balpha}_{\cC/T}(\tilde{\bt},r,d)}.
\end{equation}
Take any point $x\in M^{\balpha}_{\cC/T}(\tilde{\bt},r,d)$
which corresponds to a $(\bt,\blambda)$-parabolic connection
$(E,\nabla,\{l^{(i)}_j\})$.
A tangent vector
$v\in \Theta_{M^{\balpha}_{\cC/T}(\tilde{\bt},r,d)/
T\times\Lambda^{(n)}_r(d)}|_x$
corresponds to a flat family of $(\bt,\blambda)$-parabolic connections
$(\tilde{E},\tilde{\nabla},\{\tilde{l}^{(i)}_j\})$
on $\cC_x\times\Spec\C[\epsilon]/(\epsilon^2)$ over $\C[\epsilon]/(\epsilon^2)$
which is a lift of $(E,\nabla,\{l^{(i)}_j\})$.
We can see that
$\omega|_x(v\otimes v)\in \bH^0(\Omega^{\bullet}_{\cC_x})
\xrightarrow[\sim]{\Tr^{-1}}\bH^2(\cF^{\bullet}_x)$
is just the obstruction class for the lifting of
$(\tilde{E},\tilde{\nabla},\{\tilde{l}^{(i)}_j\})$
to a flat family over $\C[\epsilon]/(\epsilon^3)$.
Since the moduli space $M^{\balpha}_{\cC_x}(\bt,\blambda)$ is smooth,
this obstruction class vanishes: $\omega|_x(v\otimes v)=0$.
Thus $\omega$ is skew symmetric and determines a $2$-form.
For each point $x\in M^{\balpha}_{\cC/T}(\tilde{\bt},r,d)$, let
\[
 \bH^1(\cF^{\bullet}_x)\stackrel{\xi}\longrightarrow
 \bH^1(\cF^{\bullet}_x)^{\vee}
\]
be the homomorphism induced by $\omega$.
Then we have the following exact commutative diagram
\[
 \begin{CD}
  H^0(\cF^0_x) @>>> H^0(\cF^1_x) @>>> \bH^1(\cF^{\bullet}_x)
  @>>> H^1(\cF^0_x) @>>> H^1(\cF^1_x) \\
  @Vb_1VV @Vb_2VV @V\xi VV @Vb_3VV @Vb_4VV \\
  H^1(\cF^1_x)^{\vee} @>>> H^1(\cF^0_x)^{\vee} @>>>
  \bH^1(\cF^{\bullet}_x)^{\vee} @>>> H^0(\cF^1_x)^{\vee}
  @>>> H^0(\cF^0_x)^{\vee},
 \end{CD}
\]
where $b_1,\ldots,b_4$ are isomorphisms induced by the isomorphisms
$\cF^0\cong(\cF^1)^{\vee}\otimes\Omega^1_{\cC/T}$,
$\cF^1\cong(\cF^0)^{\vee}\otimes\Omega^1_{\cC/T}$ and Serre duality.
So $\xi$ becomes an isomorphism by five lemma.
Hence $\omega$ is non-degenerate.
\end{proof}

\begin{Proposition}\label{d-closed}
 For the $2$-form constructed in Proposition \ref{non-deg-form},
 we have $d\omega=0$.
\end{Proposition}

\begin{proof}
It suffices to show that the restriction of $\omega$ to
any general fiber $M^{\balpha}_{\cC_x}(\bt,\blambda)$ over
$T\times\Lambda^{(n)}_r(d)$ is $d$-closed.
We may assume $\blambda$ generic, and so
$\lambda^{(i)}_j-\lambda^{(i)}_{j'}\not\in\mathbf{Z}$
for any $i$ and $j\neq j'$.
We also denote the restriction by $\omega$.
Let $(\tilde{E},\tilde{\nabla},\{\tilde{l}^{(i)}_j\})$
be a universal family on $\cC_x\times M^{\balpha}_{\cC_x}(\bt,\blambda)$
over $M^{\balpha}_{\cC_x}(\bt,\blambda)$.
We set
\begin{gather*}
 \cF^0:=\left\{ s\in {\mathcal End}(\tilde{E}) \left|
 s|_{t_i\times M^{\balpha}_{\cC_x}(\bt,\blambda)}(\tilde{l}^{(i)}_j)
 \subset \tilde{l}^{(i)}_j \quad \text{for any $i,j$} \right\}\right. \\
 \cF^1:=\left\{ s\in {\mathcal End}(\tilde{E})
 \otimes\Omega^1_{\cC_x}(t_1+\cdots+t_n) \left|
 \res_{t_i\times M^{\balpha}_{\cC_x}(\bt,\blambda)}(s)(\tilde{l}^{(i)}_j)
 \subset \tilde{l}^{(i)}_{j+1} \quad \text{for any $i,j$} \right\}\right. \\
 \nabla_{\cF^{\bullet}}:\cF^0\longrightarrow \cF^1; \quad
 \nabla_{\cF^{\bullet}}(s)=\tilde{\nabla}\circ s-s\circ\tilde{\nabla}. \\
\end{gather*}
We also put
\[
 \V:=\ker\tilde{\nabla}^{an}|_{(\cC_x\setminus\{t_1,\ldots,t_n\})
 \times M^{\balpha}_{\cC_x}(\bt,\blambda)}.
\]
Since there is a canonical isomorphism
\[
 \ker\nabla_{\cF^{\bullet}}^{an}|_{(\cC_x\setminus\{t_1,\ldots,t_n\})
 \times M^{\balpha}_{\cC_x}(\bt,\blambda)}
 \stackrel{\sim}\longrightarrow {\mathcal End}(\V),
\]
we obtain a canonical homomorphism
\begin{equation}\label{local-solution}
 \ker\nabla_{\cF^{\bullet}}^{an}|_{\cC_x\times
 M^{\balpha}_{\cC_x}(\bt,\blambda)}
 \longrightarrow j_*({\mathcal End}(\V)),
\end{equation}
where
$j:(\cC_x\setminus\{t_1,\ldots,t_n\})\times M^{\balpha}_{\cC_x}(\bt,\blambda)
\hookrightarrow \cC_x\times M^{\balpha}_{\cC_x}(\bt,\blambda)$
is the canonical inclusion.
We will show by local calculations that
the homomorphism (\ref{local-solution}) is in fact
an isomorphism.
We can see by [\cite{Deligne}, II, Proposition 5.4] that
there is an isomorphism
\[
 (\tilde{E}^{an}_{(t_i,y)},\tilde{\nabla}^{an}_{(t_i,y)})
 \cong
 \left(({\mathcal O}^{an}_{C\times M^{\balpha}_{\cC_x}(\bt,\blambda),
 (t_i,y)})^{\oplus r},\nabla^{\blambda}\right),
\]
on the stalks at $(t_i,y)\in\cC_x\times M^{\balpha}_{\cC_x}(\bt,\blambda)$,
where
$\nabla^{\blambda}:\left({\mathcal O}^{an}
_{\cC_x\times M^{\balpha}_{\cC_x}(\bt,\blambda),(t_i,y)}\right)^{\oplus r}
\rightarrow\left({\mathcal O}^{an}_{\cC_x\times
M^{\balpha}_{\cC_x}(\bt,\blambda),(t_i,y)}\right)^{\oplus r}
\otimes\Omega^1_{\cC_x\times
M^{\balpha}_{\cC_x}(\bt,\blambda)/M^{\balpha}_{\cC_x}(\bt,\blambda)}
(t_i\times M^{\balpha}_{\cC_x}(\bt,\blambda))$ is given by
\[
 \nabla^{\blambda}
 \begin{pmatrix}
  f_1 \\
  \vdots \\
  f_r
 \end{pmatrix}
 =
 \begin{pmatrix}
  df_1 \\
  \vdots \\
  df_r
 \end{pmatrix}
 +
 \begin{pmatrix}
  \lambda^{(i)}_{r-1}z_i^{-1}dz_i & 0 & 0 \\
  0 & \ddots & 0 \\
  0 & 0 & \lambda^{(i)}_0z_i^{-1}dz_i
 \end{pmatrix}
 \begin{pmatrix}
  f_1\\
  \vdots \\
  f_f
 \end{pmatrix},
\]
where $z_i$ is a generator of the maximal ideal of
${\mathcal O}_{\cC_x,t_i}$.
With respect to this isomorphism,
we have
\begin{align*}
  (\cF^0)^{an}_{(t_i,y)}&=
  \left.\left\{
   \begin{pmatrix}
    f_{1,1} & \cdots & f_{1,r} \\
    \vdots & \ddots & \vdots \\
    f_{r,1} & \cdots & f_{r,r}
   \end{pmatrix}
  \right|
  \begin{array}{l}
  \text{$f_{p,q}\in z_i{\mathcal O}^{an}_{\cC_x\times
  M^{\balpha}_{\cC_x}(\bt,\blambda),(t_i,y)}$ for $p>q$ and} \\
  \text{$f_{p,q}\in{\mathcal O}^{an}_{\cC_x\times
  M^{\balpha}_{\cC_x}(\bt,\blambda),(t_i,y)}$ for $p\leq q$}
  \end{array}
 \right\} \\
 (\cF^1)^{an}_{(t_i,y)}&=
 \left.\left\{
  \begin{pmatrix}
   g_{1,1} & \cdots & g_{1,r} \\
   \vdots & \ddots & \vdots \\
   g_{r,1} & \cdots & g_{r,r} 
  \end{pmatrix}\right|
  \begin{array}{l}
    \text{$g_{p,q}\in{\mathcal O}^{an}_{\cC_x\times
    M^{\balpha}_{\cC_x}(\bt,\blambda),(t_i,y)}dz_i$
    for $p\geq q$ and} \\
    \text{$g_{p,q}\in{\mathcal O}^{an}_{\cC_x\times
    M^{\balpha}_{\cC_x}(\bt,\blambda),(t_i,y)}dz_i/z_i$
    for $p<q$}
   \end{array} 
   \right\}
\end{align*}
and
$(\nabla_{\cF^{\bullet}})^{an}_{(t_i,y)}:(\cF^0)^{an}_{(t_i,y)}
\rightarrow (\cF^1)^{an}_{(t_i,y)}$ is given by
\begin{gather*}
 \begin{pmatrix}
  f_{1,1} & f_{1,2} & \cdots & f_{1,r} \\
  f_{2,1} & f_{2,2} & \cdots & f_{2,r} \\
  \vdots & \vdots & \ddots & \vdots \\
  f_{r,1} & f_{r,2} & \cdots & f_{r,r}
 \end{pmatrix}
 \mapsto
 \begin{pmatrix}
  df_{1,1} & df_{1,2} & \cdots & df_{1,r} \\
  df_{2,1} & df_{2,2} & \cdots & df_{2,r} \\
  \vdots & \vdots & \ddots & \vdots \\
  df_{r,1} & df_{r,2} & \cdots & df_{r,r}
 \end{pmatrix} \\
 +
 \begin{pmatrix}
  0 & (\lambda^{(i)}_{r-1}-\lambda^{(i)}_{r-2})f_{1,2}dz_i/z_i
  & \cdots & (\lambda^{(i)}_{r-1}-\lambda^{(i)}_0)f_{1,r}dz_i/z_i \\
  (\lambda^{(i)}_{r-2}-\lambda^{(i)}_{r-1})f_{2,1}dz_i/z_i & 0 &
  \cdots & (\lambda^{(i)}_{r-2}-\lambda^{(i)}_0)f_{2,r}dz_i/z_i \\
  \vdots & \vdots & \ddots & \vdots \\
  (\lambda^{(i)}_0-\lambda^{(i)}_{r-1})f_{r,1}dz_i/z_i &
  (\lambda^{(i)}_0-\lambda^{(i)}_{r-2})f_{r,2}dz_i/z_i &
  \cdots & 0
 \end{pmatrix}.
\end{gather*}
Take any element
\[
 \begin{pmatrix}
  f_{1,1} & f_{1,2} & \cdots & f_{1,r} \\
  f_{2,1} & f_{2,2} & \cdots & f_{2,r} \\
  \vdots & \vdots & \ddots & \vdots \\
  f_{r,1} & f_{r,2} & \cdots & f_{r,r}
 \end{pmatrix}
 \in j_*({\mathcal End}(\mathbf{V}))_{(t_i,y)}=
 j_*\left(\ker\nabla^{an}_{{\mathcal F}^{\bullet}}|
 _{(\cC_x\setminus\{t_1,\ldots,t_n\})\times
 M^{\balpha}_{\cC_x}(\bt,\blambda)}\right)_{(t_i,y)}
\]
For $j\neq k$, we have
$f_{j,k}=c_{j,k}z_i^{\lambda^{(i)}_{r-k}-\lambda^{(i)}_{r-j}}$
for some function
$c_{j,k}\in{\mathcal O}^{an}_{M^{\balpha}_{\cC_x}(\bt,\blambda),y}$
Since $f_{j,k}$ is a single valued function on
$(U\setminus\{t_i\})\times M$
for some open neighborhood $U$ of $t_i$ in $\cC_x$
and some open neighborhood $M$ of $y$ in
$M^{\balpha}_{\cC_x}(\bt,\blambda)$
and $\lambda^{(i)}_{r-k}-\lambda^{(i)}_{r-j}\notin\mathbf{Z}$,
we have $c_{j,k}=0$, namely $f_{j,k}=0$.
On the other hand, we have $df_{j,j}=0$
for $j=1,\ldots,r$.
So we have
$f_{j,j}\in{\mathcal O}^{an}_{M^{\balpha}_{\cC_x}(\bt,\blambda),y}$.
Thus we have
\[
 \begin{pmatrix}
  f_{1,1} & f_{1,2} & \cdots & f_{1,r} \\
  f_{2,1} & f_{2,2} & \cdots & f_{2,r} \\
  \vdots & \vdots & \ddots & \vdots \\
  f_{r,1} & f_{r,2} & \cdots & f_{r,r}
 \end{pmatrix}
 \in({\mathcal F}^0)^{an}_{(t_i,y)}\cap
 j_*(\ker\nabla^{an}_{{\mathcal F}^{\bullet}}
 |_{\cC_x\setminus\{t_1,\ldots,t_n\}})_{(t_i,y)}
 =\left(\ker\nabla^{an}_{{\mathcal F}^{\bullet}}
 |_{\cC_x\times M^{\balpha}_{\cC_x}(\bt,\blambda)}\right)_{(t_i,y)}
\]
So the homomorphism
(\ref{local-solution}) is surjective.
The injectivity of (\ref{local-solution}) is obvious.
Thus (\ref{local-solution}) is isomorphic.

Next we will show that
$\nabla_{\cF^{\bullet}}^{an}:(\cF^0)^{an}\rightarrow(\cF^1)^{an}$
is surjective.
For this it is essential to show the surjectivity
of the morphism
\[
 (\nabla_{\cF^{\bullet}})^{an}_{(t_i,y)}:
 (\cF^0)^{an}_{(t_i,y)}\longrightarrow
 (\cF^1)^{an}_{(t_i,y)}
 \quad ((t_i,y)\in\cC_x\times M_{\cC_x}^{\balpha}(\bt,\blambda))
\]
at the stalks.
Take any member
\[
 \begin{pmatrix}
  g_{1,1} & \cdots & g_{1,r} \\
  \vdots & \ddots & \vdots \\
  g_{r,1} & \cdots & g_{r,r}
 \end{pmatrix}
 \in(\cF^1)^{an}_{(t_i,y)}.
\]
For $p>q$ we can write
$g_{p,q}=\sum_{k=0}^{\infty}b_{p,q}^{(k)}z_i^kdz_i$,
where
$b_{p,q}^{(k)}\in{\mathcal O}^{an}_{M^{\balpha}_{\cC_x}
(\bt,\blambda),y}$
for $k\geq 0$.
We define a power series
$f_{p,q}=\sum_{k=0}^{\infty}a_{p,q}^{(k)}z_i^k$ by
$a_{p,q}^{(0)}:=0$ and
$a_{p,q}^{(k)}:=(k+\lambda^{(i)}_{r-p}-\lambda^{(i)}_{r-q})^{-1}
b_{p,q}^{(k-1)}$
for $k\geq 1$.
Since $\sum_{k=0}^{\infty}b_{p,q}^{(k)}z_i^k$
is a convergent power series,
we can see that
$\sum_{k=0}^{\infty}a_{p,q}^{(k)}z_i^k$
is also convergent.
So we have
$f_{p,q}:=\sum_{k=0}^{\infty}a_{p,q}^{(k)}z_i^k
\in{\mathcal O}^{an}_{\cC_x\times
M^{\balpha}_{\cC_x}(\bt,\blambda),(t_i,y)}$
and
$df_{p,q}+(\lambda^{(i)}_{r-p}-\lambda^{(i)}_{r-q})f_{p,q}dz_i/z_i
=g_{p,q}$.
For $g_{p,p}\in{\mathcal O}^{an}_{\cC_x\times
M^{\balpha}_{\cC_x}(\bt,\blambda),(t_i,y)}dz_i$,
we can find 
$f_{p,p}\in{\mathcal O}^{an}_{\cC_x\times
M^{\balpha}_{\cC_x}(\bt,\blambda),(t_i,y)}$
such that $df_{p,p}=g_{p,p}$.
For $p<q$, we can write
$g_{p,q}=\sum_{k=0}^{\infty}b_{p,q}^{(k)}z_i^{k-1}dz_i$
with $b_{p,q}^{(k)}\in{\mathcal O}^{an}_{M^{\balpha}_{\cC_x}
(\bt,\blambda),y}$
for $k\geq 0$.
For each $k\geq 0$, we put
$a_{p,q}^{(k)}:=(k+\lambda^{(i)}_{r-p}-\lambda^{(i)}_{r-q})^{-1}
b_{p,q}^{(k)}$.
Since $\sum_{k=0}^{\infty}b_{(p,q)}^{(k)}z_i^k$ is convergent,
we can see that
$f_{p,q}:=\sum_{k=0}^{\infty}a_{p,q}^{(k)}z_i^k$
is also convergent.
So we have
$f_{p,q}\in{\mathcal O}^{an}_{\cC_x\times
M^{\balpha}_{\cC_x}(\bt,\blambda),(t_i,y)}$
and
$df_{p,q}+(\lambda^{(i)}_{r-p}-\lambda^{(i)}_{r-q})f_{p,q}dz_i/z_i
=g_{p,q}$.
By construction, we have
\[
 \begin{pmatrix}
  f_{1,1} & \cdots & f_{1,r} \\
  \vdots & \ddots & \vdots \\
  f_{r,1} & \cdots & f_{r,r}
 \end{pmatrix}
 \in(\cF^0)^{an}_{(t_i,y)}
\]
and
\[
 (\nabla_{{\mathcal F}^{\bullet}})^{an}_{(t_i,y)}
 \begin{pmatrix}
  f_{1,1} & \cdots & f_{1,r} \\
  \vdots & \ddots & \vdots \\
  f_{r,1} & \cdots & f_{r,r}
 \end{pmatrix}
 =
 \begin{pmatrix}
  g_{1,1} & \cdots & g_{1,r} \\
  \vdots & \ddots & \vdots \\
  g_{r,1} & \cdots & g_{r,r}
 \end{pmatrix}.
\]
Hence we obtain the surjectivity of
$\nabla_{\cF^{\bullet}}^{an}:(\cF^0)^{an}\rightarrow(\cF^1)^{an}$.

So we obtain an isomorphism
\[
 \R^1(\pi_{M^{\balpha}_{\cC_x}(\bt,\blambda)})_*((\cF^{\bullet})^{an})
 \stackrel{\sim}\longrightarrow
 R^1(\pi_{M^{\balpha}_{\cC_x}(\bt,\blambda)})_*(j_*({\mathcal End}(\V)))
\]
of analytic sheaves.
Then $\omega$ corresponds to a canonical pairing
\begin{gather*}
 R^1(\pi_{M^{\balpha}_{\cC_x}(\bt,\blambda)})_*(j_*({\mathcal End}(\V)))\times
 R^1(\pi_{M^{\balpha}_{\cC_x}(\bt,\blambda)})_*(j_*({\mathcal End}(\V)))
 \longrightarrow
 R^2(\pi_{M^{\balpha}_{\cC_x}(\bt,\blambda)})_*
 (\pi_{M^{\balpha}_{\cC_x}(\bt,\blambda)}^{-1}
 (\cO_{M^{\balpha}_{\cC_x}(\bt,\blambda)})) \\
 ([\{c_{\alpha\beta}\}],[\{c'_{\alpha\beta}\}])\mapsto
 [\{\Tr(c_{\alpha\beta}\circ c'_{\beta\gamma})\}].
\end{gather*}
Note that
$R^2(\pi_{M^{\balpha}_{\cC_x}(\bt,\blambda)})_*
 (\pi_{M^{\balpha}_{\cC_x}(\bt,\blambda)}^{-1}
 (\cO_{M^{\balpha}_{\cC_x}(\bt,\blambda)}))\cong
\cO_{M^{\balpha}_{\cC_x}(\bt,\blambda)}$.

Take any small open subset $M\subset M^{\balpha}_{\cC_x}(\bt,\blambda)$
and an open covering $\{U_{\alpha}\}$ of $\cC_x$ such that
$\sharp\{\alpha | t_i\in U_{\alpha} \}=1$ for any $i$ and that
$\sharp(U_{\alpha}\cap\{t_1,\ldots,t_n\})\leq 1$ for any $\alpha$.
If we replace $U_{\alpha}$ and $M$ sufficiently smaller,
there exist a sheaf $E_{\alpha}$ on
$U_{\alpha}$ such that
$E_{\alpha}|_{U_{\alpha}\cap U_{\beta}}\cong
\C^{\oplus r^2}_{U_{\alpha}\cap U_{\beta}}$
for any $\beta\neq\alpha$ and an isomorphism
$\phi_{\alpha}:j_*(\V)|_{U_{\alpha}\times M}
\stackrel{\sim}\rightarrow E_{\alpha}\otimes\cO_M$.
For each $\alpha,\beta$, we put
\[
 \varphi_{\alpha\beta}:=\phi_{\alpha}\circ\phi_{\beta}^{-1}:
 E_{\beta}\otimes\cO_M|_{(U_{\alpha}\cap U_{\beta})\times M}
 \stackrel{\phi_{\beta}^{-1}}\longrightarrow
 j_*(\V)|_{(U_{\alpha}\cap U_{\beta})\times M}
 \stackrel{\phi_{\alpha}}\longrightarrow
 E_{\alpha}\otimes\cO_M|_{(U_{\alpha}\cap U_{\beta})\times M}.
\]
So $j_*(\V)|_{\cC_x\times M}$
is given by the transition functions $\{\varphi_{\alpha\beta}\}$.
Next we consider a vector field $v\in H^0(M,\Theta_M)$.
Then $v$ corresponds to a derivation $D_v:\cO_M\rightarrow\cO_M$
which naturally induces a morphism
\[
 D_v:{\mathcal Hom}(E_{\beta}|_{U_{\alpha}\cap U_{\beta}},
 E_{\alpha}|_{U_{\alpha}\cap U_{\beta}})\otimes\cO_M
 \longrightarrow
 {\mathcal Hom}(E_{\beta}|_{U_{\alpha}\cap U_{\beta}},
 E_{\alpha}|_{U_{\alpha}\cap U_{\beta}})\otimes\cO_M.
\]
We note that $v$ also corresponds to a morphism
$f_v:\Spec\cO_M[\epsilon]\rightarrow M$,
where $\epsilon^2=0$.
Then $(1_{\cC_x}\times f_v)^*(j_*(\V))$ corresponds to the glueing data
$\{(1_{\cC_x}\times f_v)^*(\varphi_{\alpha\beta})\}$.
We can see the following equality:
\begin{gather*}
 (1_{\cC_x}\times f_v)^*(\varphi_{\alpha\beta})=
 \varphi_{\alpha\beta}+\epsilon D_v(\varphi_{\alpha\beta}):
 E_{\beta}|_{U_{\alpha}\cap U_{\beta}}\otimes\cO_M[\epsilon]
 \stackrel{\sim}\longrightarrow
 E_{\alpha}|_{U_{\alpha}\cap U_{\beta}}\otimes\cO_M[\epsilon] \\
 \hspace{200pt} a+\epsilon b \mapsto
 \varphi_{\alpha\beta}(a)+
 \epsilon(\varphi_{\alpha\beta}(b)+D_v(\varphi_{\alpha\beta})(a)).
\end{gather*}
So the isomorphism
$\Theta_M\cong R^1(\pi_M)_*(j_*({\mathcal End}(\V))_M)$
is given by
\[
 \Theta_M\ni v \mapsto [\{
 \phi_{\alpha}^{-1}\circ D_v(\varphi_{\alpha\beta})\circ\phi_{\beta}
 \}]\in R^1(\pi_M)_*(j_*({\mathcal End}(\V_M))).
\]
and
\[
 \omega(u,v)=[\{
 \Tr(D_u(\varphi_{\alpha\beta})\circ D_v(\varphi_{\beta\gamma})
 \circ\varphi_{\gamma\alpha})\}]
 \in R^2(\pi_M)_*(\pi_M^{-1}(\cO_M))\cong \cO_M.
\]
Thus we have
\begin{align*}
 d\omega(u,v,w)&=D_u(\omega(v,w))+D_v(\omega(w,u))+D_w(\omega(u,v))
  +\omega(w,[u,v])+\omega([u,w],v)+\omega(u,[v,w]) \\
 &=\Big\{ \Tr\Big(
 D_u(D_v(\varphi_{\alpha\beta})\circ D_w(\varphi_{\beta\gamma})
 \circ\varphi_{\gamma\alpha})
 +D_v(D_w(\varphi_{\alpha\beta})\circ D_u(\varphi_{\beta\gamma})
 \circ\varphi_{\gamma\alpha}) \\
 &\quad +D_w(D_u(\varphi_{\alpha\beta})\circ D_v(\varphi_{\beta\gamma})
 \circ\varphi_{\gamma\alpha})
 +D_w(\varphi_{\alpha\beta})\circ(D_uD_v-D_vD_u)(\varphi_{\beta\gamma})
 \circ\varphi_{\gamma\alpha} \\
 & \quad +(D_uD_w-D_wD_u)(\varphi_{\alpha\beta})\circ
 D_v(\varphi_{\beta\gamma})\circ\varphi_{\gamma\alpha}
 +D_u(\varphi_{\alpha\beta})\circ(D_vD_w-D_wD_v)(\varphi_{\beta\gamma})
 \circ\varphi_{\gamma\alpha}
 \Big)\Big\} \\
 &=\Big\{\Tr\Big(
 D_uD_v(\varphi_{\alpha\beta})\circ D_w(\varphi_{\beta\gamma})
 \circ\varphi_{\gamma\alpha}
 +D_w(\varphi_{\alpha\beta})\circ D_uD_v(\varphi_{\beta\gamma})
 \circ\varphi_{\gamma\alpha} \\
 & \quad +D_uD_w(\varphi_{\alpha\beta})\circ D_v(\varphi_{\beta\gamma})
 \circ\varphi_{\gamma\alpha}
 +D_v(\varphi_{\alpha\beta})\circ D_uD_w(\varphi_{\beta\gamma})
 \circ\varphi_{\gamma\alpha} \\
 & \quad +D_vD_w(\varphi_{\alpha\beta})\circ D_u(\varphi_{\beta\gamma})
 \circ\varphi_{\gamma\alpha}
 +D_u(\varphi_{\alpha\beta})\circ D_vD_w(\varphi_{\beta\gamma})
 \circ\varphi_{\gamma\alpha} \\
 & \quad +D_u(\varphi_{\alpha\beta})\circ D_v(\varphi_{\beta\gamma})
 \circ D_w(\varphi_{\gamma\alpha})
 +D_v(\varphi_{\alpha\beta})\circ D_w(\varphi_{\beta\gamma})
 \circ D_u(\varphi_{\gamma\alpha})   \\
 & \quad +D_w(\varphi_{\alpha\beta})\circ D_u(\varphi_{\beta\gamma})
 \circ D_v(\varphi_{\gamma\alpha})
 \Big)\Big\}.
\end{align*}
On the other hand, applying $D_uD_vD_w$ to
$\varphi_{\alpha\beta}\circ\varphi_{\beta\gamma}=\varphi_{\alpha\gamma}$,
we obtain the equality
\begin{align*}
 & -D_uD_vD_w(\varphi_{\alpha\beta})\varphi_{\beta\gamma}
 +D_uD_vD_w(\varphi_{\alpha\gamma})
 -\varphi_{\alpha\beta}D_uD_vD_w(\varphi_{\beta\gamma}) \\
 & =D_uD_v(\varphi_{\alpha\beta})\circ D_w(\varphi_{\beta\gamma})
 +D_w(\varphi_{\alpha\beta})\circ D_uD_v(\varphi_{\beta\gamma}) \\
 & \quad +D_uD_w(\varphi_{\alpha\beta})\circ D_v(\varphi_{\beta\gamma})
 +D_v(\varphi_{\alpha\beta})\circ D_uD_w(\varphi_{\beta\gamma}) \\
 & \quad +D_vD_w(\varphi_{\alpha\beta})\circ D_u(\varphi_{\beta\gamma})
 +D_u(\varphi_{\alpha\beta})\circ D_vD_w(\varphi_{\beta\gamma})
\end{align*}
which implies
\begin{align*}
 -d\{\Tr(D_uD_vD_w(\varphi_{\alpha\beta})\varphi_{\beta\alpha})\}
 &= \Big\{ \Tr\Big(
 D_uD_v(\varphi_{\alpha\beta})\circ D_w(\varphi_{\beta\gamma})
 \circ\varphi_{\gamma\alpha}
 +D_w(\varphi_{\alpha\beta})\circ D_uD_v(\varphi_{\beta\gamma})
 \circ\varphi_{\gamma\alpha} \\
 &\quad +D_uD_w(\varphi_{\alpha\beta})\circ D_v(\varphi_{\beta\gamma})
 \circ\varphi_{\gamma\alpha}
 +D_v(\varphi_{\alpha\beta})\circ D_uD_w(\varphi_{\beta\gamma})
 \circ\varphi_{\gamma\alpha} \\
 & \quad +D_vD_w(\varphi_{\alpha\beta})\circ D_u(\varphi_{\beta\gamma})
 \circ\varphi_{\gamma\alpha}
 +D_u(\varphi_{\alpha\beta})\circ D_vD_w(\varphi_{\beta\gamma}))
 \circ\varphi_{\gamma\alpha}
 \Big) \Big\}.
\end{align*}
We can also check the equality
\begin{align*}
 & \Tr\Big(
 D_u(\varphi_{\alpha\beta})\circ D_v(\varphi_{\beta\gamma})
 \circ D_w(\varphi_{\gamma\alpha})
 +D_v(\varphi_{\alpha\beta})\circ D_w(\varphi_{\beta\gamma})
 \circ D_u(\varphi_{\gamma\alpha}) \\
 & \hspace{200pt}
 +D_w(\varphi_{\alpha\beta})\circ D_u(\varphi_{\beta\gamma})
 \circ D_v(\varphi_{\gamma\alpha}) \Big) \\
 & = \Tr\Big(
 -D_u(\varphi_{\alpha\beta})D_v(\varphi_{\beta\alpha})
 D_w(\varphi_{\alpha\beta})\varphi_{\beta\alpha}
 +D_u(\varphi_{\alpha\gamma})D_v(\varphi_{\gamma\alpha})
 D_w(\varphi_{\alpha\gamma})\varphi_{\gamma\alpha} \\
 & \hspace{200pt}
 -\varphi_{\alpha\beta}D_u(\varphi_{\beta\gamma})D_v(\varphi_{\gamma\beta})
 D_w(\varphi_{\beta\gamma})\varphi_{\gamma\alpha}
 \Big)
\end{align*}
which means
\begin{align*}
 &-d\{ \Tr( D_u(\varphi_{\alpha\beta})D_v(\varphi_{\beta\alpha})
 D_w(\varphi_{\alpha\beta})\varphi_{\beta\alpha}) \} \\
 & = \Big\{ \Tr\Big(
 D_u(\varphi_{\alpha\beta})\circ D_v(\varphi_{\beta\gamma})
 \circ D_w(\varphi_{\gamma\alpha})
  +D_v(\varphi_{\alpha\beta})\circ D_w(\varphi_{\beta\gamma})
 \circ D_u(\varphi_{\gamma\alpha})  \\
 & \hspace{200pt}
 +D_w(\varphi_{\alpha\beta})\circ D_u(\varphi_{\beta\gamma})
 \circ D_v(\varphi_{\gamma\alpha})
 \Big) \Big\}.
\end{align*}
Thus we can see that
$d\omega(u,v,w)=0$.
\end{proof}

\begin{Remark}\rm
(1) Put $(C,\bt)=(\cC_x,\tilde{\bt}_x)$ for a point $x\in T$
and take $\blambda\in\Lambda^{(n)}_r(d)$.
Let $(L,\nabla_L)$ be a pair of a line bundle $L$ on $C$
and a connection $\nabla_L:L\rightarrow L\otimes\Omega^1_C(D(\bt))$
such that $\res_{t_i}(\nabla_L)=\sum_{j=0}^{r-1}\lambda^{(i)}_j$.
Consider the moduli space $M_C^{\balpha}(\bt,\blambda,L)$
of $\balpha$-stable parabolic connections with the determinant $(L,\nabla_L)$.
Then we can see by the same proof as Proposition \ref{non-deg-form}
that the restriction $\omega|_{M_C^{\balpha}(\bt,\blambda,L)}$
of $\omega$ to $M_C^{\balpha}(\bt,\blambda,L)$ is also non-degenerate.
So $M_C^{\balpha}(\bt,\blambda,L)$ also has a symplectic structure. \\
(2) By Proposition \ref{splitting}, we obtain a splitting
$\Theta_{M^{\balpha}_{\cC/T}(\tilde{\bt},r,d)/\Lambda^{(n)}_r(d)}
\cong \pi^*\Theta_T \oplus
\Theta_{M^{\balpha}_{\cC/T}(\tilde{\bt},r,d)/T\times\Lambda^{(n)}_r(d)}$.
With respect to this splitting, we can lift $\omega$ to
$\tilde{\omega}\in
H^0(\Omega^2_{M^{\balpha}_{\cC/T}(\tilde{\bt},r,d)/\Lambda^{(n)}_r(d)})$
and we have $d\tilde{\omega}=0$.
This $2$-form $\tilde{\omega}$ is nothing but the $2$-form considered in
\cite{Iw-1} and \cite{Iw-2}.
Note that we have $\tilde{\omega}\cdot v=0$ for $v\in D(\pi^*\Theta_T)$,
where $D$ is the homomorphism given in Proposition \ref{splitting}.
\end{Remark}

\section{Isomonodromic deformation}
\label{isomonodromic deformation}

Let $\tilde{T}\ra T$ be a universal covering.
Then $RP_r(\cC,\tilde{\bt})\times_T\tilde{T}\ra\tilde{T}$
becomes a trivial fibration and we can consider the set of constant sections
\[
 \cF_R^{\sharp}=\left\{
 \sigma:\tilde{T}\ra RP_r(\cC,\tilde{\bt})^{\sharp}\times_T\tilde{T}
 \right\},
\]
where
\[
 RP_r(\cC,\tilde{\bt})^{\sharp}=\coprod_{(x,\ba)\in T\times\cA^{(n)}_r}
 RP_r(\cC_x,\tilde{\bt}_x)_{\ba}^{\sharp}.
\]
The pull back $\tilde{\cF}_M^{\sharp}=\{\RH^{-1}(\sigma)\}$
of this constant sections determines a foliation on
$M^{\balpha}_{\cC/T}(\tilde{\bt},r,d)^{\sharp}\times_T\tilde{T}$
where $M^{\balpha}_{\cC/T}(\tilde{\bt},r,d)^{\sharp}=
\RH^{-1}(RP_r(\cC,\tilde{\bt})^{\sharp})$.
This foliation corresponds to a subbundle determined by a splitting
\[
 D^{\sharp}:
 \tilde{\pi}^*(\Theta_{\tilde{T}}^{an})
 \ra \Theta^{an}_{M^{\balpha}_{\cC/T}
 (\tilde{\bt},r,d)^{\sharp}\times_T\tilde{T}}
\]
of the analytic tangent map
\[
 \Theta^{an}_{M^{\balpha}_{\cC/T}
 (\tilde{\bt},r,d)^{\sharp}\times_T\tilde{T}}\lra
 \tilde{\pi}^*\Theta^{an}_{\tilde{T}}  \ra 0,
\]
where
$\tilde{\pi}:M^{\balpha}_{\cC/T}(\tilde{\bt},r,d)^{\sharp}
\times_T\tilde{T} \ra \tilde{T}$
is the projection.
We will show that this splitting $D^{\sharp}$ is in fact induced by a splitting
\[
 D:\pi^*(\Theta_T)
 \ra \Theta_{M^{\balpha}_{\cC/T}(\tilde{\bt},r,d)}
\]
of the algebraic tangent bundle,
where $\pi:M^{\balpha}_{\cC/T}(\tilde{\bt},r,d)\ra T$ is the projection.
This splitting is nothing but the differential equation determined by
the isomonodromic deformation.

Usually, the equation of the isomonodromic deformation
is given by the following way:
Assume $C=\BP^1$ and take a Zariski open set
$U\subset M_{\BP^1\times T/T}^{\balpha}(\tilde{\bt},r,0)$
such that a universal family on $\BP^1\times U$ is given by
$(\cO^{\oplus r}_{\BP^1\times U},\nabla,\{\tilde{l}^{(i)}_j\})$,
where $\nabla$ is given by a connection matrix
\[
 \sum_{i=1}^n \frac{A_idz}{z-t_i}
 \quad (A_i=(a^{(i)}_{jk})_{1\leq j,k\leq r}\in M_r(\cO_U)).
\]
Then the differential equation of the isomonodromic deformation
is given by
\begin{equation}\label{schlesinger}
 \frac{\partial A_i}{\partial t_j}=-\frac{[A_i,A_j]}{t_i-t_j}
 \quad (t_i\neq t_j), \quad
 \frac{\partial A_i}{\partial t_i}=
 \sum_{k\neq i}\frac{[A_i,A_k]}{t_i-t_k}.
\end{equation}
The equation (\ref{schlesinger}) is called the Schlesinger equation
and it is needless to say that this equation plays a great role
in explicit description of many differential equations
arising from the isomonodromic deformation.
However, $U$ and the equation (\ref{schlesinger})
do not have the geometric Painlev\'e property.
For example, take a point $x$ of $U$ over $\bt\in T$
which corresponds to $(\cO^{\oplus r}_{\BP^1_U},\nabla,\{l^{(i)}_j\})$.
Take a path $\gamma$ in $M^{\balpha}_{\cC/T}(\tilde{\bt},r,d)$
starting at $x$ and comming back to a point $y$ over $\bt$
such that $\gamma$ is a solution of the differential equation
determined by the isomonodromic deformation.
Let $(F,\nabla_F,\{(l^{(i)}_F)_j\})$ be the parabolic connection
corresponding to $y$.
Then $F$ may not be trivial and so $y$ may not lie in $U$.
Thus $U$ and (\ref{schlesinger}) do not have the geometric
Painlev\'e property.
Moreover, we want to describe the isomonodromic deformation 
even on the locus contracted by $\RH$.
So the equation (\ref{schlesinger}) itself is not enough for the
explicit geometric description of the isomonodromic deformation,
because (\ref{schlesinger}) is defined only on a Zariski open set
and is even not defined in higher genus case.
In order to obtain the splitting $D$, we will go back to how
the equation (\ref{schlesinger}) arises as the isomonodromic
deformation.
Take a local constant section
$\sigma:T'\to RP_r(\cC,\tilde{\bt})^{\sharp}\times_TT'$,
where $T'$ is an analytic open subset of $T$.
Then $\sigma$ corresponds to a local system
$\V_{\sigma}$ on
$\cC_{T'}\setminus((\tilde{t_1})_{T'}+\cdots+(\tilde{t_n})_{T'})$
after shrinking $T'$.
The canonical flat connection on
$\V_{\sigma}\otimes
\cO_{\cC_{T'}\setminus((\tilde{t_1})_{T'}+\cdots+(\tilde{t_n})_{T'})}$
extends to an integrable logarithmic connection
\begin{equation}\label{flat-connection}
 \nabla^{\sigma}:\tilde{E}_{\RH^{-1}(\sigma(T'))}\longrightarrow
 \tilde{E}_{\RH^{-1}(\sigma(T'))}\otimes\Omega^1_{\cC_{T'}}
 (\log((\tilde{t_1})_{T'}+\cdots+(\tilde{t_n})_{T'}))
\end{equation}
whose induced relative connection is just
$\tilde{\nabla}|_{\cC\times_T\RH^{-1}(\sigma(T'))}$.
The equation (\ref{schlesinger}) comes from the integrability
condition for $\nabla^{\sigma}$.
Taking the tangent direction of $\sigma(T')$, we obtain
the splitting $D$.
Precisely, we have the following proposition:

\begin{Proposition}\label{splitting}
 There exists an algebraic splitting
 \[
  D:\pi^*(\Theta_T)
 \ra \Theta_{M^{\balpha}_{\cC/T}(\tilde{\bt},r,d)}
 \]
 of the tangent map
 $\Theta_{M^{\balpha}_{\cC/T}(\tilde{\bt},r,d)}
 \rightarrow \pi^*(\Theta_T)$
 whose pull-back to
 $M_{\cC/T}^{\balpha}(\tilde{\bt},r,d)^{\sharp}\times_T\tilde{T}$
 coincides with $D^{\sharp}$.
\end{Proposition}

\begin{proof}
We will define the corresponding homomorphism
$\Theta_T\rightarrow\pi_*\Theta_{M^{\balpha}_{\cC/T}(\tilde{\bt},r,d)}$.
Take any affine open set $U\subset T$
and a vector field $v\in H^0(U,\Theta_T)$.
Then $v$ corresponds to a morphism
$\iota^v:\Spec\cO_U[\epsilon] \rightarrow T$
with $\epsilon^2=0$ such that the composite
$U\hookrightarrow\Spec\cO_U[\epsilon]\rightarrow T$
is just the inclusion $U\hookrightarrow T$.
We denote the restriction of the universal family to
$\cC\times_T \pi^{-1}(U)$ simply by
$(\tilde{E},\tilde{\nabla},\{\tilde{l}^{(i)}_j\})$.
Consider the fiber product
$\cC\times_T\Spec\cO_{\pi^{-1}(U)}[\epsilon]$
with respect to the canonical projection $\cC\rightarrow T$
and the composite
$\Spec\cO_{\pi^{-1}(U)}[\epsilon]\rightarrow\Spec\cO_U[\epsilon]
\stackrel{\iota^v}\rightarrow T$.
We denote the pull-back of $D(\tilde{\bt})$ by the morphism
$\cC\times_T\Spec\cO_{\pi^{-1}(U)}[\epsilon]\rightarrow \cC$ simply by
$D(\tilde{\bt})_{\cO_{\pi^{-1}(U)}[\epsilon]}$.
We call
$({\mathcal E},\nabla^{\mathcal E},\{(l_{\mathcal E})^{(i)}_j\})$
a horizontal lift of
$(\tilde{E},\tilde{\nabla},\{\tilde{l}^{(i)}_j\})$
if
\begin{enumerate}
\item $\mathcal E$ is a vector bundle on
 $\cC\times_T\Spec\cO_{\pi^{-1}(U)}[\epsilon]$,
\item
 ${\mathcal E}|_{\tilde{t}_i\times\Spec\cO_{\pi^{-1}(U)}[\epsilon]}
 =(l_{\mathcal E})^{(i)}_0\supset\cdots\supset (l_{\mathcal E})^{(i)}_r=0$
 is a filtration by subbundles for $i=1,\ldots,n$,
\item $\nabla^{\mathcal E}:{\mathcal E}\rightarrow
 {\mathcal E}\otimes\Omega^1_{\cC\times_T\Spec\cO_{\pi^{-1}(U)}[\epsilon]/
 \pi^{-1}(U)}
 \left(\log\left(D(\tilde{\bt})_{\cO_{\pi^{-1}(U)}[\epsilon]}\right)\right)$
 is a connection satisfying
 \begin{enumerate}
  \item $\nabla^{\mathcal E}(F^{(i)}_j({\mathcal E}))\subset
  F^{(i)}_j({\mathcal E})\otimes
  \Omega^1_{\cC\times_T\Spec\cO_{\pi^{-1}(U)}[\epsilon]/
  \pi^{-1}(U)}
  \left(\log D(\tilde{\bt})_{\cO_{\pi^{-1}(U)}[\epsilon]}\right)$,
  where $F^{(i)}_j({\mathcal E})$ is given by
  $F^{(i)}_j({\mathcal E}):=
  \ker\left({\mathcal E}\rightarrow
  {\mathcal E}|_{\tilde{t}_i\times_T\Spec\cO_{\pi^{-1}(U)}[\epsilon]}/
  (l^{\mathcal E})^{(i)}_j\right)$,
  \item The curvature
  $\nabla^{\mathcal E}\circ\nabla^{\mathcal E}:{\mathcal E}\rightarrow
  {\mathcal E}\otimes\Omega^2_{\cC\times_T\Spec\cO_{\pi^{-1}(U)}[\epsilon]/
  \pi^{-1}(U)}
  \left(\log\left(D(\tilde{\bt})_{\cO_{\pi^{-1}(U)}[\epsilon]}
  \right)\right)$
  is zero,
  \item $(\res_{\tilde{t}_i\times_T\Spec\cO_{\pi^{-1}(U)}[\epsilon]}
  (\tilde{\nabla}^{\mathcal E})-\tilde{\lambda}^{(i)}_j)
  ((l^{\mathcal E})^{(i)}_j)\subset (l^{\mathcal E})^{(i)}_{j+1}$
  for any $i,j$, where
  $\tilde{\nabla}^{\mathcal E}$ is the relative connection over
  $\Spec\cO_{\pi^{-1}(U)}[\epsilon]$
  induced by $\nabla^{\mathcal E}$ and
  \item $({\mathcal E},\tilde{\nabla}^{\mathcal E},
  \{(l^{\mathcal E})^{(i)}_j)\otimes \cO_{\pi^{-1}(U)}[\epsilon]/(\epsilon)
  \cong (\tilde{E},\tilde{\nabla},\{\tilde{l}^{(i)}_j\})$.
 \end{enumerate}
\end{enumerate}
Here we define the sheaf
$\Omega^1_{\cC\times_T\Spec\cO_{\pi^{-1}(U)}[\epsilon]/\pi^{-1}(U)}
\left(\log\left(D(\tilde{\bt})_{\cO_{\pi^{-1}(U)}[\epsilon]}\right)\right)$
as the coherent subsheaf of
$\Omega^1_{\cC\times_T\Spec\cO_{\pi^{-1}(U)}[\epsilon]/\pi^{-1}(U)}
\left(D(\tilde{\bt})_{\cO_{\pi^{-1}(U)}[\epsilon]}\right)$
locally generated by $\tilde{g}^{-1}d\tilde{g}$ and $d\epsilon$
for a local defining equation $\tilde{g}$ of
$D(\tilde{\bt})_{\cO_{\pi^{-1}(U)}[\epsilon]}$
and the sheaf
$\Omega^2_{\cC\times_T\Spec\cO_{\pi^{-1}(U)}[\epsilon]/\pi^{-1}(U)}
\left(\log\left(D(\tilde{\bt})_{\cO_{\pi^{-1}(U)}[\epsilon]}\right)\right)$
as the coherent subsheaf of
$\Omega^2_{\cC\times_T\Spec\cO_{\pi^{-1}(U)}[\epsilon]/\pi^{-1}(U)}
\left(D(\tilde{\bt})_{\cO_{\pi^{-1}(U)}[\epsilon]}\right)$
locally generated by $\tilde{g}^{-1}d\tilde{g}\wedge d\epsilon$.
Assume that the parabolic connection
$(\tilde{E},\tilde{\nabla},\{\tilde{l}^{(i)}_j\})$
is locally given in a small affine open subset $W$ of
$\cC\times_T\pi^{-1}(U)$
by a connection matrix
$Ag^{-1}dg$, where $g$ is a local defining equation of
$(\tilde{t}_i\times\pi^{-1}(U))\cap W$,
$A\in M_r(\cO_W)$,
$A((\tilde{t}_i\times\pi^{-1}(U))\cap W)$
is an upper triangular matrix
and the parabolic structure $\{\tilde{l}^{(i)}_j\}_W$ is given by 
$(\tilde{l}^{(i)}_j)_W=(\overbrace{*,*,\ldots,*}^{r-j},0,\ldots,0)$.
Let $\tilde{W}$ be the affine open subscheme of
$\cC\times\Spec\cO_{\pi^{-1}(U)}[\epsilon]$ whose underlying space
is $W$ and
let $\tilde{g}$ be a lift of $g$ which is a local defining equation
of $(\tilde{t}_i\times\Spec\cO_{\pi^{-1}(U)}[\epsilon])\cap \tilde{W}$.
If we denote the composite
$\cO_{\tilde{W}}\stackrel{d}\rightarrow\Omega^1_{\tilde{W}/\pi^{-1}(U)}
=\cO_{\tilde{W}}d\tilde{g}\oplus\cO_{\tilde{W}}d\epsilon
\rightarrow\cO_{\tilde{W}}d\epsilon$
by $d_{\epsilon}$,
we can take a lift $\tilde{A}\in M_r(\cO_{\tilde{W}})$ of $A$ such that
$d_{\epsilon}(\tilde{A})=0$ and
$\tilde{A}((\tilde{t}_i\times\Spec\cO_{\pi^{-1}(U)}[\epsilon])\cap \tilde{W})$
is an upper-triangular matrix.
Then the connection matrix $\tilde{A}\tilde{g}^{-1}d\tilde{g}$
gives a local horizontal lift of
$(\tilde{E},\tilde{\nabla},\{\tilde{l}^{(i)}_j\})|_W$.
We put
\begin{align*}
 \cF^0&:=\left\{ u\in {\mathcal End}(\tilde{E}) \left|
 u|_{\tilde{t}_i\times\pi^{-1}(U)}(\tilde{l}^{(i)}_j)\subset
 \tilde{l}^{(i)}_j
 \quad (1\leq i\leq n, 0\leq j\leq r) \right\}\right. \\
 \cF^1&:=\left\{ u\in {\mathcal End}(\tilde{E})\otimes\tilde{\Omega}^1
 \left|
 \begin{array}{l}
 \text{$u(F^{(i)}_j(\tilde{E}))\subset
 F^{(i)}_j(\tilde{E})\otimes\tilde{\Omega}^1$
 for any $i,j$ and the image of} \\
 \text{$F^{(i)}_j(\tilde{E})\stackrel{u}\rightarrow
 \tilde{E}\otimes\tilde{\Omega}^1
 \stackrel{\res_{\tilde{t}_i}}\rightarrow
 \tilde{E}|_{\tilde{t}_i\times\pi^{-1}(U)}$
 is contained in $\tilde{l}^{(i)}_{j+1}$ for any $i,j$}
 \end{array}
 \right\}\right. \\
 \cF^2&:=\left\{ u\in {\mathcal End}(\tilde{E})\otimes\tilde{\Omega}^2
 \left|
 \text{$u(F^{(i)}_j(\tilde{E}))\subset
 F^{(i)}_{j+1}(\tilde{E})\otimes\tilde{\Omega}^2$
 for any $i,j$}
 \right\}\right.,
\end{align*}
where we put
\begin{align*}
 \tilde{\Omega}^1 &:=
 \Omega^1_{\cC\times_T\pi^{-1}(U)/\pi^{-1}(U)}(D(\tilde{\bt}))
 \oplus \cO_{\cC\times_T\pi^{-1}(U)}d\epsilon \\
 \tilde{\Omega}^2 &:=
 \Omega^1_{\cC\times_T\pi^{-1}(U)/\pi^{-1}(U)}(D(\tilde{\bt}))
 \wedge d\epsilon \\
 F^{(i)}_j(\tilde{E}) &:=
 \ker\left( \tilde{E}\rightarrow
 \tilde{E}|_{\tilde{t}_i\times\pi^{-1}(U)}/\tilde{l}^{(i)}_j \right).
\end{align*}
Consider the complex
\[
 \cF^0 \stackrel{d^0}\longrightarrow \cF^1
 \stackrel{d^1}\longrightarrow \cF^2
\]
defined by
$d^0(u)=\tilde{\nabla}\circ u-u\circ\tilde{\nabla}+ud\epsilon$
and
$d^1(\omega+ad\epsilon)=d\epsilon\wedge\omega
+(\tilde{\nabla}\circ a-a\circ\tilde{\nabla})\wedge d\epsilon$
for $u\in\cF^0$, $\omega\in{\mathcal End}(\tilde{E})\otimes g^{-1}dg$
and $a\in\cF^0$ such that $\omega+ad\epsilon\in\cF^1$.

We can see that an obstruction class for the existence of
a horizontal lift of
$(\tilde{E},\tilde{\nabla},\{\tilde{l}_i\})$
is in ${\bf H}^2({\mathcal F}^{\bullet})$ and
the set of horizontal lifts is isomorphic to
${\bf H}^1({\mathcal F}^{\bullet})$ if it is not empty.
It is obvious that $d^0$ is injective and
$d^1$ is surjective.
Take any member
$\omega+ad\epsilon\in\ker(d^1)$.
Then we have
$d^1(\omega+ad\epsilon)=d\epsilon\wedge\omega
+(\tilde{\nabla}\circ a-a\circ\tilde{\nabla})\wedge d\epsilon
=0$.
So we have
$\omega=\tilde{\nabla}\circ a-a\circ\tilde{\nabla}$ and
$d^0(a)=\tilde{\nabla}\circ a-a\circ\tilde{\nabla}+ad\epsilon
=\omega+ad\epsilon$.
Thus ${\mathcal F}^{\bullet}$ is exact at ${\mathcal F}^1$.
Hence ${\mathcal F}^{\bullet}$
is an exact complex.
Thus we have ${\bf H}^2({\mathcal F}^{\bullet})=0$,
${\bf H}^1({\mathcal F}^{\bullet})=0$ and
there is a unique horizontal lift
$({\mathcal E},\nabla^{\mathcal E},
\{(l^{\mathcal E})^{(i)}_j\})$
of $(\tilde{E},\tilde{\nabla},\{\tilde{l}_i\})$.
We can check that
$v\mapsto({\mathcal E},\tilde{\nabla}^{\mathcal E},
\{(l^{\mathcal E})^{(i)}_j\})$
gives an ${\mathcal O}_T$-homomorphism
$\Theta_T\rightarrow\pi_*(\Theta_{M^{\balpha}_{\cC/T}(\tilde{\bt},r,d)})$.
Thus we obtain an
${\mathcal O}_{M^{\balpha}_{\cC/T}(\tilde{\bt},r,d)}$-homomorphism
\[
  D:\pi^*(\Theta_T)
 \ra \Theta_{M^{\balpha}_{\cC/T}(\tilde{\bt},r,d)}
\]
which is a splitting of the surjection
$\Theta_{M^{\balpha}_{\cC/T}(\tilde{\bt},r,d)}
 \rightarrow \pi^*(\Theta_T)$.
By construction we can see that the pull-back of $D$ to
$M^{\balpha}_{\cC/T}(\tilde{\bt},r,d)^{\sharp}\times\tilde{T}$
coincides with $D^{\sharp}$.
\end{proof}

\begin{Remark}\rm
(1) We can see by construction that the image $D(\pi^*\Theta_T)$
is in fact contained in the relative tangent bundle
$\Theta_{M^{\balpha}_{\cC/T}(\tilde{\bt},r,d)/\Lambda^{(n)}_r(d)}$. \\
(2) As is explained in section 2,
the subbundle
$D(\pi^*\Theta_T)\subset\Theta_{M^{\balpha}_{\cC/T}(\tilde{\bt},r,d)}$
satisfies an integrability condition and determines a foliation
$\cF_M$ on $M^{\balpha}_{\cC/T}(\tilde{\bt},r,d)$.
By the proof of Proposition \ref{splitting}, we can see that
$\{\RH({\mathcal L})|{\mathcal L}\in{\mathcal F}_M\}$
coincides with ${\mathcal F}_R$.
\end{Remark}

\noindent
{\bf Proof of Theorem \ref{gpp-thm}.}
Take any path
$\gamma:[0,1]\rightarrow T$
and a point $x\in M^{\balpha}_{\cC/T}(\tilde{\bt},r,d)$
such that $\pi(x)=\gamma(0)$.
We take a lift
$\tilde{x}\in M^{\balpha}_{\cC/T}(\tilde{\bt},r,d)\times_T\tilde{T}$
of $x$.
We can easily lift $\gamma$ to a path
$\tilde{\gamma}:[0,1]\rightarrow
(RP_r(\cC,\tilde{\bt})\times_T\tilde{T})
\times_{\cA^{(n)}_r}\Lambda^{(n)}_r(d)$
such that $\tilde{\gamma}[0,1]$ is contained in a leaf of
${\mathcal F}_R\times_{\cA^{(n)}_r}\Lambda^{(n)}_r(d)$
and $\RH(\tilde{x})=\tilde{\gamma}(0)$.
Since a leaf of $\tilde{\mathcal F}_M$ passing through $\tilde{x}$
is locally isomorphic to $\tilde{T}$, we can take a unique local lift
$\delta_{\epsilon}:[0,\epsilon)\rightarrow 
M^{\balpha}_{\cC/T}(\tilde{\bt},r,d)\times_T\tilde{T}$
of $\tilde{\gamma}$ for small $\epsilon>0$ such that
$\delta_{\epsilon}(0)=\tilde{x}$
and $\delta_{\epsilon}([0,\epsilon])$ is contained in a leaf of
$\tilde{\cF}_M$.
Consider the set
\[
 A:=\left\{t\in(0,1]\left|
 \begin{array}{l}
 \text{there exists a lift $\delta_t:[0,t)\rightarrow
 M^{\balpha}_{\cC/T}(\tilde{\bt},r,d)\times_T\tilde{T}$
 of $\tilde{\gamma}|_{[0,t)}$} \\
 \text{such that $\delta(0)=\tilde{x}$ and $\delta([0,t))$
 is contained in a leaf of $\tilde{{\mathcal F}}_M$}
 \end{array}
 \right\}\right.
\]
and put
$a:=\sup A$.
Assume that $a<1$.
Then there is a sequence $\{t_n\}_{n\geq 0}$ in $A$ such that
$\lim_{n\to\infty}t_n=a$.
Since
$\RH:M^{\balpha}_{\cC/T}(\tilde{\bt},r,d)\times_T\tilde{T}\rightarrow
(RP_r(\cC,\tilde{\bt})\times_T\tilde{T})\times_{\cA^{(n)}_r}\Lambda^{(n)}_r(d)$
is a proper morphism, $\RH^{-1}(\tilde{\gamma}[0,a])$ is compact.
Thus the sequence $\{\delta_{t_n}(t_n)\}$ in $\RH^{-1}(\tilde{\gamma}[0,a])$
has a subsequence $\{\delta_{t_{n_k}}(t_{n_k})\}$ which is convergent
in $\RH^{-1}(\tilde{\gamma}[0,a])$.
We put
$\tilde{x}_a:=
\lim_{k\to\infty}\delta_{t_{n_k}}(t_{n_k})\in\RH^{-1}(\tilde{\gamma}[0,a])$.
Then $\RH(\tilde{x}_a)=\tilde{\gamma}(a)$ and we can take a local lift
$\tilde{\delta}_a:(a-\epsilon,a+\epsilon)\rightarrow
M^{\balpha}_{\cC/T}(\tilde{\bt},r,d)\times_T\tilde{T}$
of $\tilde{\gamma}|_{(a-\epsilon,a+\epsilon)}$
for small $\epsilon>0$ such that
$\tilde{\delta}_a(a)=\tilde{x}_a$ and
$\tilde{\delta}_a((a-\epsilon,a+\epsilon))$ is contained in a
leaf of $\tilde{\mathcal F}_M$.
Gluing $\tilde{\delta}_a$ and $\delta_{t_n}$ for some $t_n\in A$
with $t_n>a-\epsilon$, we can obtain a lift
$\delta_a:[0,a+\epsilon)\rightarrow
M^{\balpha}_{\cC/T}(\tilde{\bt},r,d)\times_T\tilde{T}$
of $\tilde{\gamma}|_{[0,a+\epsilon)}$ whose image
is contained in a leaf of $\tilde{\mathcal F}_M$.
Thus $a+\epsilon\in A$ which contradicts the choice of $a$.
So we have $a=1$ and we can obtain a lift
$\delta_1:[0,1]\rightarrow
M^{\balpha}_{\cC/T}(\tilde{\bt},r,d)\times_T\tilde{T}$
of $\tilde{\gamma}$ whose image is contained in a leaf of
$\tilde{\mathcal F}_M$.
Taking the image by the projection
$M^{\balpha}_{\cC/T}(\tilde{\bt},r,d)\times_T\tilde{T}
\rightarrow M^{\balpha}_{\cC/T}(\tilde{\bt},r,d)$,
we obtain an $\cF_M$-horizontal lift $\delta$ of $\gamma$
such that $\delta(0)=x$.
\hfill $\square$

\end{document}